\documentclass[11pt]{amsart}
\usepackage{amssymb}
\usepackage{amsmath,amscd}
\usepackage{amsfonts,amssymb}
\usepackage{latexsym}
\usepackage{amsbsy}

\theoremstyle{definition}

\renewcommand{\Box}{\hfill\qed}
\newcommand{\bcen}{\begin{center}}      \newcommand{\ecen}{\end{center}}
\newcommand{\bay}{\begin{array}}      \newcommand{\eay}{\end{array}}
\newcommand{\beq}{\begin{eqnarray*}}      \newcommand{\eeq}{\end{eqnarray*}}
\def\lr#1{\langle #1\rangle}
\newcommand{\cal}{\mathcal}
\def\az{\alpha}    
\def\bz{\beta}  \def\dt{\Delta}   \def\bd{{\bf d}}
\def\gz{\gamma}
  \def\ooz{\Omega} 
\def\sz{\sigma}
\def\oz{\omega}

\def\dz{\delta}  

\def\llz{\Lambda}    \def\vhi{\varphi}  
\def\lz{\lambda}    \def\sset{\subseteq}
\def\cl{{\cal L}}     \def\cp{{\cal P}} \def\cc{{\cal C}}
\def\ch{{\cal H}}   \def\ct{{\cal T}} \def\co{{\cal O}}
  \def\cd{{\cal D}} \def\cl{{\cal L}}
\def\ci{{\cal I}}   
\def\ca{{\cal A}}
\def\cb{{\cal B}}
\def\cz{{\cal Z}}
\def\ck{{\cal K}}
\def\ct{{\cal T}}
\def\ce{{\cal E}}
\def\cm{{\cal M}}
\def\cn{{\cal N}}
\def\bbn{{\mathbb N}}  \def\bbz{{\mathbb Z}}  \def\bbq{{\mathbb Q}}
\def\bbf{{\mathbb F}}   \def\bbe{{\mathbb E}} \def\bbf{{\mathbb F}}
\def\fq{{\mathbb F}_q}
  \def\bbc{{\mathbb C}}
   \def\fn{\mathfrak n}

\def\te{{\tilde E}}
\def\bi{{\bf 1}}
   
\def\fk#1{{\mathfrak #1}}   
 
\def\pe{\preceq}
\def\lra{\longrightarrow}  \def\strl{\stackrel}

\def\ra{\rightarrow}
\def\hom{\operatorname{Hom}}

\def\ext{\operatorname{Ext}}
\def\aut{\operatorname{Aut}}

\def\dim{\operatorname{dim}}
\def\cdim{\operatorname{codim}}
\def\min{\operatorname{min}}
\def\udim{\operatorname{\underline {dim}}}
\def\ed{\operatorname{End}}
\def\ol{\overline}
\def\mod{\operatorname{mod}}  
    \def\aut{\operatorname{Aut}}

\def\uq2{U_q(\hat{sl}_2)}
\def\fq{\bbf_q}

\def\wc{{\bf c}}
\def\wm{{\bf m}}
\def\bl{{\bf L}}
\def\bj{{\bf J}}
\def\bfz{{\bf Z}}
\def\bb{{\bf b}}
\def\ba{{\bf a}}
\def\be{{\bf e}}
\def\bc{{\bf c}}
\def\bt{{\bf t}}
\def\bd{{\bf d}}
\def\nd{{\noindent}}
\def\mk{{\medskip}}
\def\hs{{\hskip 1cm}}
\def\mh{{\mathcal{H}}}
\def\mr{{\mathcal{R}}}
\def\mf{{\mathcal{F}}}

\title[Affine canonical bases]{Representations of tame quivers and affine canonical bases}
\thanks{The main results
in this article were reported by J. Xiao at the Fields Institute,
Toronto, Canada (Workshop and Conference on Infinite Dimensional
Lie Theory and its Application, July 17-25, 2003), at the Algebra
seminar of the University of Sydney, Australia, October 10, 2003
and at Snowbird, Utah, USA (AMS-IMS-SIAM Summer Research
Conferences, Representations of Algebraic Groups, Quantum Groups,
and Lie Algebras, July 11-15, 2004).}
\thanks{  The research was
supported in part by NSF of China and by the 973 Project of the
Ministry of Science and Technology of China and by NSF grant
DMS-0200673. }
\author{Zongzhu Lin}
\address{Department of Mathematics\\
Kansas State University\\
Manhattan, KS 66506, U.S.A.}
\email{zlin@math.ksu.edu}
\author{Jie
Xiao}
\address{Department of Mathematical Sciences\\
Tsinghua University\\
Beijing 10084, P.~R.~China}
\email{jxiao@math.tsinghua.edu.cn}
\author{Guanglian Zhang}
\address{Department of Mathematical Sciences\\
Tsinghua University\\
Beijing 10084, P.~R.~China}
\email{zhangguanglian@mails.tsinghua.edu.cn}
\dedicatory{ Dedicated to Claus Michael Ringel on the occasion
of his 60th birthday}

\date{}

\begin{document}

\begin{abstract} An integral PBW-basis of type $A_1^{(1)}$ has been
constructed by Zhang [Z] and  Chen [C] using  the Auslander-Reiten
quiver of the Kronecker quiver. We associate  a geometric order to
elements in this basis  following an idea of Lusztig [L1] in the
case of finite type. This leads to an algebraic realization of a
bar-invariant  basis of $\uq2$.  For any affine symmetric type, we
obtain an integral PBW-basis of the generic composition algebra,
by using an algebraic construction of the integral basis for a
tube in [DDX], an embedding of the module category of the
Kronecker quiver into the module category of the tame quiver, and
a list of the root vectors of indecomposable modules according to
the preprojective, regular, and preinjective components of the
Auslander-Reiten quiver of the tame quiver. When the basis
elements are ordered to be compatible with the geometric order
given by the dimensions of the orbit varieties and the extension
varieties, we can show that the transition matrix between the
PBW-basis and a monomial basis is triangular with diagonal entries
equal to $1$. Therefore we  obtain a bar-invariant basis. By a
orthogonalization for the PBW-basis with the inner product, we
finally give an algebraic way to realize the canonical bases of
the quantized enveloping algebras of all symmetric affine
Kac-Moody Lie algebras.
\end{abstract}

\maketitle

\bigskip
\bigskip

\centerline{\bf 0. Introduction}

\bigskip
\nd{\bf 0.1} Let $U^+$ be the positive part of the quantized
enveloping algebra of $U$ associated to a Cartan datum.
For a finite type root system,
Lusztig's construction of the canonical basis of $U^+$ [L1]
involves three ingredients. The first one can be understood
 as purely combinatorial. By applying  Lusztig's symmetries
and the induced actions of the braid group on $U^+,$ one may have
a complete list of root vectors of $U^+.$ Associated to each
reduced expression of the longest element of the Weyl group, there
is a PBW-basis of $U^+$ with a specific order and  a monomial
basis on the Chevalley generators such that the transition matrix
between these two bases is triangular with diagonal entries equal
to $1.$ (See [L1, 7.8-7.9].) The second  is the quiver approach. Each isomorphism
class of the Dynkin quiver corresponds to a PBW-type basis element $E^c$,
($c\in\bbn^{\Phi^+}$) of $U^+.$ Now the representations of a fixed
dimension vector of the quiver are the orbits of an algebraic
group action on an affine variety. The geometric dimension of
these orbits can be applied to give an order in
$\{E^c|c\in\bbn^{\Phi^+}\}.$ This ordered basis relates to a
monomial basis by a triangular transition matrix with diagonal
entries equal to $1.$ By a standard linear algebra method one can
easily obtain the canonical basis.  The third  is the
geometric approach by using perverse sheaves and intersection
cohomology.  There is also a different approach to
construct the global crystal basis of $U^+$  in the Kashiwara's
work [K]. Now it is well known that Lusztig has generalized his
geometric method to construct the canonical bases of $U^+$ for all
infinite type (see[L2] and [L3]).

\mk\nd{\bf 0.2}  Although  most knowledge on the canonical basis
in finite type can be carried out in a pure combinatorial way, it
is obvious to see that  the definition of the canonical basis was
introduced by Lusztig in a framework of  representations of
quivers. Specifically,  Lusztig has extended the Gabriel's theorem
to build up a PBW type  basis for  $U^+$, which is ordered by the
geometric properties of the corresponding orbit varieties. The
representation category of a tame quiver has been completely
described by a generalization of the Gabriel's theorem and its
Auslander-Reiten quiver (see[DR]). The objective of  this paper is
to provide a process to construct a PBW type basis and
characterize the canonical basis of $U^+$ of affine type by using
Ringel-Hall algebra and the knowledge of the representations of
tame quivers. We hope that the approach we adopt here  is closer
to Lusztig's original idea of [L1].

 \mk\nd{\bf 0.3} For infinite type root systems,
there is no longest elements and the braid group action does not
construct PBW-type basis. A natural question is to seek an algebraic
construction of PBW type basis, monomial basis  and the canonical
basis, just like Luszitg did for finite type cases. For affine
types, a PBW type basis was attained first by Beck, Chari and
Pressley in [BCP] for the quantized enveloping algebra of untwisted
affine type, and then was improved and extended by Beck and Nakajima
in [BN] to all twisted and untwisted affine types. Their approach is
to give the real root vectors by applying Lusztig's symmetries on
the generators and to construct the imaginary root vectors by using
Schur functions in the Heisenberg generators; and then use these
PBW-bases with the almost orthonormal property to obtain the crystal
bases. However we like to point out that the order of the PBW-basis
elements from the representations of tame quivers is different from
theirs. A detailed analysis for this order enables us to  construct
the PBW-basis, also the monomial basis and a triangular transition
matrix with diagonal entries equal to $1.$ Then we can use the
standard linear method, which was used by Lusztig for finite type
cases, to obtain the canonical basis.

\mk\nd{\bf 0.4} In Section~1 we recall the definition of Hall
algebras of quivers by Ringel and by Lusztig respectively, and
point out that the two constructions coincide essentially for the
representations of a quiver over a finite field. Section~2
presents the basic geometric properties of the orbit varieties and
extension varieties for the representations of quivers. In
Section~3 we construct an integral PBW basis of $A^{(1)}_1$ type
by using the representations of the Kronecker quiver. Most results
in this section are already known for some experts (see [Z] and
[C]). The category $\mod\llz$  of the Kronecker quiver has a
strong representation-directed property [DR]. This enables us in
Section~4 to arrange the positive roots in a special order. In
addition, by the basic properties of the orbit varieties, we find
a monomial basis whose  transition matrix with the PBW basis is
triangular with diagonal entries equal to $1.$ Section~5 is taken
from [DDX], in which the integral basis and the canonical basis of
$A^{(1)}_n$ type were given in terms of the nilpotent
representations of the cyclic quivers. In Section~6 we consider
the $\cz$-submodule of $U^+$ generated by $\lr{u_M}$ for $M$ being
preprojective or preinjective. It is a $\cz$-subalgebra of $U^+.$
An integral basis for this $\cz$-subalgebra can be listed in an
order with respect to the representation-directed property of the
preprojective (resp. preinjective ) component. We verify that the
basis elements are products of images of Chevalley generators
under the action of the sequences of  Lusztig's symmetries in an
admissible order. So the situation in Section~6 resembles the
construct of PBW-type basis in the finite type case. In Section~7,
we show that the subalgebras corresponding to the preprojective
component, preinjective component, non-homogeneous tubes, and an
embedding of the module category of the Kronecker quiver can be
put together, according to the representation-directed property of
the tame quiver. This gives rise to an integral basis of $U^+$
over $\bbq[v,v^{-1}].$ In  Section 8, we again find a monomial
basis, which has a unipotent triangular relation with the integral
PBW type basis we obtained. But this needs a little more subtle
analysis of the orbit varieties and the extension varieties.
Finally, a bar-invariant basis $\{\ce^{\bc}|\bc\in\cm\}$ of $U^+$
can be constructed in an elementary and algebraic way. The last
section is new to an old version of the paper. By a detailed
calculation of the inner product on the PBW-basis in  the
orthogonalization  process using the properties of the Schur
functions,  we can answer Nakajima's question in [N]
affirmitively, that is, we show that the basis
$\{\ce^{'\bc}|\bc\in\cm\},$ which is a modified form of the basis
$ \{\ce^{\bc}|\bc\in\cm\},$ exactly equal to the canonical basis
in [L2].

In a preprint [Li], Y.Li has given a geometric characterization of
the monomial basis $\{{m_{\bc}|\bc\in\cm}\},$ which is constructed
 by us in Section 8, and he has proved that the transition matrix
 between $\{{m_{\bc}|\bc\in\cm}\}$ and the canonical basis is
 triangular with diagonal entries equal to $1$ and entries above the
 diagonal in $\bbz_{\geq 0}[v,v^{-1}].$
Finally we like to remark that our construction of the
canonical bases is independent of the assumption for the existence
of  Lusztig's canonical basis, or the existence
of Kashiwara's global crystal basis.

\mk\nd{\bf Acknowledgments.}  (1) We are very grateful to
O.~Schiffmann, B.~Deng and J.~Du for getting our attention to the
preprint [H] by A.~Hubery, in which an integral PBW basis for the
composition algebras of affine type are constructed according to
the representations of tame quivers. (2) We thank H.Nakajima very
much for his encouragement and suggestions for us to solve his
question in [N], in particular, for his suggestion to use Schur
functions to modify $E_{n\delta}$'s. (3) We thank F.Xu for his
great help in writing the present version of this paper.

\bigskip
\centerline{\bf 1. Ringel-Hall algebras }
\bigskip
\nd {\bf 1.1} A quiver $Q=(I,H,s,t)$ consists of  a
vertex set $I$,  an arrow set $H$, and two maps $s,t:H\rightarrow I$
such that an arrow $\rho\in H$ starts at $s(\rho)$ and
terminates at $t(\rho).$

Throughout the paper, $\bbf_q$ denotes a finite field with $q$
elements,  $Q=(I,H,s,t)$ is a fixed connected
quiver, and $\llz=\bbf_q Q$ is the path algebra of $Q$ over $\bbf_q.$ By
$\mod\llz$ we denote the category of all finite dimensional left
$\llz$-modules, or equivalently finite modules. It is well-known
that $\mod\llz$ is equivalent to the category of finite
dimensional representations of $Q$ over $\bbf_q.$ We shall simply
identify $\llz$-modules with representations of $Q.$

\mk \nd{\bf 1.2 Ringel-Hall algebra. }Given three modules $L,M,N$ in
$\mod\llz,$ let $g^L_{M N}$ denote the number of $\llz$-submodules $W$ of
$L$ such that $W\simeq N$ and $L/W\simeq M$ in $\mod \llz$. More generally, for
$M_1,\cdots,M_t, L\in \mod\llz,$ let $g^{L}_{M_1\cdots M_t}$ denote
the number of the filtrations $0=L_0\subseteq L_1\subseteq\cdots
\subseteq L_t=L$ of $ \llz$-submodules such that $L_{i}/L_{i-1}\simeq M_i$ for
$i=1,\cdots, t.$ Let $v_q=\sqrt q\in \bbc$ and $\cp$
 be the set of isomorphism classes of finite dimensional nilpotent  $\llz$-modules.
Then the Ringel-Hall algebra $\ch(\llz)$ of $\llz$ is by definition
the $\bbq(v)$-space with basis $\{u_{[M]}|[M]\in\cp\}$ whose
multiplication is given by
$$u_{[M]}u_{[N]}=\sum_{[L]\in\cp}g^{L}_{MN}u_{[L]}. $$
Note that $g^L_{M N}$  depends only on the isomorphism classes of $M,N$
and $L$, and for fixed isomorphism classes of $M,N$ there are only
finitely many isomorphism classes $[L]$ such that $g^L_{M N}\neq
0.$ It is clear that $\ch(\llz)$ is associative $\bbq(v_q)$-algebra
with unit $u_0$, where $0$ denotes the zero module.

The set of isomorphism classes of (nilpotent) simple
$\llz$-modules is naturally indexed by the set $I$ of vertices of
$Q.$ Then the Grothendieck group $G(\llz)$ of $\mod\llz$ is the
free Abelian group $\bbz I$. For each nilpotent $\llz$-module $M$,
the dimension vector $\udim M=\sum_{i\in I}(\dim M_i)i$ is an
element of $ G(\llz)$.  The Ringel-Hall algebra $\ch(\llz)$
is graded by $\bbn I,$ more precisely, by dimension vectors of
modules.

The Euler form $\lr{-,-}$ on $G(\llz)=\bbz I$ is defined by

\beq\hs \lr{\az,\bz}=\sum_{i\in I}a_ib_i-\sum_{\rho\in
H}a_{s(\rho)}b_{t( \rho)}\eeq for $\az=\sum_{i\in I}a_i i$ and
$\bz=\sum_{i\in I}b_i i$ in $\bbz I.$ For
any nilpotent $\llz$-modules $M$ and $N$ one has
$$\lr{\udim M, \udim
N}=\dim_{\fq}\hom_{\llz}(M,N)-\dim_{\fq}\ext_{\llz}(M,N).$$
The symmetric Euler form is defined as
$$(\az,\bz)=\lr{\az,\bz}+\lr{\bz,\az}\ \ \text{for}\ \
\az,\bz\in\bbz I.$$ This gives rise to a symmetric generalized
Cartan matrix $C=(a_{ij})_{i,j\in I}$ with $a_{ij}=(i,j).$ It is
easy to see that $C$ is independent of the field $\fq$ and the
orientation of $Q.$


The twisted Ringel-Hall algebra $\ch^*(\llz)$ is defined by
setting $\ch^*(\llz)=\ch(\llz)$ as $\bbq(v_q)$-vector space, but the
multiplication is defined by
$$u_{[M]}\ast
u_{[N]}=v_q^{\lr{\udim M, \udim N}}\sum_{[L]\in\cp}g^{L}_{MN}u_{L}.$$
Following [R3], for any $\llz$-module $M$, we denote $\lr{M}=v^{-\dim
 M+\dim\ed_{\llz}(M)}u_{[M]}.$  Note that $\{ \lr{M} \;|\; M \in
\cp\}$  a $\bbq(v_q)$-basis of $\ch^*(\llz)$.

The $\bbq(v_q)$-algebras $ \ch^*(\llz)$ and $\ch(\llz)$ depends on $q$.
We will use $ \ch_q^*(\llz)$ and $\ch_q(\llz)$ indicate the dependence
on $q$ when such a need arises.

\nd {\bf 1.3 A construction by Lusztig.} For any finite dimensional $I$-graded
$\fq$-vector space $V=\sum_{i\in I}V_i,$ let
$\bbe_V$ be the subset of $\oplus_{\rho\in H}\hom(V_{s(\rho)},V_{t(\rho)})$
defining  nilpotent
representations of $Q.$ Note that
$\be_V=\oplus_{\rho\in H}\hom(V_{s(\rho)},V_{t(\rho)})$ when $Q$
has no oriented cycles. The  group $G_V=\prod_{i\in
I}GL(V_i)$ acts naturally on  $\bbe_V$ by
$$(g,x)\mapsto g\bullet x=x'\ \ \text{where}\ \
x'_{\rho}=g_{t(\rho)}x_{\rho}g^{-1}_{s(\rho)}\ \ \text{for all}\ \
\rho\in H.$$ Let $\bbc_G(\bbe_V)$ be the space of $G_V$-invariant
functions $\bbe_V\ra\bbc.$ For $\gz\in\bbn I,$ we fix a $I$-graded
$\fq$-vector space $V_{\gz}$ with $\udim V_{\gz}=\gz.$ There is no
danger of confusion if we denote by $\bbe_{\gz}=\bbe_{V_{\gz}}$
and $G_{\gz}=G_{V_{\gz}}.$ For $\az,\bz\in\bbn I$ and $\gz=\az
+\bz,$ we consider the diagram
$$\bbe_{\az}\times\bbe_{\bz}\strl{p_1}{\longleftarrow}\bbe'\strl{p_2}{\lra}\bbe''\strl{p_3}{
\lra}\bbe_{\gz}.$$
Here $\bbe''$ is the set of all pairs $(x,W)$,
consisting of $x\in\bbe_{\gz}$ and an $x$-stable $I$-graded
subspace $W$ of $V_{\gz}$ with $\udim W=\bz$, and $\bbe'$ is the
set of all quadruples $(x,W,R',R'')$, consisting of $(x,W)\in\bbe''$
and two invertible linear maps $R':\fq^{\bz}\ra W$ and
$R'':\fq^{\az}\ra\fq^{\gz}/W.$ The maps are defined in obvious way as
follows: $p_2(x,W,R',R'')=(x,W),$
$p_3(x,W)=x,$ and $p_1(x,W,R',R'')=(x',x''),$ where
$x_{\rho}R'_{s(\rho)}=R'_{t(\rho)}x'_{\rho}$ and
$x_{\rho}R''_{s(\rho)}=R''_{t(\rho)}x''_{\rho}$ for all $\rho\in
H.$

For any map $p: X\rightarrow Y$ of finite sets,
$p^*:\bbc(Y)\rightarrow \bbc(X)$ is defined by
$p*(f)(x)=f(p(x))$ and $p_!: \bbc(X)\rightarrow \bbc(Y)$
is defined by $ p_!(h)(y)=\sum_{x \in p^{-1}(y)}h(x)$
integration along the fibers).
Given $f\in\bbc_G(\bbe_\az)$ and $g\in\bbc_G(\bbe_\bz)$,
there is a unique $h\in \bbc_G(\bbe'')$ such that $p_2^*(h)=p_1^*(f\times g).$ Then define
$$f\circ g=(p_3)_!(h)\in\bbc_G(\bbe_{\gz}).$$

 Let $${\wm}(\az,\bz)= \sum_{i\in
I}a_ib_i+\sum_{\rho\in H}a_{s(\rho)}b_{t( \rho)}.$$
We again
define the multiplication in the $\bbc$-space ${\bf
K}=\oplus_{\az\in\bbn I}\bbc_G(\bbe_{\az})$ by
$$f\ast g=v_q^{-{\wm}(\az,\bz)}f\circ g$$ for all $f\in\bbc_G(\bbe_\az)$ and
$g\in\bbc_G(\bbe_\bz).$ Then $({\bf K},\ast)$ becomes an
associative $\bbc$-algebra.

\mk\nd{\bf Convention.}  {\sl Although we are working over finite $ \fq$, we
will regularly use $G_V$ and $\bbe_V$ for the  algebraic group and
the algebraic variety which are defined over $\fq$ and use the features of
algebraic geometry without introducing extra notations,  i. e., the set
of $\fq$-rational points and the algebraic variety are denoted by
the same notation. This should not cause any confusion and  in
particular, the concept of $G_V$-orbits will be consistent in both cases due to
Lang's theorem for this group $G_V$ acting on $ \bbe_V$.
For $M \in \bbe_V$, we will use  $M$  to denote the representation of
$Q$ on $V$ defined by $M$.}

For $M\in\bbe_{\az}$, let  $\co_M\subset\bbe_{\az}$ be the
$G_{\az}$-orbit of $M.$ We take $\bi_{[M]}\in\bbc_G(V_{\az})$ to
be the characteristic function of $\co_M,$ and set
$f_{[M]}=v_q^{-\dim\co_M}\bi_{[M]}.$ We consider the subalgebra
$({\bf L},\ast)$ of $({\bf K},\ast)$ generated by $f_{[M]}$ over
$\bbq(v_q),$ for all $M\in\bbe_{\az}$ and all $\az\in\bbn I.$ In
fact ${\bf L}$ has a $\bbq(v_q)$-basis $\{f_{[M]}|M\in\bbe_{\az},
\az\in\bbn I\},$ since we have the relation $\bi_{[M]}\circ
\bi_{[N]}(W)=g^W_{M N}$ for any $W\in\bbe_{\gz}.$

\mk \nd{\bf Proposition~1.1} {\sl The linear map
$\varphi:({\bf L},\ast)\lra \ch^*(\llz)$ defined by
$$\varphi(f_{[M]})=\lr{M},\ \ \ \text{for all}\ [M]\in\cp$$
is an isomorphism of the associative $\bbq(v_q)$-algebras.}

\mk\nd {\bf Proof.} Note that $\varphi$ is a linear isomorphism.
For $[M],[N]\in\cp$ with $\udim M=\az$ and  $\udim
N=\bz$, since  $\bi_{[M]}\circ\bi_{[N]}=\sum_{[L]}g^L_{M
N}\bi_{[L]}$ in ${\bf L}$,  we have
\begin{eqnarray*}
f_{[M]}\ast f_{[N]}\!&=&
\sum_{[L]\in \cp}v_q^{-\dim\co_M-\dim\co_N-{\wm}(\az,\bz)+\dim\co_L}g^L_{M
 N}f_{[L]}.
\end{eqnarray*}
Note that   $\dim\co_M=\dim G_{\az}-\dim\ed_{\llz}(M)$  and
$\dim G_{\az+\bz}-\dim G_{\az}-\dim G_{\bz}= \lr{\az,\bz}+{\wm}(\az,\bz)$. In
$\ch^*(\llz)$ we have
\begin{eqnarray*}\lr{M}\!\ast\!\lr{N}
&=&v_q^{-\dim
M+\dim\ed_{\llz}(M) -\dim
N+\dim\ed_{\llz}(N)+\lr{\az,\bz}} u_{[M]}\circ u_{[N]}\\
&=&\sum_{L}v_q^{\dim\ed_{\llz}(M)+\dim\ed_{\llz}(N)-\dim\ed_{\llz}(L)+\lr{\az,\bz}}g^L_{M
N}\lr{L}\\
&=&\sum_{L}v_q^{\dim G_{\az}-\dim\co_M+\dim G_{\bz}-\dim\co_N-(\dim
G_{\az+\bz}-\dim\co_L)+\lr{\az,\bz}}g^L_{M N}\lr{L}\\
&=&\sum_{L}v_q^{-\dim\co_M-\dim\co_N+\dim\co_L-{\wm}(\az,\bz)}
g^L_{M N}\lr{L}. \ \Box
\end{eqnarray*}

\mk\nd{\bf 1.4} The free abelian group $G(\llz)=\bbz I$
 with the symmetric Euler form $(-,-)$ defined in~1.2 is a Cartan
datum in the sense of Lusztig [L5]. Associated to $(\bbz I,(-,-))$ is
the Drinfeld-Jimbo quantized enveloping algebra $U=U^-\otimes
U^0\otimes U^+$ defined over $\bbq(t),$ where $t$ is
transcendental over $\bbq.$ It is generated by  the Chevalley
generators $E_i, F_i, K_i^{\pm}$ $(i\in I)$ with respect to the
quantum Serre relations.  Let $\cz=\bbz[t,t^{-1}].$ The Lusztig form $U^+_{\cz}$ of
$U^+$ is the $\cz$-subalgebra in $U^+$ generated by $E_i^{(m)}=\frac{E_i^m}{[m]!}$ ($m\geq 0 $ and $i\in I$).
 For  $v=v_q \in \bbc$, let $ \cz_v$ be the subring of $ \mathbb{C}$
as the image of $\cz$ under the map $ \cz\rightarrow \mathbb{C} $
with $t\mapsto v$. Let $\cc^*(\llz)_{\cz_v}$ be the
$\cz_v$-subalgebra of  $\ch^*(\llz)$ generated by
$u_{[S_i]}^{(\ast m)}=\frac{u_{[S_i]}^{\ast m}}{[m]_v!}$ ($ i\in I$),
where
$$[n]=\frac{t^n-t^{-n}}{t-t^{-1}},\ [n]!=\Pi_{r=1}^{n}[r], \
\Bigl[\begin{array}{c}n\\ r\end{array}\Bigr]=\frac{[n]!}{[r]![n-r]!}$$ and $
[n]_v \in \cz_v $ is the image of $[n]$ in $ \cz_v$.

It follows from the works of Ringel [R1], Green [G],
and Sevenhant-Van den Bergh
[SV] that $\cc^*(\llz)_{\cz_v}$ is isomorphic to
$U^+_{\cz}\otimes_{\cz}\cz_{v}$ by sending $u_i^{(\ast  m)}$ to
$E_i^{(m)}$.

We will denote $\cc^*(\llz)_\cz $ for $ U^+_{\cz}$ and call it the
integral generic composition algebra. In fact, following
Ringel's point of view,  $ \cz$ can be identified with the subring of
$\prod_{q}\cz_{v_{q}}$ generated by $t^{\pm 1}=(v_{q}^{\pm 1})$ and $ \cc^*(\llz)_\cz
$ as a $ \cz$-subalgebra of $\prod_{q}\ch^*_q(\llz)$ generated by
$(u_{[S_i\otimes \fq]}^{(\ast m)}),$ $m\geq1.$ Here the product is
taken over all $q$ (though infinitely many will be enough.

In this paper,
 computations in $\prod_{q}\ch^*_q(\llz)$ will  performed  in each
component. When an expression in each component is
written as an element of $ \mathbb{Z}[v_q,v_q^{-1}]$ with coefficients
in $\mathbb{Z}$ independent of the choice of the field $
\mathbb{F}_q$, we say that the expression is invariant (or generic) as
$\mathbb{F}_q$ varies. In this case replacing  $v_q$ by $t$ will get a
formula in $\prod_{q}\ch^*_q(\llz)$. We will not repeat this replacement
each time and simply write $v=v_q$ and call it generic in this expression.
In stead of write $t$, we will also use $v$ and this will not cause any confusion.

There is bar involution $\ol{(\ )}:U^+\ra U^+$ (of $\bbz$-algebras)  defined by $\ol{t}=t^{-1},$ $\ol{E_i}=E_i$ and
$\ol{E_i^{(m)}}=E_i^{(m)}.$ Then  $\overline{U^+_{\cz}}=
U^+_{\cz}$.

\mk\nd{\bf 1.5}  In general, if we take a special value
$v=\sqrt{q}$ for the finite field $\fq,$ it is easy to see that

\mk\nd{\bf Lemma~1.2} {\sl Given any monomial $\fk{m}$  of
$u_{S_i}^{(m)},\  i\in I, m \in \mathbb{N}$ we have
$\fk{m}=\sum_{M\in\cp}f_{M,q}\lr{M}$ in $\ch^*(\llz)$  with
$f_{M,q}\in \cz_v$. Then for each $ M$, there is an integer $b$
such that  $v^bf_{M,q}\in \bbz[v]$ (the subring of algebraic
integers)  and $b$ is independent of $\fq.$} $\Box$

\bigskip

\centerline{\bf 2. The variety of representations}

\bigskip
\nd We need slightly more knowledge about the geometry of
representations of quivers over algebraically closed field
 $k=\ol{\bbf}_{q}$.  In this section  we only consider finite
quivers $Q$ without  oriented cycles. Take $\llz=kQ$ and all $\hom$ and $\ext$ are taken in $\llz$-mod.

\mk\nd{\bf 2.1}  For $\az\in\bbn I,$, the $I$-graded $k$-vector space $\oplus_{i\in I}k^{\az_i}$ defines the affine algebraic $k$-variety
$\bbe_{\az}$ on which  the algebraic group $G_{\az}$ acts
in a similar way as in~1.3. For any
$x\in\bbe_{\az}$ we have the corresponding representation $M(x)$
of $Q$ over $k.$  The following properties are well-known
(see[CB]).

\mk\nd {\bf Lemma~2.1} {\sl For any $\az\in\bbn I$ and
$M\in\bbe_{\az},$ we have
\begin{itemize}
\item[(1)]
  $\dim\bbe_{\az}-\dim\co_M=\dim\ed(M)-(\az,\az)/2=\dim\ext^1(M,M).$
\item[(2)] $\co_M$ is open in $\bbe_{\az}$ if and only if $M$ has no
self-extension.
\item[(3)] There is at most one orbit $\co_M$ in $\bbe_{\az}$ such that
$M$ has no self-extension.
\item[(4)] If $0\ra M\ra L\ra N\ra 0$ is a non-split exact sequence, then
$\co_{M\oplus N}\subseteq{\ol{\co}_{L}}\setminus\co_{L}$.
\item[(5)] If $\co_L$ is an orbit in $\bbe_{\az}$ of maximal dimension
and $L=M\oplus N,$ then $\ext^1(M,N)=0$. $\Box$
\end{itemize}}

\mk\nd For subsets $\ca\subset\bbe_{\az}$ and $\cb\subset\bbe_{\bz},$ we define
the extension set $\ca\star\cb$ of $\ca$ by $\cb$  to be
\begin{eqnarray*} \ca\star\cb&=&\{z\in\bbe_{\az+\bz}|\ \text{there exists an exact
sequence}\ \\
&& \quad 0\ra M(x)\ra M(z)\ra M(y)\ra 0 \ \text{with}\
x\in\cb,\ y\in\ca\}.
\end{eqnarray*}
Set
$\cdim\ca=\dim\bbe_{\az}-\dim\ca.$ It follows from [Re] that

\mk\nd{\bf Lemma~2.2} {\sl Given any  $\az, \bz\in\bbn I,$ if  $\ca\subset\bbe_{\az}$ and $\cb\subset\bbe_{\bz}$  are irreducible algebraic varieties and are stable
under the action of $G_{\az}$ and $G_{\bz}$ respectively, then
$\ca\star\cb$ is irreducible and stable under the action of
$G_{\az+\bz},$ too. Moreover,
$$\cdim \ca\star\cb =\cdim\ca+\cdim\cb-\lr{\bz,\az}+r,$$
where $0\leq r\leq\min\{\dim_k\hom(M(y),M(x))|y\in\cb,x\in\ca\}$.}
$\Box$

\mk\nd{\bf 2.2} For any $\az,\bz\in\bbn I,$ we consider the
diagram if algebraic $k$-varieties
$$\bbe_{\az}\times\bbe_{\bz}\strl{p_1}{\longleftarrow}\bbe'\strl{p_2}{\lra}\bbe''\strl{p_3}{
\lra}\bbe_{\az+\bz}$$ defined by a similar way
as in~1.3. It follows from the definition that
$\ca\star\cb=p_3p_2(p_1^{-1}(\ca\times \cb)$. Thus we have
$\overline{\ca}\star \overline{\cb} \subseteq \overline{\ca\star\cb}$ since $ p_1 $ is
a locally trivial fibration (see Lemma~2.3).
 For any $M\in\bbe_{\az}, N\in\bbe_{\bz}$ and
$L\in\bbe_{\az+\bz}$ we define
$$\bfz=p_2p_1^{-1}(\co_M\times\co_N), \ \bfz_{L,M,N}=\bfz\cap p_3^{-1}(L).$$
Then it follows from [L1] that

 \mk\nd{\bf Lemma~2.3} {\sl For
the diagram above and $M\in\bbe_{\az}, N\in\bbe_{\bz}$ and
$L\in\bbe_{\az+\bz},$ we have the following properties.
\begin{itemize}
\item[(1)] The
map $p_2$ is a principal $G_{\az}\times G_{\bz}$ fibration.
\item[(2)] The map $p_1$ is a locally trivial fibration with smooth
connected fibres of dimension
$$\sum_{i\in I}a_i^2+\sum_{i\in I}b_i^2 +{\wm}(\az,\bz).$$
\item[(3)] The map $p_3$ is proper.
\item[(4)] The variety $\bfz$ is smooth and irreducible of dimension
$$\dim\bfz=\dim(\co_M)+\dim(\co_N) +{\wm}(\az,\bz).$$
\item[(5)] If $L$ is an extension of $M$ by $N,$ then
$$\dim(\co_L)\leq\dim(\co_M)+\dim(\co_N) +{\wm}(\az,\bz).$$
\item[(6)] If $\co_L$ is dense in $p_3\bfz,$ then
$$\dim(\co_L)=\dim(\co_M)+\dim(\co_N) +{\wm}(\az,\bz)-\dim\bfz_{L,M,N}.$$
\item[(7)] Assume that $\ext(M,N)=0 $
and $\hom(N,M)=0.$ If $M'\in\ol{\co}_M$ and  $N'\in\ol{\co}_N$
such that either $M'\in\ol{\co}_M\setminus\co_M$ or
$N'\in\ol{\co}_N\setminus\co_N,$ then $X\in\ol{\co}_{M\oplus
N}\setminus\co_{M\oplus N}$ for all
$X\in\ol{\co}_{M'}\star\ol{\co}_{N'}$.  In particular,
$\dim\co_X<\dim\co_{M\oplus N}.$ $\Box$\end{itemize} }

 \mk\nd As a consequence of Lemma~2.2 we have

\nd{\bf Lemma~2.4} {\sl Given any two representations $M$ and $N$
of $Q$ over $k,$ if $\ext(M,N)=0,$ then
$\ol{\co}_M\star\ol{\co}_N=\ol{\co}_{M\oplus N},$ i.e.,
$\co_{M\oplus N}$ is open and dense  in
$\ol{\co}_M\star\ol{\co}_N.$} $\Box$

\mk\nd{\bf Lemma~2.5} {\sl Let $M,N,X\in\mod\llz.$ Then $\co_X$ is
open in $\co_M\star\co_N$ if and only if $\co_X$ is open in
$\ol{\co}_M\star\ol{\co}_N.$ In that case for any $Y\in
\ol{\co}_M\star\ol{\co}_N$ we have $\dim\co_Y\leq\dim\co_X.$}

\nd{\bf Proof.} This follows from  $ \co_X \subseteq \co_M\star \co_N
\subseteq \ol{\co}_M\star \ol{\co}_N \subseteq \ol{\co_M\star
\co_N}$ and Lemma~2.2. $\Box$

\bigskip

\centerline{\bf 3. The integral bases from the Kronecker quiver}

\bigskip

Most results in this section can be found in [Z] and [C] while
others can be found in [BK]. For completeness, we give some proofs
here.

\nd{\bf 3.1} Let $\fq$ be the finite field with $q$ elements and $Q$
the Kronecker quiver with $I=\{1,2\}$ and $H=\{\rho_1, \rho_2\}$
such that $s(\rho_1)=s(\rho_2)=2$ and $t(\rho_1)=t(\rho_2)=1$. Let $\llz=\fq Q$ be the path algebra. It is
known that the structure of the preprojective and preinjective
components of $\mod\llz$ is the same as those of $\mod kQ$ for $k$
being an algebraically closed field. However the regular components
of $\mod\llz$ is different with that of $\mod kQ.$

The set of dimension vectors of  indecomposable
representations is
$$\Phi^+=\{(l+1,l),(m,m),(n,n+1)|l\geq0,m\geq 1, n\geq 0\}. $$
The dimension vectors  $(n+1,n)$ and $(n,n+1)$ correspond to preprojective and
preinjective indecomposable representations respectively and are call real roots. For each real root $ \az$,
 there is
only one isoclass of indecomposable representation with dimension vector $\az$ which
will be denoted by $ V_{\az}$.  Define a total order $\prec$ on $\Phi^+$ by
\begin{eqnarray*} && (1,0)\prec\cdots\prec (m+1,m)
\prec (m+2,m+1)\prec\cdots\prec(k,k)\prec (k+1,k+1)\\
&& \prec\cdots\prec (n+1,n+2)\prec
(n,n+1)\prec\cdots\prec (0,1).
\end{eqnarray*}

The strong representation-directed property implies that there is
no non-zero homomorphism from an indecomposable module of
dimension vector $\az$ to an indecomposable module of dimension
 vector $\beta$ if $\beta\prec \az$. This property will  used
 frequently in the computation.

Any $\llz$-module is given by the date $(V_1,V_2;\sigma,\tau),$
where $V_1$ and $V_2$ are finite dimensional vector space over
$\fq,$ $\sigma$ and $\tau$ are $\fq$-linear maps from $V_2$ to
$V_1.$

\medskip
\nd {\bf Proposition 3.1.} {\sl The isomorphism classes of the
regular quasi-simple modules in $\mod\llz$ are indexed by
$spec(\fq[x]).$ That is, each regular quasi-simple module is
isomorphic to $(V_1,V_2;\sigma,\tau),$ where
$V_1=V_2=\fq[x]/(p(x))$ for an irreducible polynomial $p(x)$ in
$\fq[x],$ $\sigma$ is the identity map and $\tau$ is given by the
multiplication by $x,$ except $(\fq,\fq;0,1)$ which corresponds to the zero ideal. }
\medskip

\nd{\bf 3.2} In this
section, let $\cp$ be the set of  isomorphism classes of finite
dimensional $\llz$-modules, $\ch=\ch_q$ be the Ringel-Hall algebra of
$\llz$ over $\bbq(v),$ where $v^2=q,$  and $\ch^*$ be the twisted
 form of $\ch.$
If $\bd\in\bbn I$ be a dimension vector, we set in $\ch$
$$R_{\bd}=\sum_{\substack{[M]\in\cp, M \text{ regular}\\
\udim M=\bd}}u_{[M]}.$$
 For an element
$x=\sum_{[M]\in\cp}c_{[M]}u_{[M]}\in\ch,$ we call $u_{[M]}$ to be
a (non-zero) term of $x$ if $c_{[M]}\neq 0.$ Furthermore,
$$R(x)=\sum_{[M]\in\cp, M\text{ regular}}c_{[M]}u_{[M]}$$
 is called the regular part of $x.$ According to our notation, we
denote  $u_{\az}=u_{[V_{\az}]}$ for $ \az=(n-1,n)$ or $(n,n+1)$ being
real roots.

Let  $\az_1=(1,0)$ and $\az_2=(0,1)$ be the simple root vectors. The
orientation of $Q$ implies $ \lr {\az_1, \az_2}=0 $ and $ \lr{\az_2,
\az_1}=-2$. Thus for  $ \delta=(1,1)$ we have  $\lr{\dz,\az_1}=-1,$ $\lr{\az_1,\dz}=1,$
$\lr{\dz,\az_2}=1$ and $\lr{\az_2,\dz}=-1.$

\medskip
\nd{\bf 3.3} In this section,  the multiplication in $ \ch$ will be simply written as  $ xy$
instead of $x\circ y$.  The following can be computed easily  as in~[Z].

\nd{\bf Lemma~3.2.} {\sl Let $i$ and $j$ be two positive integers.
Then
$$u_{(j-1,j)}u_{(i,i-1)}=R(
u_{(j-1,j)}u_{(i,i-1)})+q^{i+j-2}u_{(i,i-1)}u_{(j-1,j)}.\Box$$}

\medskip
\nd{\bf Lemma~3.3} {\sl \begin{eqnarray*}R_{\dz}&=&u_{(0,1)}u_{(1,0)}-
u_{(1,0)}u_{(0,1)},\\
u_{(n+1,n)}&=&\frac{1}{q+1}(R_{\dz}u_{(n,n-1)}-qu_{(n,n-1)}R_{\dz}),\\
u_{(n,n+1)}&=&\frac{1}{q+1}(u_{(n-1,n)}R_{\dz}-qR_{\dz}u_{(n-1,n)}). \Box\end{eqnarray*}}

\medskip
\nd{\bf Lemma~3.4} {\sl Let $i$ and $j$ be two positive integers
and $n=i+j-1.$ Then
$$R(u_{(j-1,j)}u_{(i,i-1)})=R(u_{(n-1,n)}u_{(1,0)})=R(u_{(0,1)}u_{(n,n-1)}).\Box $$}

\medskip
\nd{\bf Lemma~3.5} {\sl Let $m,n\geq1.$ Then
\begin{eqnarray*}u_{(m-1,m)}R_{n\dz}&=&\sum_{0\leq i\leq
n}\frac{q^i-q^{n+1}}{1-q}R_{i\dz}u_{(m+n-i-1,m+n-i)},\\
R_{n\dz}u_{(m,m-1)}&=&\sum_{0\leq i\leq
n}\frac{q^i-q^{n+1}}{1-q}u_{(m+n-i,m+n-i-1)}R_{i\dz}.\Box \end{eqnarray*}}

\medskip
\nd{\bf 3.4} We will introduce a new set of elements in $\ch^*$ to describe a basis that resembles PBW basis for enveloping algebra of a Lie algebra. We give here some  quantum commutative relations in
$\ch$ and in $\ch^{*}.$  We define (cf. 1.2)
\begin{eqnarray*}E_{(n+1,n)}&=&\lr{u_{(n+1,n)}}=v^{-2n}u_{(n+1,n)},\\
E_{(n,n+1)}&=&\lr{u_{(n,n+1)}}=v^{-2n}u_{(n,n+1)}.
\end{eqnarray*}
We will call $E_1=E_{(1,0)}, E_2=E_{(0,1)}$ the Chevalley generators.
For $n\geq 1,$ define in $\ch^*$
$$\te_{n\dz}=E_{(n-1,n)}\ast E_1-v^{-2}E_1\ast E_{(n-1,n)}.$$

In the following we give a sequence of computations we will
need. Most of them  are known.
\medskip

\nd{\bf Lemma~3.6} {\it
$\te_{n\dz}=v^{-3n+1}R(u_{(n-1,n)}u_{(1,0)}).$}

\nd{\bf Proof.} By taking  $u_1=u_{(1,0)}$ we have
\begin{eqnarray*}\te_{n\dz}
&=&v^{-2(n-1)}(v^{\lr{(n-1)\dz+\az_2,\az_1}}u_{(n-1,n)}u_1-v^{-2}v^{\lr{\az_1,(n-1)\dz+\az_2}}u_1
u_{(n-1,n)})\\
&=&v^{-3n+1}(u_{(n-1,n)}u_1-v^{2(n-1)}u_1 u_{(n-1,n)})\\
&=&v^{-3n+1}R(u_{(n-1,n)}u_1)\ \ \ \ \text{by Lemma~3.2}. \quad \Box\end{eqnarray*}

\medskip

\nd{\bf Lemma~3.7} {\sl In $\ch^{*}$ we have
\begin{eqnarray*} &&[\te_{\dz},E_{(n+1,n)}] = [2]_v E_{(n+2,n+1)},\\
 && [E_{(n,n+1)},\te_{\dz}]=[2]_v E_{(n+1,n+2)}.
\end{eqnarray*}}

\nd{\bf Proof.} We only check the first equation. By definition and
Lemma~3.3, we
\begin{eqnarray*} {} [\te_{\dz},E_{(n+1,n)}]
&=&v^{-2(n+1)}v^{\lr{\dz,n\dz+\az_1}}R_{\dz}u_{(n+1,n)}-v^{-2(n+1)}v^{\lr{n\dz+\az_1,\dz}}
u_{(n+1,n)}R_{\dz}\\
&=&v^{-2(n+1)}v^{-1}((q+1)u_{(n+2,n+1)}+q
u_{(n+1,n)}R_{\dz})-v^{-2(n+1)}v u_{(n+1,n)}R_{\dz}\\
&=&v^{-2(n+1)}(v+v^{-1})u_{(n+2,n+1)}=[2]_v E_{(n+2,n+1)}.\Box
\end{eqnarray*}

\medskip

\nd{\bf Lemma~3.8} {\sl $E_{(2,1)}\ast E_1=v^2 E_1\ast
E_{(2,1)} $ and $
E_2\ast E_{(1,2)}=v^2 E_{(1,2)}\ast E_2.$}

\nd{\bf Proof.} Let $M=V_{(1,0)}\oplus V_{(2,1)}$. Then $E_{(2,1)}\ast
E_1=v^{-2}v^{\lr{\dz+\az_1,\az_1}}u_{(2,1)}u_1=v^2u_{\bz} $
and
 $E_1\ast E_{(2,1)}=v^{-2}v^{\lr{\az_1,\dz+\az_1}}u_1
u_{(2,1)}=u_{\bz}.$ This proves the first equality. The second
equality follows from a similar computation. $\Box$

\mk\nd{\bf Lemma~3.9} {\sl For any non-negative integers $r$ and
$s,$ we have in $\ch^{*}$
$$\te_{(r+s+1)\dz}=E_{(r,r+1)}\ast
E_{(s+1,s)}-v^{-2}E_{(s+1,s)}\ast E_{(r,r+1)}.$$}

\nd{\bf Proof.}  Using Lemma~3.2, Lemma~3.4, and Lemma~3.6, we have
\begin{eqnarray*}&&E_{(r,r+1)}\ast E_{(s+1,s)}-v^{-2}E_{(s+1,s)}\ast
E_{(r,r+1)}\\
&&\quad =v^{-3(r+s)-2}u_{(r,r+1)}u_{(s+1,s)}-v^{-(r+s)-2}u_{(s+1,s)}u_{(r,r+1)}\\
&&\quad =v^{-3(r+s)-2}(R(u_{(r,r+1)}u_{(s+1,s)})+q^{r+s}u_{(s+1,s)}u_{(r,r+1)})
-v^{-(r+s)-2}u_{(s+1,s)}u_{(r,r+1)}\\
&&\quad =v^{-3(r+s)-2}R(u_{(r,r+1)}u_{(s+1,s)})=v^{-3(r+s)-2}R(u_{(r+s,r+s)}u_1)=\te_{(r+s+1)\dz}.
\Box\end{eqnarray*}
\mk\nd{\bf Lemma~3.10} {\sl There exist
$a^{(r)}_{h}(t)\in\bbz[t,t^{-1}]$ for all $r\in\bbn\setminus\{0\}$
and $h\in\{0,1,\cdots,\lfloor{\frac{r}{2}}\rfloor\}$ such that for all $n>m$ in $\bbn$,
\begin{eqnarray*} E_{(n+1,n)}\ast
E_{(m+1,m)}&=&\sum^{\lfloor{\frac{n-m}{2}}\rfloor}_{h=0}a_{h}^{(n-m)}(v)E_{(m+h+1,m+h)}\ast
E_{(n-h+1,n-h)}\\
E_{(m,m+1)}\ast
E_{(n,n+1)}&=&\sum^{\lfloor{\frac{n-m}{2}}\rfloor}_{h=0}a_{h}^{(n-m)}(v)E_{(n-h+1,n-h)}\ast
 E_{(m+h,m+h+1)}
\end{eqnarray*}}
\nd{\bf Proof.} Using the strong representation-directed property, we have
\begin{eqnarray*} E_{(n+1,n)}\ast E_{(m+1,m)}&=&v^{-2(n+m)}v^{\lr{n\dz+\az_1,m\dz+\az_1}}u_{(n+1,n)}
u_{(m+1,m)}\\
&=&v^{-3n-m+1)}\sum^{\lfloor{\frac{n-m}{2}}\rfloor}_{h=0}g_{V_{(n+1,n)}
V_{(m+1,m)}}^{M_h}u_{[M_h]}.
\end{eqnarray*}
where $M_h=V_{(n-h+1,n-h)}\oplus V_{(n+h+1,n+h)}$.  Since $ n-h\geq
m+h$,  the strong
representation-directed property again implies
$$ E_{(m+h+1,m+h)}\ast
E_{(n-h+1,n-h)} =v^{-3m-n-2h+1}g_{V_{(m+h+1,m+h)}V_{(n-h+1,n-h)}}^{M_{h}}u_{[M_h]}.$$
Thus $g_{V_{(m+h+1,m+h)}V_{(n-h+1,n-h)}}^{M_{h}}u_{[M_h]}=v^{3m+n-1+2h}=E_{(m+h+1,m+h)}\ast
E_{(n-h+1,n-h)} $. Substitution implies that $ a^{(r)}_h(t)=t^{-2(r-h)}\in \bbz[t, t^{-1}]$.   To verify the second identity, one uses the strong representation-directed property again and carry out similar computation. The computation will give the same $a^{(r)}_{h}(t)=t^{-2r+2h}$. Thus the same set of $ a_{h}^{(r)}(t) $ works for both identities.\Box

\mk For $k\geq 0$, we inductively define
$$E_{0\dz}=1,\ \ E_{k\dz}=\frac{1}{[k]}\sum_{s=1}^k
v^{s-k}\te_{s\dz}\ast E_{(k-s)\dz}.$$
\mk\nd{\bf Lemma~3.11} {\sl We have $E_{k\dz}=v^{-2k}R_{k\dz}.$}

\nd{\bf Proof.} If  $k=1,$ $E_{\dz}=\te_{\dz}=v^{-2}R_{\dz}.$  We
assume that the assertion is true for all numbers $t<k.$ Then using
Lemma~3.6, and [Z] (Lem~3.7, Thm~4.1,  Lem~4.7),  we have

\begin{eqnarray*} E_{k\dz}
&=&\frac{1}{[k]}\sum^k_{s=1}v^{s-k}v^{-3s+1}R(u_{(s-1,s)}u_1)\ast
v^{-2(k-s)}R_{(k-s)\dz}\\
&=&\frac{1}{[k]}\sum^k_{s=1}v^{-3k+1}R(u_{(s-1,s)}u_1)\ast
R_{(k-s)\dz}\\
&=&\frac{1}{[k]}\sum^k_{s=1}v^{-3k+1}a_s(R_\dz,R_{2\dz},\cdots,R_{s\dz})R_{(k-s)\dz}\\
&=&\frac{v^{-3k+1}}{[k]}\frac{1-q^k}{1-q}R_{k\dz}=v^{-2k}R_{k\dz}. \quad \Box
\end{eqnarray*}

\mk\nd{\bf Lemma~3.12} {\sl For $m,n\in\bbn$ we have in $\ch^{*}$
\begin{eqnarray*} E_{n\dz}\ast
E_{(m+1,m)}&=&\sum_{k=0}^{n}[n+1-k]E_{(m+n+1-k,m+n-k)}\ast
E_{k\dz};\\
E_{(m,m+1)}\ast
E_{n\dz}&=&\sum_{k=0}^{n}[n+1-k]E_{k\dz}\ast E_{(m+n-k,m+n-k+1)}.\end{eqnarray*}}

\nd{\bf Proof.} Again it only needs to verify the first equation. By
Lemma~3.11, we have
\begin{eqnarray*} E_{n\dz}\ast E_{(m+1,m)}&=&v^{-2n}R_{n\dz}\ast
v^{-2m}u_{(m+1,m)}
=v^{-2(n+m)}v^{-n}R_{n\dz}u_{(m+1,m)}\\
&=&v^{-3n-2m}\sum_{k=0}^n
\frac{q^k-q^{n+1}}{1-q}u_{(m+n-k+1,m+n-k)}R_{k\dz} \ \ \ \text{(by
Lemma~3.5)} \\
&=&\sum_{k=0}^n\frac{v^{-3n-2m}}{1-v^2}(v^{2k}-v^{2n+2})v^{2k}v^{-k}v^{2(n+m-k)}
E_{(m+n-k+1,m+n-k)}\ast E_{k\dz}\\
&=&\sum_{k=0}^n [n+1-k]
E_{(m+n-k+1,m+n-k)}\ast E_{k\dz}. \Box \end{eqnarray*}

\mk\nd{\bf 3.5} Let  $\cl_v$ be the
$\cz_v=\bbz[v,v^{-1}]$-subalgebra of $ \ch_q^*$  generated by the set
$$\{E^{(*s)}_{(m+1,m)}, E_{k\dz}, E^{(*t)}_{(n,n+1)}|m\geq 0, n\geq
0, s\geq 1, t\geq 1,k\geq 1\}.$$ It contains the divided powers
$E_1^{(*s)}, E_2^{(*t)},$ $s,t\in\bbn,$ of the Chevalley generators.
We have obtained an integral $\cz_v$-basis of $\cl_v$ consisting of the monomials
$$\{\prod_{m\geq 0}E^{(s_m)}_{(m+1,m)}\prod_{k\geq 1}E^{r_k}_{k\dz}\prod_{n\geq 0} E^{(t_n)}_{(n,n+1)}\;|\;m\geq 0, n\geq
0, s_m\geq 0, t_n\geq 0,k\geq 1, r_k\geq 0\}$$ with the product taken with respect to the
order given in~3.1 and there are only finitely many non-zero $s_m$,
$t_n$, and $r_k$ in each monomial.  This follows easily  from
the facts: (1) the commutation relations in the above lemmas imply
that the $\cz_v$-span of the monomials above is closed under the
multiplication in $\ch^*$ and that  $\cl$ contains all  monomials
we defined above; (2) those monomials
are linearly independent over $\cz$ (even over $\bbq(v)$) by the
definition of Ringel-Hall algebras.

\mk\nd{\bf Remark} The formulae in the lemmas are unchanged when we vary
$v=\sqrt{q}$. The statement of the lemmas can be stated in
$ \prod_{q}\ch^*_q$ with $v$ replaced by $t=(v_q)$ in
$\prod_{q}\cz_q$ and $E_{*,*}$ replaced by $E_{(*,*)}=(E_{(*,*),q})$.
We then denote $ \cl$ as the $\cz=\bbz[t,t^{-1}]$-algebra with
a $\cz$-basis consisting monomials described above.

 As
remarked in 1.4,   Lusztig's integral  $\cz$-form $\cc^*_{\cz}$, which we called the generic composition algebra, can be viewed as a
$\cz$-subalgebra of of $\prod_{q}\ch^*_{q}$ by the Ringel-Green
theorem (see [G],[R1]).  Using this identification, we can view
$\cc^*_{\cz}$ as a $ \cz$-subalgebra of  $ \cl$. In the rest
of this section, we will construct a $\cz$-basis of $\cc^*_{\cz}$.


For any $n> m\geq 0$, let  $P_{(n,m)}$ (resp.  $ I_{(m,n)}$) be an
isomorphism class of preprojective (resp. preinjective)
modules  with $ \udim P_{(n,m)}=(n,m)$ (resp. $\udim
I_{(m,n)}=(m,n)$).  In  the following formulas, the summation is taken
over all nonzero preprojective and preinjective representations of
the indicated
dimension vectors.

\mk\nd {\bf Lemma~3.13}\setcounter{equation}{0} {\sl In the following
formulas all $P$ and $ I$ are non-zero.
 \begin{multline}  E_2^{(\ast  n)}\ast E_1^{(\ast  (n+1))}=E_{(n+1,n)}+\sum_{1\leq l\leq
n}v^{-l-1}E_{(n-l+1,n-l)}\ast E_{l\dz}  \\
+\sum_{\substack{0\leq l\leq n-1
\\ p\geq 1, s\geq0,t\geq 0 \\ s+t+l+(p-1)=n}}
v^{-\dim\ed(P)-dim\ed(I)}v^{-p(l+t)-(s+l)(p-1)}\lr{P_{(s+p,s)}}\ast
E_{l\dz}\ast\lr{I_{(t,t+p-1)}}; \end{multline}
  \begin{multline}  E_2^{(\ast  (n+1))}\ast E_1^{(\ast  n)}=E_{(n,n+1)}+\sum_{1\leq l\leq
n}v^{-l-1}E_{l\dz}\ast E_{(n-l,n-l+1)}  \\
+\sum_{\substack{0\leq l\leq n-1\\
p\geq 1, s\geq0,t\geq 0 \\ s+t+l+(p-1)=n}}v^{-\dim\ed(P)-dim\ed(I)}
v^{-p(l+s)+(t+l)(p-1)}\lr{P_{(s+p-1, s)}}\ast
E_{l\dz}\ast\lr{I_{(t,t+p)}}; \end{multline}

   \begin{multline} E_2^{(\ast  n)}\ast E_1^{(\ast  n)}=E_{n\dz}\\+
\sum_{\substack{0\leq l\leq n-1, p\geq 1\\ s\geq0,t\geq 0,
s+t+l+p=n}}v^{-\dim\ed(P)-\dim\ed(I)} v^{-p(s+2l+t)}\lr{P_{(s+p,s)}}\ast
E_{l\dz}\ast\lr{I_{(t,t+p)}}. \end{multline}}

\nd{\bf Proof.} We only verify (1) and
others can be verified in a similar way. We have the following
relation in $\ch$ (see [R3]).
$$u_2^n u_1^{n+1}=\psi_n(q)\psi_{n+1}(q)(u_{(n+1,n)}+\sum_{1\leq
l\leq n}u_{(n-l+1,n-l)}R_{l\dz} +\sum_{\substack{0\leq l\leq n-1, p\geq
1 \\ s\geq0,t\geq 0, s+t+l+(p-1)=n}}u_{[P]}R_{l\dz}u_{[I]})$$
where $P$ is a non-zero preprojective module with $\udim
P=(s+p,s)$ and $I$ is a non-zero preinjective module with
$\udim I=(t,t+p-1)$ and
$$\psi_n(q)=\frac{(1-q)\cdots (1-q^n)}{(1-q)^n}.$$
Then by a routine calculation according to the relation in~3.2, we
have the relation  (1).  $\Box$

\mk\nd{\bf 3.6} Note  that the dimensions of $P,I,\ed I$ and $\ed
P$ over $\fq$ are invariant as  $\fq$ varies.  By induction using
Lemma~3.13,  the set
$$\{E_{(m+1,m)}, E_{k\dz}, E_{(n,n+1)}|m\geq 0, n\geq
0,k\geq 1\}$$ is contained in $\cc^*_\cz.$ If $M$ is
indecomposable preprojective or preinjective, then, by [R3],
$$\lr{u_{[M]}}^{(*s)} =\lr{u_{[ M^{\oplus s}]}}\in\cc^*(\llz)_{\cz}\  \text{for  any }\ s
\geq 1.$$
Using this and the strong representation-directed property on preprojectives and preinjectives, we have, for $0\leq n_1<n_2<\cdots <n_l$ and $s_1, s_2,  s_l \geq 1 $, we have
\begin{eqnarray*} E_{(n_1+1,n_1)}^{(\ast s_1)}\ast \cdots E_{(n_l+1,n_l)}^{(\ast s_l)}
&=&v^a \lr{[\oplus_{i=1}^{l}V_{(n_i+1,n_i)}^{\oplus s_i}]},\\
E_{(n_l,n_1+1)}^{(\ast s_l)}\ast \cdots E_{(n_1+1,n_1)}^{(\ast s_1)}
&=&v^b \lr{[\oplus_{i=1}^{l}V_{(n_i,n_i+1)}^{\oplus s_i}]},
\end{eqnarray*}
where $a$ and $b$ are integers depends only on the sequences
$n_1< \cdots <n_l$ and $s_1, s_2,  s_l \geq 1 $. Hence the subset
$$\{E^{(s)}_{(m+1,m)}, E_{k\dz}, E^{(t)}_{(n,n+1)}|m\geq 0, n\geq
0, s\geq 1, t\geq 1,k\geq 1\}$$ is also contained  in $\cc^*_\cz.$ Therefore,
$\cl=\cc^*_\cz.$

 Let ${\bf P}(n)$ the set of all partitions of $n$. Recall that
there are no nontrivial extensions between homogeneous regular
representation. For any
$w=(w_1,w_2,\cdots,w_m)\in{\bf P}(n),$  we define
$$E_{w\dz}=E_{w_1\dz}\ast E_{w_2\dz}\ast\cdots\ast E_{w_m\dz}.$$

\mk\nd{\bf Proposition~3.14} {\sl The set
$$\{\lr{P}\ast
E_{w\dz}\ast\lr{I}\|P\in\cp\text{preprojective},w\in{\bf
P}(n),I\in\cp\text{preinjective},n\in\bbn\}$$ is a $\cz$-basis
of $\cc^*_\cz.$} $\Box$

\mk\nd{\bf Remarks.} (1) It has been proved by Zhang in [Z] that these
monomials are $\bbq(v)$-bases of $U^+,$ then improved by Chen in
[C] that they are $\cz$-bases of $U^+_{\cz}.$

(2) It is not difficult to see that the root vectors provided here
exactly correspond to the root vectors of $\uq2$ provided by
Damiani in [Da] and by Beck in [Be].

(3) It can be proved that the set in Proposition~3.14 is an
integral basis of $\cc^*$ over $\ca=\bbq[v,v^{-1}]$ by an easier
way, see the proofs of Proposition~7.2 and~7.3 below.

\bigskip

\centerline{\bf 4. A bar-invariant basis from the Kronecker
quiver}

\bigskip
With the PBW type basis constructed for $\cc^*_{\cz}$ We now can
construct a bar-invariant basis following the approach in [L1,
7.8-7.11].

\nd{\bf 4.1}   In this section, we work in $\cc^*=\cc^*_{\cz}$.
Recall from~3.1 that $\Phi^+$ is  the positive root system of
$\hat{sl}_2.$  A function $\wc:\Phi^+\ra\bbn$ is called
support-finite, if  $\wc(\az)\neq 0$ only for finitely many
$\az\in\Phi^+.$ Let $\bbn_f^{\Phi^+} $ be the set of all
support-finite $\bbn$-valued functions.  We will use the order in
$\Phi^+$ given in~3.1.

For $\wc\in\bbn_f^{\Phi^+},$ if $ \{\az\in\Phi^+|\wc(\az)\neq
0\}=\{\bz_1\prec\bz_2\prec\cdots\prec\bz_k\},$  we set
$$E^{\wc}=E_{\bz_1}^{(\ast  \wc(\bz_1))}\ast
E_{\bz_2}^{(\ast  \wc(\bz_2))}\ast\cdots\ast E_{\bz_k}^{(\ast
\wc(\bz_k))},$$
where $ E_{\bz_k}^{(\ast \wc(\bz_k))}\triangleq
 E_{\bz_k}^{\ast \wc(\bz_k)}$ if $\beta_i=m\delta$.
Then Proposition~3.14 is equivalent to the statement:
\begin{center} {\it The set $\{E^{\wc}|\wc\in\bbn_f^{\Phi^+}\}$ is a
$\cz$-basis of $\cc^*.$}
\end{center}

 For
$\bd=(d_1,d_2)\in\bbn_f^2,$ we denote
$$E(\bd)=E_2^{(\ast  d_2)}\ast E_1^{(\ast  d_1)}.$$
Similarly we defined
$$E(\wc)=E(\wc(\bz_1)\bz_1)\ast E(\wc(\bz_2)\bz_2)\ast \cdots\ast
E(\wc(\bz_k)\bz_k).$$
Note that $E(\wc)\in\cc^*$ since it is  a monomial on the Chevalley generators
$E_1$ and $E_2$ in the form of divided powers. Moreover, by
definition, we  $ \overline{E(\bd)}=E(\bd)$. Thus $
\overline{E(\wc)}=E(\wc)$.

\mk\nd{\bf 4.2} The rest  of this section is devoted to
giving  a triangular relation between the PBW-basis and the monomial
basis.

For any $\wc\in\bbn_f^{\Phi^+}$ we assume that
$E^{\mathbf{c}}=\lr{P}*E_{\omega\delta}*\lr{I},$ where $P$ is a
preprojective module and $I$ is a preinjective module. For any partition $\omega=(\omega_1, \omega_2\dots, \omega_m)$, write
$$E_{w\dz}=E_{w_1\dz}\ast
E_{w_2\dz}\ast\cdots\ast E_{w_m\dz}=\sum_{V}a^\omega_{V}\lr{V}$$
over a fixed field $ \bbf_q$. We choose\footnote{This selection is not unique, in fact we
may require that $V_{m\dz}$ is defined over the prime field and  absolutely indecomposable in a
homogeneous tube.}
$V_{\omega\delta}$ to be a module such that
$$\dim
\co_{V_{\omega\delta}}=max\{\dim \co_{V}\;|\; a^{\omega}_{V}\neq 0 \}.$$
 Set
    $$V_{\wc}=P\oplus V_{\omega\delta}\oplus I,\ \ \text{and}
    \ \ \co_{\wc}=\co_{P\oplus V_{\omega\delta}\oplus I}.$$

For any $\wc\in\bbn_f^{\Phi^+}$ and any real root $\az\in\Phi^+$, define
$u_{\wc(\az)\az}=u_{[V_{\az}\oplus\cdots\oplus V_{\az}]},$ where
$V_{\az}$ is the indecomposable representation with $\udim
V_{\az}=\az$.

\mk\nd{\bf Lemma~4.2} {\sl For any  $\wc\in\bbn_f^{\Phi^+}$ and
any real root $\bz\in\Phi^+$, we have  in $\cc^*$
$$
E(\wc(\bz)\bz)=\lr{u_{\wc(\bz)\bz}}+
\sum_{\substack{\wc'\in\bbn_f^{\Phi^+}\\ \dim\co_{\wc'}<\dim\co_{\wc(\bz)\bz}}}v^{-\lz(\wc')}E^{\wc'},$$
where $\lz(\wc')\in\bbn.$}

\nd{\bf Proof.} Let $\wc(\bz)\bz=(m,n).$ In $\ch_q$ (for any fixed $\bbf_q$) we have
$$u_2^m u_1^n=\psi_m(q)\psi_n(q)\sum_{\udim N=(m,n)}u_{[N]}.$$
By Lemma~2.1
$$u_2^m
u_1^n=\psi_m(q)\psi_n(q)u_{[V_{\bz}\oplus\cdots\oplus
V_{\bz}]}+\psi_m(q)\psi_n(q) \sum{u_{[P]}R_{l\dz}u_{[I]}},$$ where
$P$ is preprojective,  $I$ is preinjective,  $\udim
P+l\dz+\udim I=(m,n)$, and $\dim \co_{P\oplus V_{l\dz}\oplus
I}<\dim\co_{\wc(\bz)\bz}.$ Although the number of terms of $R_{l\dz}=\sum_{[M]}u_{[M]}$ in $\ch_q$ depends on $q$, Lemma~3.11 shows that $R_{l\dz}$ has a generic form in $\cc*_{\cz}$ with each component in $\ch_q$ being $R_{l\dz}$.  Then in $\cc^*_\cz,$
\begin{eqnarray*} u_2^{(\ast  m)}* u_1^{(\ast  n)}&=&\frac{v^{m(m-1)/2}v^{n(n-1)/2}}{[m]![n]!}v^{-2mn}u_2^m u_1^n \\
&=&v^{m^2-m+n^2-n-2mn}u_{[V_{\bz}\oplus\cdots\oplus V_{\bz}]}+v^{m^2-m+n^2-n-2mn}
\sum u_{[P]}R_{l\dz}u_{[I]}\\
 &=&\lr{u_{\wc(\bz)\bz}}+
\sum_{\substack{\wc'\in\bbn_f^{\Phi^+}\\ \dim\co_{\wc'}<\dim\co_{\wc(\bz)\bz}}}v^{-\lz(\wc')}E^{\wc'}.\
\Box
\end{eqnarray*}

\mk\nd{\bf Lemma~4.3} {\sl Let $\az,\bz\in\Phi^+$ be real roots
and $\az\prec \bz.$ We have in $\cc^*$
$$E(\az)\ast
E(\bz)=\lr{u_{[V_\az]}}\ast\lr{u_{[V_\bz]}}
+\sum_{\substack{\wc\in\bbn_f^{\Phi^+}\\ \dim
\co_{\wc}<\dim\co_{V_{\az}\oplus
V_{\bz}}}}h_{\wc}^{\az\bz}E^{\wc},$$ where $ h_{\wc}^{\az\bz}\in
\cz.$}

\nd{\bf Proof.} By Lemma~3.13 we have
\begin{eqnarray*}E(\az)&=&\lr{u_{\az}}
+\sum_{\substack{\wc'\in\bbn_f^{\Phi^+}\\
\dim\co_{\wc'}<\dim\co_{V_{\az}}}}
v^{-\lz(\wc')}E^{\wc'},\\
E(\bz)&=&\lr{u_{\bz}}+\sum_{\substack{\wc''\in\bbn_f^{\Phi^+}\\ \dim\co_{\wc''}<\dim\co_{V_{\bz}}}}
v^{-\lz(\wc'')}E^{\wc''}.
\end{eqnarray*}
 Since $\hom_{\Lambda}(V_\beta, V_\alpha)=\ext_\Lambda(V_{\alpha},
V_{\beta})=0$ and  $\dim
Z_{V_{\az}\oplus V_{\bz}, V_{\az},V_{\bz}}=0$, by Lemma~2.3(7),
$$\dim\co_{\wc}<\dim\co_{V_{\az}\oplus V_{\bz}}$$
for any extension $V_{\wc}$  of $V_{\wc'}$ by $V_{\wc''}$ with the
property: $$\co_{\wc'}\sset\overline{\co_{V_{\az}}}\setminus
\co_{V_{\az}},\text{or}\
\co_{\wc''}\sset\overline{\co_{V_{\bz}}}\setminus \co_{V_{\bz}}.$$
Therefore, the conclusion follows from  Proposition~3.14.  $\Box$

\mk\nd{\bf Lemma~4.4} {\sl Let $\az=(n+1,n)$, $\bz=(l,l)=l\dz$, and $
\gamma=(m,m+1)$ be in $\Phi^+$. The for all $s\geq 1$
\setcounter{equation}{0}
\begin{eqnarray}
E(s\az)\ast E(\bz)&=&\lr{u_{[sV_{\az}]}}\ast
E_{l\dz}+\sum_{\substack{\wc\in\bbn_f^{\Phi^+}\\ \dim\co_{\wc}<\dim\co_{sV_{\az}\oplus
V_{l\dz}}}}h(\wc)E^{\wc}, \\
E(\bz)\ast E(s\gamma)&=&E_{l\dz}\ast\lr{u_{[sV_\gamma]}}
+\sum_{\substack{\wc'\in\bbn_f^{\Phi^+}\\ \dim\co_{\wc'}<\dim\co_{V_{l\dz}\oplus
sV_{\gamma}}}}h(\wc')E^{\wc'}.\end{eqnarray}
Here $h(\wc'), \ h(\wc)\in\cz.$}

\nd{\bf Proof.} Using Lemma~2.3(7) and Lemma~3.13(3), the proof
is same as that of Lemma~4.3. $\Box$

\mk\nd{\bf Lemma~4.5} {\sl Let $V$ be an indecomposable regular
module with $\udim V=n\delta $. $M=P\oplus M'\oplus I$ with $P\neq
0,M',I\neq 0$ are respectively preprojective,regular and preinjective
modules and $\udim M=n\delta $. Then $\dim \co_{V}>\dim \co_{M}$.}

\nd{\bf Proof.} By Lemma~2.1(1), we only need to prove that $\dim \ed(V)< \dim \ed(M) .$ By
Proposition~3.1, we have $\dim
\ed(V)=n.$ Suppose
$$P=P_1\oplus P_2\oplus\cdots\oplus P_r,\ \ \text{and}\ \ I=I_1\oplus I_2\oplus\cdots\oplus I_t,$$
where $P_j,$ and $I_j(j\geqslant 1)$ are respectively indecomposable
preprojective and preinjective modules with
$\udim P_j=(n_j+1,n_j)$,  $\udim I_j=(m_j,m_j+1),$ and
$\udim M'=s\delta$. Thus
$r=t$ and  $n=\sum_{j=1}^{t}(n_j+1)+s+\sum_{j=1}^{t}(m_j).$  Note that
\begin{align*} \dim \ed(P)\geqslant & t, \quad \dim \ed(M')=s, \quad \dim \ed(I)\geqslant t,\\
\dim \hom(P,M')=&\lr{\udim P,\udim
M'}=st, \quad \
\dim \hom(M',I)=st,\\
\dim\hom(P,I)=&\lr{\udim P,\udim I}=t(\sum_{j=1}^tn_j+\sum_{j=1}^tm_j). \end{align*}
Using the direct sum decomposition of $ M$, one computes
\begin{eqnarray*}
\dim \ed(M) \geqslant
t+s+t+st+t(\sum_{j=1}^tn_j+\sum_{j=1}^tm_j)+st\geqslant
2t+\sum_{j=1}^tn_j+\sum_{j=1}^tm_j>n.
\end{eqnarray*}
 This implies that $\dim \co_{V}>\dim \co_{M}$. $\Box$

 \mk\nd{\bf Lemma~4.6} {\sl Let $n\geqslant1,m\geqslant 1.$ Then
 \begin{eqnarray*}
E(n\delta)*E(m\delta)&=&E_{n\delta}*E_{m\delta}
+\sum_{\substack{\wc\in\bbn_f^{\Phi^+}\\ \dim\co_{\wc}<\dim\co_{V_{(n,m)\delta}}=
\dim\co_{V_{(n+m)\delta}} }}h(\wc)E^{\wc},\end{eqnarray*}
 where $V_{(n,m)\dz}$ is defined in~{\rm 4.2} and $h(\wc)\in\cz.$}

\nd{\bf Proof.} By Lemma~3.13(3), we have
\begin{eqnarray*}
E(n\delta)&=&E_{n\dz}+
\sum_{P\neq 0,I\neq 0}v^{l(\lr{P}\ast E_{l\dz}\ast\lr{I})}
\lr{P}\ast E_{l\dz}\ast\lr{I},\\
E(m\delta)&=&E_{m\dz}+
\sum_{P\neq 0,I\neq 0}v^{l(\lr{P}\ast E_{l\dz}\ast\lr{I})}
\lr{P}\ast E_{l\dz}\ast\lr{I},\end{eqnarray*}
 where $l(\lr{P}\ast E_{l\dz}\ast\lr{I})\in \mathbb{Z}.$

We then have
$$E(n\delta)*E(m\delta)=E_{n\delta}*E_{m\delta}+\sum h(\wc)E^{\wc}.$$
To prove the lemma, it is sufficient to prove that $V_{\wc},$ which
is defined in~4.2, is decomposable. This is easy to see from the
structure of the  AR-quiver of Kronecker quiver.
 $\Box$

\nd{\bf Remark.} By Lemma~4.6 we can get
$$E(\omega\delta)=E(\omega_1\delta)*\cdots*E(\omega_m\delta)=E_{\omega\delta}+
\sum_{\dim\co_{\wc}<\dim\co_{\omega\delta} }h(\wc)E^{\wc},$$ where
$h(\wc)\in\cz.$

 \mk\nd Let
$\vhi:\bbn_f^{\Phi^+}\ra\bbn^2$ be defined by
$\vhi(\wc)=\sum_{\az\in\Phi^+}\wc(\az)\az.$ Then for any
$\bd\in\bbn^2,$ $\vhi^{-1}(\bd)$ is a finite set. We define a
(geometric) order in $\vhi^{-1}(\bd)$ as follows: {\it
$\wc'\preceq\wc$ if and only if $\wc'=\wc$ or $\wc'\neq\wc$ but
$\dim\co_{\wc'}<\dim\co_{\wc}.$}\footnote{This is independent of the choices of $V_{\wc'}$ and $V_{\wc}$ as in~4.2 such that
$\co_{V_{\wc'}}\subseteq\overline\co_{V_\wc}\setminus\co_{V_\wc}.$ }

{}From Lemma~2.3(7) and above lemmas, we may summarize our results of
this subsection as follows resembling [L1, 7.8].

\mk\nd {\bf Proposition~4.7} {\sl For any $\wc\in\bbn_f^{\Phi^+}$, we
have
$$E(\wc)=\sum_{\wc'\in\vhi^{-1}(\vhi(\wc))}h^{\wc}_{\wc'}E^{\wc'}$$
such  that {\rm(1)} $h^{\wc}_{\wc'}\in\cz,$ {\rm(2)}
$h^{\wc}_{\wc}=1,$ {\rm(3)} if $h^{\wc}_{\wc'}\neq 0$  then
$\wc'\preceq\wc,$ (4) $\overline{E(\wc)}=E(\wc)$. }  $\Box$

\mk\nd For any $ \wc, \wc' \in \bbn_f^{\Phi^+}$ we define  $\oz^{\wc}_{\wc'}\in\cz$ such that
$$\ol{E^{\wc}}=\sum_{\wc'\in\bbn_f^{\Phi^+}}\oz^{\wc}_{\wc'}E^{\wc'}.$$
The following Proposition resembles [L1, Prop. 7.9].

\mk\nd{\bf Proposition~4.8} {\sl $\oz^{\wc}_{\wc}=1$ and, if $\oz^{\wc}_{\wc'}\neq 0$ and
$\wc'\neq\wc$ then $\wc'\prec\wc.$}

\nd{\bf Proof.} Using  $ \overline{E(\wc)}=E(\wc)$ and
the fact that $ \{E^{\wc}\;|\; \wc\in\bbn_f^{\Phi^+}\}$ is a
 $\cz$-bases of $\cc^*$, we have
$$h^{\wc}_{\wc''}=\sum_{\wc'}\ol{h^{\wc}_{\wc'}}\oz^{\wc'}_{\wc''},\
\ \text{for} \ \wc,\wc''\in\vhi^{-1}(\bd).$$ By Lemma~4.5,
the matrices $(h^{\wc}_{\wc''})$ as well as
$(\ol{h^{\wc}_{\wc'}}),$ where the index set is $\vhi^{-1}(\bd),$
are triangular with $1$ on diagonal. Hence, by the equation above,
the matrix $(\oz^{\wc}_{\wc''})$ has the same property.  $\Box$

\mk\nd Consider the bar involution $\ol{(\ )}:\cc^*\ra\cc^*.$ For any $\wc\in\bbn_f^{\Phi^+},$
$$E^{\wc}=\ol{\ol{E^{\wc}}}=\ol{\sum_{\wc'}\oz^{\wc}_{\wc'}E^{\wc'}}=\sum_{\wc',\wc''}
\ol{\oz^{\wc}_{c'}}\oz^{\wc'}_{\wc''}E^{\wc''}.$$
implies the orthogonal relation
$$\sum_{\wc'}\ol{\oz^{\wc}_{\wc'}}\oz^{\wc'}_{\wc''}=\dz_{\wc
\wc''}.$$
Therefore one can solve uniquely the system of equations
$$\zeta^{\wc}_{\wc'}=\sum_{\wc'\preceq\wc''\preceq\wc}\oz^{\wc''}_{\wc'}
\ol{\zeta^{\wc}_{\wc''}}$$ with unknowns
$\zeta^{\wc}_{\wc'}\in\bbz[v^{-1}],$ $\wc'\preceq\wc$ and $\wc'\ ,
\wc\in\vhi^{-1}(\bd),$ such that
$$\zeta^{\wc}_{\wc}=1\ \ \it{and}\ \ \zeta^{\wc}_{\wc'}\in v^{-1}\bbz[v^{-1}]\ \ \it{for \ all}\
\wc'\prec\wc.$$
\mk\nd For any $\bd\in\bbn^2$ and $\wc\in\vhi^{-1}(\bd)$,   we set
$$\ce^{\wc}=\sum_{\wc'\in\vhi^{-1}(\bd)}\zeta^{\wc}_{\wc'}E^{\wc'}\\
\text{ and } \
\bj=\{\ce^{\wc}|\wc\in\vhi^{-1}(\bd),\bd\in\bbn^2\}.$$ Let
$$\cl=\text{span}_{\bbz[v^{-1}]}\{E^{\wc}|\wc\in\bbn_f^{\Phi^+}\}$$
We verify the following two properties of $\bj.$ The first is
\begin{eqnarray*} \ol{\ce^{\wc}}&=&\sum_{\wc'}\ol{\zeta^{\wc}_{\wc'}}\ol{E^{\wc'}}
=\sum_{\wc'}\ol{\zeta^{\wc}_{\wc'}}\sum_{\wc''}\oz^{\wc'}_{c''}E^{\wc''}\\
&=&
\sum_{\wc''}(\sum_{\wc'}\ol{\zeta^{\wc}_{\wc'}}\oz^{\wc'}_{c''})E^{\wc''}
=\sum_{\wc''}\zeta^{\wc}_{\wc''}E^{\wc''}=\ce^{\wc}.\end{eqnarray*}
 So the
elements $\ce^{\wc}$ are bar-invariant. The second, obviously the
set $\bj$ is a $\bbz[v^{-1}]$-basis of the lattice $\cl.$
Therefore we have

\mk\nd{\bf Proposition~4.9} {\sl The set $\bj$ is a  basis of
$\cc^*_{\cz}.$} which satisfies that  $\ol{\ce^{\wc}}=\ce^{\wc}$
 and
 $\pi(\ce^{\wc})=\pi(E^{\wc})$ for any $\ce^{\wc}\in\cl,$ where $\pi$ is  the canonical projection
 $\cl\ra\cl/v^{-1}\cl.$

\bigskip

\centerline{\bf 5. The integral and canonical bases arising from a
tube }
\bigskip

The main results we present in this section are taken from [DDX],
in which the canonical bases of $U_q(\hat{sl}_n)$ and
$U_q(\hat{gl}_n)$ are constructed by a linear algebra method from
the category of finite dimensional nilpotent representations of a
cyclic quiver, i.e, from a tube. However in an preliminary version
of the present paper we assumed the existence and
the structure of Lusztig's canonical basis for the composition
algebra of a tube from [L3] and [VV].

\medskip

\nd{\bf 5.1} Let $\dt=\dt(n)$ be the cyclic quiver with vertex set
$\dt_0=\bbz/n\bbz=\{1,2.,\cdots,n\}$ and arrow set $\dt_1=\{i\ra
i+1|i\in\bbz/n\bbz\}.$ We consider the category $\ct=\ct(n)$ of
finite dimensional nilpotent representations of $\dt(n)$ over
$\fq.$ For the reason of the shape of its Auslander-Reiten quiver,
$\ct(n)$ is called a tube of rank $n.$ Let $S_i,\ i\in\dt_0$ be
the irreducible objects in $\ct(n)$ and $S_i[l]$ the (unique) absolutely
indecomposable object in $\ct(n)$ with top $S_i$ and length $l.$ Note that $S_i[l]$ is independent of $q$.
Again in this section, we let $\cp$ be the set of isomorphism
classes of objects in $\ct(n)$, $\ch$ the Ringel-Hall algebra
of $\ct(n)$, $\ch^*$ the twisted Ringel-Hall algebra, and $\bl$ the
Lusztig form of the Hall algebra of $\ct(n)$ (cf. 1.3). Because the Hall
polynomials always exist in this case (see [R2]), we may regard
the algebras $\ch$, $\ch^*$ and $\bl$ in their generic form. So
they all are defined generically over $\bbq(t),$ where $t$ is an
indeterminate. By Proposition~1.1, we may identify $\bl$ with
$\ch^*$ via the morphism $\varphi.$

In this section, all properties we obtain are
generic and independent of the base field $\fq,$ although the computations will be performed over $\fq$ (for each $q$). We will omit  the subscript $q$
for simplicity. Since the number $n$ is
fixed, sometimes it is omitted too, e.g., $\ct=\ct(n)$.

\medskip

\nd{\bf 5.2} Let $\Pi$ be the set of $n$-tuples of partitions
$\pi=(\pi^{(1)},\pi^{(2)},\cdots,\pi^{(n)})$ with each component  $
\pi^{(i)}=(\pi^{(i)}_1\geq \pi^{(i)}_2\geq \dots)$ being a partition
of an integer. For each $ \pi\in \Pi$, we define an object in $\ct$
$$M(\pi)=\bigoplus_{\substack{i\in \dt_0 \\  j\geq 1}}  S_i[{ \pi}^{(i)}_j].$$
  In this way we obtain a
bijection  between $\Pi$ and the set $\cp.$ We will simply denote
by $u_{\pi}, \pi\in\Pi$ for $u_{[M(\pi)]}$ in $\ch.$

An $n$-tuple $\pi=(\pi^{(1)},\pi^{(2)},\cdots\pi^{(n)})$ of
partition in $\Pi$ is called aperiodic (in the sense of Lusztig
[L3]), or separated (in the sense of Ringel [R2]), if for each
$l\geq 1$ there is some $i=i(l)\in\dt_0$ such that
${\pi}_j^{(i)}\neq l$ for all $j\geq 1.$ By $\Pi^a$ we
denote the set of aperiodic $n$-tuples of partitions. An object
$M$ in $\ct$ is called aperiodic if $M\simeq M(\pi)$ for some
$\pi\in\Pi^a.$ For any dimension vector $\az\in\bbn^n (=\bbn I)$,
we let
$$\Pi_{\az}=\{\lz\in\Pi|\udim M(\lz)=\az\}\ \ \text{and}\ \
\Pi^a_{\az}=\Pi^a\cap \Pi_{\az}.$$

Given any two modules $M,N$ in $\ct,$ there exists a unique (up to
isomorphism) extension $L$ of $M$ by $N$ with minimal
$\dim\ed(L)$ [Re]. This extension $L$ is called the generic extension
of $M$ by $N$ and is denoted by $L=M\diamond N.$ If we define the
operation in $\cp$ by $[M]\diamond [N]=[M\diamond N],$ then
$(\cp,\diamond)$ is a monoid with identity $[0].$

Let $\ooz$ be the set of all words on the alphabet $\dt_0.$
For each $w=i_1i_2\cdots i_m\in\ooz,$ we set
$$M(w)=S_{i_1}\diamond S_{i_2}\diamond\cdots\diamond S_{i_m}.$$
Then there is a unique $\pi\in\Pi$ such that $M(\pi)\simeq M(w),$
we define $\wp(w)=\pi.$ It has been proved in [R2] that
$\pi=\wp(w)\in\Pi^a$ and $\wp$ induces a surjection
$\wp:\ooz\twoheadrightarrow\Pi^a.$

We have a (geometric) partial order on $\cp,$ or equivalently in $\Pi,$ as
follows: for $\mu,\lz\in\Pi,$ $\mu\pe\lz$ if and only if
$\co_{M(\mu)}\subseteq\ol{\co}_{M(\lz)},$ or equivalently,
$\dim\hom(M,M(\lz))\leq\dim\hom(M,M(\mu))$ for all modules $M$ in
$\ct.$

For each module $M$ in $\ct$ and integer $s\geq 1,$,   we
denote  by $sM$ thedirect sum of  $s$ copies of $M.$ For $ w \in \ooz$, write $w$
in  tight
form  $w=j_1^{e_1} j_2^{e_2}
\cdots j_t^{e_t}\in\ooz$ with $j_{r-1}\neq j_r$ for all $r$ and define
 $\mu_r\in\Pi$
such that $M(\mu_r)=e_r S_{j_r}$. For any  $\lz\in\Pi_{\sum_{r=1}^t
e_r j_r}$,  write $g^{\lz}_w$ for the Hall polynomial $g^{M(\lz)}_{M(\mu_1),\dots,M(\mu_t)}$.
A word $w$ is called distinguished if the Hall polynomial
$g^{\wp(w)}_w=1.$ This means that $M(\wp(w))$ has a unique reduced
filtration of type $w,$ i.e., a filtration
$$M(\wp(w))=M_0\supset M_1\supset\cdots\supset M_{t-1}\supset
M_t=0$$ with  $M_{r-1}/M_r\simeq e_r S_{j_r}$ for all
$r.$

\mk\nd{\bf Proposition 5.1} {\sl For any $\pi\in\Pi^a,$ there
exists a distinguished word $w_{\pi}=j_1^{e_1} j_2^{e_2} \cdots
j_t^{e_t}\in \wp^{-1}(\pi)$ in tight form.}

\mk In $\ch^*$, let $u_i^{(\ast  m)}=E_i^{(\ast  m)}= \frac{u_i^{\ast
m }}{[m]!}, i\in\dt_0, m\geq 1.$ The $\cz$-subalgebra
$\cc^*=\cc^*_{\cz}$ of $\prod_{q}\ch^*_q$ generated by $u_i^{(\ast  m)}, i\in\dt_0,
m\geq 1$, is  the twisted composition algebra of $\ct$ (cf. 1.4).

\mk\nd{\bf 5.3} For each $w=j_1^{e_1} j_2^{e_2} \cdots
j_t^{e_t}\in\ooz$ in tight form, define in  $\cc^*$ the monomial
$${\fk{m}}^{(w)}=E_{j_1}^{(\ast  e_1)}\ast\cdots\ast
E_{j_t}^{(\ast  e_t)}.$$

For each $\pi\in\Pi^a_{\az},$ we from now on fix a
distinguished word $w_{\pi}\in\wp^{-1}(\pi).$  Thus we have a section
 $\cd=\{w_{\pi}|\pi\in\Pi^a\}$ of $ \wp $ over $\Pi^a$. $\cd$  is
called a section of distinguished words  in [DDX].

For each $\pi \in \Pi^a$ with the fixed distinguished word $w_\pi=j_1^{e_1} j_2^{e_2} \cdots
j_t^{e_t}$ in tight form, define  $L_0=e_{j_1}S_{j_1},L_1=e_{j_1}S_{j_1}\diamond
e_{j_2}S_{j_2},L_2=L_1\diamond
e_{j_3}S_{j_3},\cdots,L_{t-1}=L_{t-2}\diamond e_{j_t}S_{j_t} $.  Set $\az=\udim L_{t-1}$.  As
$L_i$ is the generic extension of $L_{i-1}$ by
$e_{j_{i+1}}S_{j_{i+1}}$ and thus $\dim\ed(L_i)$ is minimal, we
have $M({\pi})\simeq L_{t-1}$. Since
$$1=g^{\pi}_{w_{\pi}}=g^{L_1}_{e_{j_1}S_{j_1},e_{j_2}S_{j_2}}g^{L_2}_{L_1,e_{j_3}S_{j_3}}\cdots
g^{\pi}_{L_{t-2},e_{j_t}S_{j_t}},$$
we get
$g^{L_i}_{L_{i-1},e_{j_{i+1}}S_{j_{i+1}}}=1,1\leq i\leq t-2.$
Furthermore, by Lemma~2.3(6) and Proposition~1.1, we have
$$\lr{L_{i-1}}\ast\lr{e_{j_{i+1}}S_{j_{i+1}}}=
\lr{L_i}+\sum_{X,\dim\mathcal{O}_X<\dim\mathcal{O}_{L_i}}a_X\lr{X},$$
with $a_X\in\cz.$ Recall from 1.2 that $\lr{M}=v^{-\dim M+\dim\ed(M)}u_{[M]}.$ Thus
$${\fk{m}}^{(w_{\pi})}=\lr{M(\pi)}
 +\sum_{\lz\prec\pi}\xi^{\lz}_{w_{\pi}}\lr{M(\lz)},$$
where $\xi^{\lz}_{w_{\pi}}\in\cz$. Note that
$\xi^{\lz}_{w_{\pi}}\neq 0$ implies  $\udim M(\lz)=\udim M(\pi)=\az.$
 Although $\lr{M}$ are in $ \ch^*$, they are not necessarily in $\cc^*$. Define  $E_{\pi}$ inductively by the relation (noting that $ v^2=q$ in each component)
$$E_{\pi}={\fk{m}}^{(w_{\pi})}-\sum_{\lz\prec \pi,
\lz\in\Pi^a_{\az}}v^{-\dim M(\pi)+\dim\ed M(\pi)+\dim
M(\lz)-\dim\ed M(\lz)}g^{\lz}_{w_{\pi}}(v^2) E_{\lz}.$$
if $ \pi \in \Pi_{\az}^{a}$ is minimal, then $E_{\pi}=\fk{m}^{(w_{\pi})} \in \cc^*$. By inductions on the partial order, we have $ E_{\lz}\in \cc^*$ for all $ \lz \in \Pi^a$.
Therefore
we have the relations
$$E_{\pi}=\lr{M(\pi)}+\sum_{\lz\in\Pi_{\az}\setminus\Pi^a_{\az},
\lz\prec \pi}\eta^{\pi}_{\lz}\lr{M(\lz)}$$ with
$\eta^{\pi}_{\lz}\in\cz.$

\mk\nd{\bf Proposition 5.2} {\sl Let $\cd=\{w_{\pi}|\pi\in\Pi^a\}$
be a section of distinguished words of $\ooz$ over $\Pi^a.$ Then
both $\{\fk{m}^{(w_\pi)}|\pi\in\Pi^a\}$ and $\{E_{\pi}|\pi\in\Pi^a\}$
are  $\cz$-bases of $\cc^*_{\cz}.$ Furthermore, for any $ \pi \in
\Pi^a_{\az}$,
$${\fk{m}}^{(w_{\pi})}=E_{\pi}+\sum_{\lz\in\Pi^a_\az, \lz\prec\pi}
v^{-\dim M(\pi)+\dim\ed M(\pi)+\dim M(\lz)-\dim\ed
M(\lz)}g^{\lz}_{w_{\pi}}(v^2)E_{\lz}.$$ }

\nd{\bf Remark.} The definition of the basis
$\{E_{\pi}|\pi\in\Pi^a\}$ depends on  the choice of the
section $\cd$ of distinguished words, but eventually it has been
proved in [DDX] that it is independent of the selection of the
sections of distinguished words.

We will call $\{\fk{m}^{(w_\pi)}|\pi\in\Pi^a\}$  a monomial
$\cz$-basis of $\cc^*_{\cz}$ and $\{E_{\pi}|\pi\in\Pi^a\}$ as a
``PBW''-basis of $\cc^*_{\cz}.$ With the triangular relation
between the two bases, we can follow the approach
of Lusztig [L1, 7.8-7.11], as we did in Section~4, to
obtain the canonical bases $\{\ce_{\pi}|\pi\in\Pi^a\}$ of
$\cc^*_{\cz}$  in the sense of [L1, 3.1] by
$$\ce_{\pi}=\sum_{\lz\pe\pi, \lz\in\Pi^a_{\az}}p_{\lz \pi}E_{\lz}, \
\ \text{for}\ \ \pi\in\Pi^a_{\az}, \az\in\bbn^n,$$ with $p_{\lz
\lz}=1$ and $p_{\lz \pi}\in v^{-1}\bbz[v^{-1}]$ for $\lz\prec\pi.$

\bigskip
\centerline{\bf 6. Integral bases arising from preprojective and
preinjective components}
\bigskip
In this  section consider a connected  tame quiver $Q$ without
oriented cycles. For the preprojective and preinjective components, the argument  in this
section is essentially as same as in the case of finite type.

\mk \nd{\bf 6.1}  Let $U$ be
the quantized affine enveloping algebra associated to the quiver $Q$, with the
Chevalley generators: $E_i, F_i $ and $K_i^{\pm}.$ Lusztig in [L5]
introduced the symmetries $T''_{i,1}: U\ra U$ for $i\in I,$
as algebra  automorphisms of $U$ defined by  relations:
\begin{eqnarray*}
T''_{i,1}(K_{\bz})&=&K_{s_i(\bz)},\quad T''_{i,1}(E_i)=-F_iK_i,\ \ \quad T''_{i,1}(F_i)=-K_iE_i,\\
T''_{i,1}(E_j)&=&\sum_{r+s=-a_{ij}}(-1)^rv^{-r}E_i^{(s)}E_jE_i^{(r)}\
\ \text{for}\ j\neq i\ \text{in}\ I,\\
T''_{i,1}(F_j)&=&\sum_{r+s=-a_{ij}}(-1)^rv^{r}F_i^{(r)}F_jF_i^{(s)}\
\ \text{for}\ j\neq i\ \text{in}\ I.
\end{eqnarray*}
Here $a_{ij}=(i,j)$ for $i, j\in I,$ and $\bz\in\bbz I$ and $s_i
(\bz)=\bz-(\bz,i)i.$  For each $i\in I$,  define
$$ U^+[i]=\{x\in U^+|T''_{i,1}(x)\in U^+\},$$
which is subalgebra $U^+$. Then $T''_{i,1}: U^+[i]\ra U^+[i]$ is
an automorphism. Moreover, if we consider the Lusztig form
$U^+_{\cz}$ and let $U^+_{\cz}[i]=U^+_{\cz}\cap U^+[i],$ then
$T''_{i,1}:U^+_{\cz}[i]\ra U^+_{\cz}[i]$ is an automorphism.

\mk\nd{\bf 6.2} We define $\sz_iQ$ to be the quiver obtained from
$Q$ by reversing the direction of every arrow connected to the
vertex $i.$ If $i$ is a sink of $Q,$ one may define the BGP
reflection functor (see[BGP] or[DR]):
$$\sz_i^+:\mod\llz\lra\mod\sz_i\llz$$
where $\llz=\fq(Q)$ and $\sz_i\llz=\fq(\sz_iQ)$ are path algebras.
Therefore we have an algebra
homomorphism: $$\sz_i: \ch^*(\llz)[i]\lra\ch^*(\sz_i\llz)[i]$$
defined by
$$\sz_i(u_{[M]})=u_{[\sz^+_i(M)]}\ \ \text{for any }\
M \in\mod\llz[i].$$
Here $\mod\llz[i]$ is the subcategory of
all representations which do not have $S_i$ as a direct summand
and $\ch^*(\llz)[i]$ is the subalgebra of $\ch^*(\llz)$ generated
by $u_{[M]}$ with $M \in\mod\llz[i].$ Note  that
$\cc^*(\llz)_{\cz}$ is canonically isomorphic to
$\cc^*(\sz_i\llz)_{\cz}$ by fixing the Chevalley generators which
correspond to the simple modules of $\llz$ and $\sz_i\llz$
respectively. Furthermore, we may regard that the functor
$\sz_i^{+}$ induces the homomorphism:
$$\sz_i:\cc^*(\llz)_{\cz}[i]\lra \cc^*(\llz)_{\cz}[i],$$
where
$\cc^*(\llz)_{\cz}[i]=\{x\in\cc^*(\llz)_{\cz}|\sz_i(x)\in\cc^*(\llz)_{\cz}\}.$
 It is known that $\sz_i=T''_{i,1}$ under the identification
 $\cc^*(\llz)=U^+$ (for example, see[XY]).

Dually, if $i$ is a source of $Q$, we have the similar results.

  We call an indecomposable $\llz$-module $M$
exceptional if $\ext_\llz^1(M,M)=0.$
It is proved in [CX] that
$$\lr{s M}\in\cc^*(\llz)_{\cz}\  \text{for  any }\ s \geq 1$$ if $M$ is
exceptional. In fact,
$$\lr{M}^{(\ast  s)}
=\frac{1}{[s]!}v^{-s\dim M+s}u_{[M]}^{\ast s}=\frac{1}{[s]!}
v^{-s\dim M+s}(v^{\binom{s}{2}}\psi_s(q) )u_{[sM]} =v^{-s\dim
M+s^2}u_{[s M]}=\lr{s M}.
$$

\mk\nd We denote by $Prep$ and $Prei$, respectively,  the isomorphism classes of
indecomposable preprojective and
preinjective $\llz$-modules. In
particular,  $\cc^*_{\cz}$ contains the set
$$\{\lr{u_{[s M]}}|M\  \text{is \ indecomposable in \ } Prep\
\text{or}\ Prei\ \text{and}\ s\geq1\}.$$

\mk\nd{\bf 6.3} Let $i_m,\cdots,i_1$ be an admissible sink
sequence of $Q,$ i.e., $i_m$ is a sink of $Q$ and for any
$1\leq t<m,$ the vertex $i_t$ is a sink for the orientation
$\sz_{i_{t+1}}\cdots\sz_{i_m}Q.$ Let $M$ be  in
$Prei.$ Then there exists an admissible sink sequence of $Q$ such
that
$$M=\sz_{i_1}^+ \cdots \sz_{i_m}^+(S_{i_{m+1}}),$$
where $S_{i_{m+1}}$ is a simple representation in $\mod\sz_{i_m}
\cdots \sz_{i_1}\llz.$

\mk\nd{\bf Lemma~6.1} {\sl Let $M$ be an indecomposable preinjective
representation. Then
$$\lr{M}=T''_{i_1,1}\cdots T''_{i_m,1}(E_{i_{m+1}}),$$ where
$M=\sz_{i_1}^+ \cdots \sz_{i_m}^+(S_{i_{m+1}}),$ for an admissible
sink sequence $i_m,\cdots,i_1$ of $Q.$}

\nd{\bf Proof.} See [R3].$\Box$

\mk\nd Since $Prei$ is representation-directed, we can total order on  the set
$$\Phi^+_{Prei}=\{\cdots,\bz_3,\bz_2,\bz_1\}$$ of  all positive real roots
appearing in $Prei$ with  $\{\cdots,M(\bz_3), M(\bz_2),
M(\bz_1)\}$ being the corresponding  indecomposable $\llz$-modules   such that
$$\hom(M(\bz_i),M(\bz_j))\neq 0\ \text{implies}\ \bz_i\preceq \bz_j\
\text{and}\ i\geq j.$$ Then such an ordering has the property
$$\lr{\bz_i,\bz_j}>0\ \text{implies}\ \bz_i\preceq\bz_j\
\text{and}\  i\geq j$$ and
$$\lr{\bz_i,\bz_j}<0 \ \text{implies}\  \bz_j\prec \bz_i\ \text{and}\
i<j$$ and $$\ext(M(\bz_i),M(\bz_j))=0 \ \text{for} \ i\geq j.$$
Therefore $\bz_i\preceq\bz_j$ if and only if $i\geq j.$
\mk\nd Similarly, since $Prep$ is representation-directed, we define a
total ordering on the set
$$\Phi^{+}_{Prep}=\{\az_1,\az_2,\az_3,\cdots\}$$ of
 of all positive real roots
appearing in $Prep,$ with   $\{M(\az_1), M(\az_2), M(\az_3),\cdots
\}$ be the corresponding  indecomposable modules in $Prep$   such that
$$\hom(M(\az_i),M(\az_j))\neq 0\ \text{implies}\ \az_i\preceq \az_j\
\text{and}\ i\leq j.$$ Then such an ordering has the property
$$\lr{\az_i,\az_j}>0\ \text{implies}\ \az_i\preceq\az_j\
\text{and}\  i\leq j$$ and
$$\lr{\az_i,\az_j}<0 \ \text{implies}\  \az_j\prec \az_i\ \text{and}\
j<i$$ and $$\ext(M(\az_i),M(\az_j))=0 \ \text{for}\ i\leq j.$$

We denote by $\bbn_f^{Prei}$ the set of all support-finite functions
$\bb: \Phi^+_{Prei}\ra\bbn$. Each $ \bb \in \bbn_f^{Prei}$ defines  a
preinjective representation
$$M(\bb)=\bigoplus_{\bz_i\in \Phi^+_{Prei}}\bb(\bz_i)M(\bz_i)$$
and  any preinjective representation is
isomorphic to one of the form.  By Ringel (Proposition 1' of
[R3]) we have

 \mk\nd{\bf Lemma~6.2} {\sl For any $\bb\in\bbn_f^{Prei},$
 $$\lr{M(\bb)}=\lr{\bb(\bz_{i_m})M(\bz_{i_m})}\ast\cdots\ast
 \lr{\bb(\bz_{i_1})M(\bz_{i_1})},$$
 where $\{\bz_{i_m}\prec \bz_{i_{m-1}}\prec\cdots\prec \bz_{i_1}\}
 $ are those $\bz\in \Phi^+_{Prei}$ such that $\bb(\bz)\neq0.$} $\Box$

 \mk\nd Thus, by~6.2,  $\lr{M(\bb)}\in\cc^*_{\cz}$ for all $\bb\in\bbn_f^{Prei}.$
We now  define $\cc^*(Prei)$ to be the
 $\cz$-submodule of $\cc^*_{\cz}$ generated by
 $$\{\lr{M(\bb)}|\bb\in\bbn_f^{Prei}\}.$$

 \mk\nd{\bf Lemma~6.3} {\sl The $\cz$-submodule $\cc^*(Prei)$ is
 an subalgebra of $\cc^*_{\cz}$ and
 $\{\lr{M(\bb)}|\bb\in\bbn_f^{Prei}\}$ is a $\cz$-basis of
 $\cc^*(Prei).$}

 \nd{\bf Proof.}  If $\bb, \bb_1,
 \bb_2\in\bbn_f^{Prei},$ then the Hall polynomial $g_{M(\bb_1)
 M(\bb_ 2)}^{M(\bb)}$ always exists (see Ringel [R5]). Then it is
 easy to see that $\cc^*(Prei)$ is closed under the multiplication
 $\ast.$  $\Box$

 With similar definitions for  $Prep$, we have

 \mk\nd{\bf Lemma~6.4} {\sl For any $\ba\in\bbn_f^{Prep},
 M(\ba)=\oplus_{\az_i\in \Phi^+_{Prep}}\ba(\az_i)M(\az_i),$ then
 $$\lr{M(\ba)}=\lr{\ba(\az_{i_1})M(\az_{i_1})}\ast\cdots\ast
 \lr{\ba(\az_{i_m})M(\az_{i_m})},$$
 where $\{\az_{i_1}\prec \az_{i_{2}}\prec\cdots\prec \az_{i_m}\}
 $ is the support of $\ba$.  $\Box$}

\mk\nd{\bf Lemma~6.5} {\sl Let $\cc^*(Prep)$ be the
 $\cz$-submodule of $\cc^*_{\cz}$ generated by
 $$\{\lr{M(\ba)}|\ba\in\bbn_f^{Prep}\}.$$ Then
$\cc^*(Prep)$ is
 an subalgebra of $\cc^*_{\cz}$ and
 $\{\lr{M(\ba)}|\ba\in\bbn_f^{Prep}\}$ is a $\cz$-basis of
 $\cc^*(Prep).$ }$\Box$

 \mk\nd{\bf 6.4} Since $Q$ is a tame quiver without oriented cycles, we
can order  a complete set  $\{S_1, S_2,\cdots, S_n\}$
 of non-isomorphic nilpotent simple modules of $\mod\llz$ such that
 $$\ext^1(S_i, S_j)=0\ \text{for}\ i\geq j.$$

 We can now identify $ I=\{1, 2, \dots,n\}$ and $\bbn I=\bbn^n$ such that $S_i$ is the simple module at the vertex $i\in I$. Any module $M$ with dimension vector
 $\bd=(d_1,d_2,\cdots,d_n)$ has a unique filtration
 $$M=M_0\supseteq M_1\supseteq \cdots\supseteq M_n=0$$ with
 factors $M_{i-1}/M_i$ isomorphic to $d_iS_i,$ since
 $\ext(S_i,S_j)=0$ for $i\geq j.$ This shows that the Hall
 polynomial $g^{M}_{d_1S_1\cdots d_nS_n}=1.$
 Then  we have, in $\ch_q$ and  $\ch_q^*$ respectively,
 \begin{eqnarray*} u_{[S_1]}^{d_1}u_{[S_2]}^{d_2}\cdots
 u_{[S_n]}^{d_n}&=&\psi_{d_1}(q)\psi_{d_2}(q)\cdots \psi_{d_n}(q)\sum
 u_{[M(\ba)\oplus M(\bt)\oplus M(\bb)]},\\
u_{[S_1]}^{(\ast d_1)}\ast u_{[S_2]}^{(\ast d_2)}\ast\cdots
 \ast u_{[S_n]}^{(\ast d_n)}&=&v^{-(d_1+d_2+\cdots+d_n)+\lr{\bd,\bd}}\sum
 u_{[M(\ba)\oplus M(\bt)\oplus M(\bb)]},
\end{eqnarray*}
 where the summation is over  the triples $(M(\ba), M(\bt), M(\bb))$
with   $M(\ba)$   preprojective,   $M(\bt)$  regular,  $M(\bb)$  preinjective, and $\udim M(\ba)+\udim M(\bt)+\udim
 M(\bb)=(d_1,\cdots,d_n)=\bd.$

 For any $\ba\in\bbn_f^{Prep},$ let
 $\{\az_{i_1}\prec\az_{i_2}\prec\cdots\prec\az_{i_m}\}$ be the support
of $\ba$ and, for  $1\leq t\leq m,$ define
 \begin{eqnarray*}
\ba_t&=&\ba(\az_{i_t})\az_{i_t}=(a_{1t},a_{2t},\cdots, a_{nt})\in \bbn^n;\\
 \fk{m}_{\ba_t}&=&u_{[S_1]}^{(\ast a_{1t})}\ast
u_{[S_2]}^{(\ast a_{2t})}\ast\cdots\ast u_{[S_n]}^{(\ast a_{nt})};\\
\fk{m}_{\ba}&=&\fk{m}_{\ba_1}\ast \fk{m}_{\ba_2}\ast\cdots\ast
\fk{m}_{\ba_m}. \end{eqnarray*}
 Similarly  for $\bb\in\bbn_f^{Prei}$ define
\begin{eqnarray*} \fk{m}_{\bb_t}&=&u_{[S_1]}^{(\ast b_{1t})}\ast
u_{[S_2]}^{(\ast b_{2t})}\ast\cdots\ast u_{[S_n]}^{(\ast b_{nt})};\\
\fk{m}_{\bb}&=&\fk{m}_{\bb_m}\ast \fk{m}_{\bb_{m-1}}\ast\cdots\ast
\fk{m}_{\bb_1}.
\end{eqnarray*}

\mk\nd{\bf Lemma~6.6} {\sl For any $ \ba \in \bbn_f^{Prep}$ and $ \bb
\in \bbn_f^{Prei}$,  in $\ch^*$, we have
\setcounter{equation}{0}
\begin{equation}\fk{m}_{\ba}=\lr{M(\ba)}+\sum_{\dim\co_{M(\ba')\oplus
M(\bt')\oplus M(\bb')}<\dim\co_{M(\ba)} }c^{\ba}_{\ba' \bt'
\bb'q}u_{[M(\ba')\oplus M(\bt')\oplus M(\bb')]}.
\end{equation}
 Here the sum
ranges over  all triples $M(\ba'), M(\bt'), M(\bb')$
with  $M(\ba')$  preprojective, $M(\bt')$ regular,
$M(\bb')$
 preinjective, and $\udim M(\ba')+\udim M(\bt')+\udim
 M(\bb')=\sum_{\az\in Prep}\ba(\az)\az$, and $c^{\ba}_{\ba' \bt'
\bb'q}\in\bbz[v,v^{-1}];$
\begin{equation} \fk{m}_{\bb}=\lr{M(\bb)}+\sum_{\dim\co_{M(\ba'')\oplus
M(\bt'')\oplus M(\bb'')}<\dim\co_{M(\bb)} }d^{\bb}_{\ba'' \bt''
\bb''q}u_{[M(\ba'')\oplus M(\bt'')\oplus M(\bb'')]},
\end{equation}
  where the
sum is over all triples $M(\ba''), M(\bt''), M(\bb'')$ with
$M(\ba'')$  preprojective, $M(\bt'')$ regular, $M(\bb'')$
 preinjective, and $\udim M(\ba'')+\udim M(\bt'')+\udim
 M(\bb'')=\sum_{\bz\in Prep}\bb(\bz)\bz$, and $d^{\bb}_{\ba'' \bt''
\bb''q}\in\bbz[v,v^{-1}].$}

\nd{\bf Proof.} (1) Since $M(\az_{i_t})$ is exceptional, then
by Lemma~2.1, $\co_{\ba(\az_{i_t})M(\az_{i_t})}$ is a unique orbit
of maximal dimension in $\bbe_{\ba(\az_{i_t})\az_{i_t}}.$  Note that all simple modules are exceptional. We have
\begin{eqnarray*}\fk{m}_{\ba_t}&=&u_{[S_1]}^{(\ast a_{1t})}\ast
u_{[S_2]}^{(\ast a_{2t})}\ast\cdots\ast u_{[S_n]}^{(\ast a_{nt})}\\
&=&\lr{a_{1t}S_1}\ast\lr{a_{2t}S_2}\ast \cdots\ast
\lr{a_{nt}S_n}\\
&=&v^{-\dim(\ba(\az_{i_t})M(\az_{i_t}))+\dim\ed(\ba(\az_{i_t})M(\az_{i_t}))}\sum_{\udim
M=\ba(\az_{i_t})\az_{i_t}}u_{[M]}\\
&=&\lr{\ba(\az_{i_t})M(\az_{i_t})}+\sum_{\dim\co_M<\dim\co_{\ba(\az_{i_t})M(\az_{i_t})}}
v^{-\dim\ext(M,M)}\lr{M}.
\end{eqnarray*}
 Because
$\ext(M(\az_{i_t}),M(\az_{i_s}))=0$ and
$\hom(M(\az_{i_s},M(\az_{i_t}))=0$ for $i_t\leq i_s,$ by Lemma
2.3(7)  and Lemma 6.4, we have
\begin{eqnarray*}\fk{m}_{\ba}&=&\fk{m}_{\ba_1}\ast \fk{m}_{\ba_2}\ast\cdots\ast
\fk{m}_{\ba_m}\\
&=&\lr{M(\ba)}+\sum_{\dim\co_{M(\ba')\oplus
M(\bt')\oplus M(\bb')}<\dim\co_{M(\ba)}}c^{\ba}_{\ba' \bt'
\bb'q}u_{[M(\ba')\oplus M(\bt')\oplus M(\bb')]},
\end{eqnarray*}
 which satisfies
the condition. The proof for (2) is dual, so the proof is
completed.  $\Box$

\mk\nd{\bf Remark.} {In Lemma~6.6, the element $v$ is equal to
$\sqrt{q},$ but the degree of $v^{-1}$  in $c^{\ba}_{\ba' \bt'
\bb'q}$ or in $d^{\bb}_{\ba'' \bt'' \bb''q}$ is bounded and
independent of $\fq.$} (See Lemma~1.2.)

\bigskip
\centerline{\bf 7. Integral bases for the generic composition
algebras}
\bigskip

\nd{\bf 7.1} In this section, we still assume that $ Q$ is connected tame quiver without oriented cycles. We first consider the  embedding of the
representation category of the Kronecker quiver into the
representation category of $Q$.

Let  $e$ be an extending vertex of $Q$ and
$\llz=\fq Q:$ the path algebra of $Q$ over $\fq.$ Let $P=P(e)$ be
projective module cover of the simple module $S_e$. Set $\fk{p}=\udim P(e).$ Clearly
$\lr{\fk{p},\fk{p}}=1=\lr{\fk{p},\dz}$ and there exists a unique
indecomposable preprojective module $L$ with $\udim L=\fk{p}+\dz.$
Moreover we have $\hom_{\llz}(L,P)=0$ and $\ext_{\llz}(L,P)=0.$
This means that $(P,L)$ is an exceptional pair. Let $\fk{C}(P,L)$
be the smallest full subcategory of $\mod\llz$ which contains $P$
and $L$ and is closed under extensions, kernels of epimorphisms
and cokernels of monomorphisms. Also we have
$\dim_{\fq}\hom_{\llz}(P,L)=2,$ therefore $\fk{C}(P,L)$ is
equivalent to the module category of the Kronecker quiver over
$\fq.$ Thus it induces an exact embedding $F:\mod
K\hookrightarrow\mod\llz,$ where $K$ is the path algebra of the
Kronecker quiver over $\fq.$ We note here that the embedding
functor $F$ is essentially independent of the field $\fq.$ This
gives rise to an injective homomorphism of algebras, still denoted by
$F:\ch^*(K)\hookrightarrow\ch^*(\llz).$ In $\ch^*(K)$ we have
defined the element $E_{m\dz_K}$ for $m\geq 1.$  Denote
 $E_{m\dz}=F(E_{m\dz_K}).$ Since $E_{m\dz_K}\in
\cc^*(K),$ and $\lr{L}, \lr{P}\in \cc^*(\llz)$, so $E_{m\dz}$ is in $\cc^*(\llz)$ and even  in
$\cc^*(\llz)_{\cz}.$ Let $\cal K$ be the subalgebra of
$\cc^*(\llz)$ generalized by $E_{m\dz}$ for $m\in\bbn,$ it is a
polynomial ring on infinitely many variables
$\{E_{m\dz}|m\geq1\},$ and its integral form is the polynomial
ring on variables $\{E_{m\dz}|m\geq1\}$ over $\cz.$

\mk\nd{\bf 7.2}
 We may list all non-homogeneous tubes $\ct_1,\ct_2,\cdots,\ct_s$
 in $\mod\llz$ (in fact, $s\leq 3$). For each $\ct_i$, let
 $r_i=r(\ct_i)$ be the period of $\ct_i,$ i.e., the number of
 quasi-simple modules in $\ct_i$. Then $r_i>1.$
 It is well-known that ( for example see [CB])

 \mk\nd{\bf Lemma~7.1} {\sl We have the equation
 $\sum_{i=1}^{s}(r_i-1)=|I|-2$  and the multiplicity of each imaginary root
 $m\dz $ is equal to $|I|-1,$ where $|I|$ is the number of
 vertices of $Q.$} $\Box$

\mk\nd{\bf 7.3}  For each non-homogeneous tube $\ct_i,$ as we did in
Section~5, we have the generic composition algebra $\cc^*(\ct_i)$ of
$\ct_i$ and its integral form $\cc^*(\ct_i)_{\cz}.$ For each
$\ct_i$ we have the set $\Pi^a_i$ of aperiodic $r_i$-tuples of
partitions such that for any $\pi_i\in \Pi^a_i,$ $M_i(\pi_i)$ is an
aperiodic module in $\ct_i.$ We have constructed in~5.3  the element
$$E_{\pi_i}=\lr{M_i(\pi_i)}+\sum_{\lz_i\in\Pi_{i}\setminus\Pi^a_{i},
\lz_i\prec \pi_i}\eta^{\pi_i}_{\lz_i}\lr{M_i(\lz_i)}.$$
Then $\{E_{\pi_i}|\pi_i\in\Pi^a_i\}$ is a $\cz$-basis of
$\cc^*(\ct_i)_{\cz}.$

Let $\cm$ be the set of quadruples $\bc=(\ba_\bc, \bb_\bc, \pi_{\bc}, w_{\bc})$ such that $\ba_{\bc}\in \bbn_f^{Prep}$, $\bb_{\bc}\in \bbn_f^{Prei}$, $\pi_{\bc}=(\pi_{1\bc}, \dots, \pi_{s\bc})\in \Pi_{1}^a\times \cdots \times\Pi_{s}^{a}$, and $w_{\bc}=(w_1\geq w_2\geq\cdots\geq w_t)$ is a partition.

Then for each $\bc\in\cm$ we define
$$E^{\bc}=\lr{M(\ba_{\bc})}\ast E_{\pi_{1 \bc}}\ast  E_{\pi_{2
\bc}}\ast\cdots\ast  E_{\pi_{s \bc}}\ast
E_{w_{\bc}\dz}\ast\lr{M(\bb_{\bc})},$$ where $\lr{M(\ba_{\bc})}$
and $\lr{M(\bb_{\bc})}$ are defined in~6.3, $E_{\pi_{i
\bc}}$ is defined above and
$E_{w_{\bc}\dz}$ is defined in~3.5. Obviously,
$\{E^{\bc}|\bc\in\cm\}$ lies in $\cc^*(\llz),$ in fact in
$\cc^*(\llz)_{\cz},$ and are linearly independent over $\bbq(v).$

\mk\nd{\bf Proposition~7.2} {\sl The set $\{E^{\bc}|\bc\in\cm\}$
is a $\bbq(t)$-basis of $\cc^*(\llz).$}

 The  proof of Proposition~7.2 will be given in~7.4.  We first  need some
preparation.

 \mk\nd{\bf Lemma~7.3} {\sl Let $\{S_j\mid 1\leq j\leq
r_i\}$ be a complete set of non-isomorphic quasi-simple modules of
a non-homogeneous tube $\ct_i$ such that $S_j=\tau^{(j-1)}S_1$ and
let $\ch^*(\ct_i)$ be the generic integral form of the
 twisted Ringel-Hall algebra of $\ct_i$ over $\cz=\bbz[t,t^{-1}].$
For any $ l \in \bbn $ and $ 1\leq j \leq r_i$, let $ \pi, \pi' \in \Pi^a_i$ such that $ S_j[l]=M(\pi) $ and $S_{j+1}[l]=M(\pi')$. Then
\setcounter{equation}{0}
\begin{eqnarray}
\label{7.3(1)}  u_{[S_j[l]]} &\equiv&
\sum_{\lambda\preceq\pi,\lambda\in\Pi_i^a}a_\lambda E_\lambda
\quad  \pmod{(t-1)\ch^*(\ct_i)} \quad \text{ if  $r_i\nmid l$},\\
\label{7.3(2)}
u_{[S_j[l]]}-u_{[S_{j+1}[l]]}&\equiv&
\sum_{\lambda\preceq\pi(\text{or}~\pi'),\lambda\in\Pi^a_i}a_\lambda
E_\lambda \quad \pmod{(t-1)\ch^*(\ct_i)} \quad  \text{if $r_i\mid l$}.
\end{eqnarray}
Here $ a_{\lambda}\in \bbq$.}

\nd{\bf Proof.} Without loss of generality, we may take  $j=1.$
   When $l=1$, we have $u_{[S_1]}=E_{S_1}.$ The conclusion follows.
   We suppose that the conclusion is true  when $1\leqslant l\leqslant r_i-2$.  Then the assumption  $$u_{[S_1[l]]} \equiv
\sum_{\lambda\preceq\pi_1,\lambda\in\Pi^a_i}a_\lambda E_\lambda\quad  \pmod{(t-1)\ch^*(\ct_i)}$$
and $u_{[S_1[l]]}u_{[S_{l+1}]}-u_{[S_{l+1}]}u_{[S_1[l]]}\equiv
u_{[S_1[l+1]]} \quad \pmod{(t-1)\ch^*(\ct_i)} $ imply
\begin{align*}
u_{[S_1[l+1]]}\equiv &
(\sum_{\lambda\preceq\pi_1,\lambda\in\Pi^a_i}a_\lambda
E_{\lambda}) E_{S_{l+1}} -E_{S_{l+1}}
(\sum_{\lambda\preceq\pi_1,\lambda\in\Pi^a_i}a_\lambda E_{\lambda}
) \pmod{(t-1)\ch^*(\ct_i)} \\
\equiv &(\sum_{\lambda\preceq\pi_1,\lambda\in\Pi^a_i}a_\lambda E_\lambda)
E_{S_{l+1}} -E_{S_{l+1}}
(\sum_{\lambda\preceq\pi_1,\lambda\in\Pi^a_i}a_\lambda
E_\lambda)=\sum_{\lambda\preceq\pi,\lambda\in\Pi^a_i}a'_\lambda
E_\lambda\end{align*}
 since  $\{E_\lambda\mid
\lambda\in\Pi^a_i\}$ is a basis of $\ct_i.$ Thus the conclusion is true
 for  $l+1.$
For $l=r_i$, by assumption, we have
$$u_{[S_2[l-1]]} \equiv
\sum_{\lambda\preceq\pi_1,\lambda\in\Pi^a_i}a_\lambda E_\lambda
\quad\pmod{(t-1)\ch^*(\ct_i)}.$$

\begin{align*}
u_{[S_1[l]]}-u_{[S_2[l]]}\equiv &u_{[S_1]}u_{[S_2[l-1]]}-u_{[S_2[l-1]]}u_{[S_1]} \quad\pmod{(t-1)\ch^*(\ct_i)}. \\
\equiv & E_{S_1}(\sum_{\lambda\preceq\pi_1,\lambda\in\Pi^a_i}a_\lambda E_\lambda)
 -(\sum_{\lambda\preceq\pi_1,\lambda\in\Pi^a_i}a_\lambda
E_\lambda)E_{S_1} \quad\pmod{(t-1)\ch^*(\ct_i)}\\
\equiv &
\sum_{\lambda\preceq\pi(\text{or}~\pi')\lambda\in\Pi^a_i}a_\lambda
E_\lambda \quad\pmod{(t-1)\ch^*(\ct_i)}.
\end{align*}
    Now we consider the general case. Let $l=kr_i+m,0<m\leq
    r_i-1,$ if $m=1$, by assumption, we have $$u_{[S_1[l-1]]}-u_{[S_2[l-1]]}\equiv
\sum_{\lambda\preceq\pi_1(\text{or}~\pi_1')\lambda\in\Pi^a_i}a_\lambda
E_\lambda\ \ \ \ (\text{mod}(t-1)\ch^*(\ct)).$$
Hence
\begin{align*}
u_{[S_1[l]]}\equiv &(u_{[S_1[l-1]]}-u_{[S_2[l-1]]})u_{[S_1]}-u_{[S_1]}(u_{[S_1[l-1]]}-u_{[S_2[l-1]]})\\
\equiv &(\sum_{\lambda\preceq\pi_1(\text{or}~\pi_1')\lambda\in\Pi^a_i}a_\lambda
E_\lambda)E_{S_1}
 -E_{S_1}(\sum_{\lambda\preceq\pi_1(\text{or}~\pi_1')\lambda\in\Pi^a_i}a_\lambda
E_\lambda)\\
\equiv &
\sum_{\lambda\preceq\pi(\text{or}~\pi')\lambda\in\Pi^a_i}a'_\lambda
E_\lambda \quad\pmod{(t-1)\ch^*(\ct_i)}.
\end{align*}
If  $~2\leq m\leq r_i-1,$ by assumption,
$$u_{[S_1[l-1]]}\equiv\sum_{\lambda\preceq\pi_1,\lambda\in\Pi^a_i}a_\lambda
E_\lambda \quad\pmod{(t-1)\ch^*(\ct_i)}. $$
Hence
\begin{align*} u_{[S_1[l]]}\equiv&
u_{[S_1[l-1]]}u_{[S_l]}-u_{[S_l]}u_{[S_1[l-1]]}\\
\equiv&(\sum_{\lambda\preceq\pi_1,\lambda\in\Pi^a_i}a_\lambda
E_\lambda)E_{S_l}-E_{S_l}(\sum_{\lambda\preceq\pi_1,\lambda\in\Pi^a_i}a_\lambda
E_\lambda)\\
\equiv & \sum_{\lambda\preceq\pi,\lambda\in\Pi^a_i}a_\lambda
E_\lambda \quad\pmod{(t-1)\ch^*(\ct_i)}.
\end{align*}
Then the conclusion is true. When $r_i\mid l$, it can be proved by
a similar
 method for $ l=r_i.$  $\Box$

\mk\nd{\bf Remark.} Of course we can replace $\ch^*(\ct)$ in Lemma
7.3 by $\cc^*(\ct)$ using the natural embedding
$\cc^*(\ct)/(t-1)\cc^*(\ct)$ into $\ch^*(\ct)/(t-1)\ch^*(\ct),$
here we consider the integral forms over $\cz.$

\mk\nd{\bf Lemma~7.4} {\sl In $\cc^*(\llz)_{\cz},$

$$\te_{n\dz}=\sum_{\substack{m_1\leq\cdots\leq
m_s\\ m_1+\cdots+m_s=n}}b_{m_1,\cdots,m_s}E_{m_1\dz}\ast\cdots\ast
E_{m_s\dz},$$ where $b_{m_1,\cdots,m_s}\in\cz.$}

\nd{\bf Proof.}  By the relation $$E_{0\dz}=1,\ \
E_{k\dz}=\frac{1}{[k]}\sum_{s=1}^k v^{s-k}\te_{s\dz}\ast
E_{(k-s)\dz},$$ we can solve the equation inductively to get the
relation in the lemma. $\Box$

\mk It is known from Ringel that the Lie subalgebra
$\fn^+\subseteq\cc^*(\llz)_{\cz}/(t-1)\cc^*(\llz)_{\cz}$ generated
by $u_{[S_i]}$ ($i\in I$) over $\bbq$ is the positive part of the
corresponding affine Kac-Moody Lie algebra over $\bbq,$ and
$\cc^*(\llz)_{\cz}/(t-1)\cc^*(\llz)_{\cz}$
 is the universal enveloping algebra of $\fn^+.$

 For each non-homogeneous tube $\ct_i$ of rank $r_i,$ we denote
 $u_{\az,i}=u_{[S_j[l]]}$ where $S_j[l]$ is indecomposable in
 $\ct_i$ and $\udim{S_j[l]}=\az$ is a real root; and
 $u_{j,m\dz,i}-u_{j+1,m\dz,i}=u_{[S_j[l]]}-u_{[S_{j+1}[l]]}$ where $S_j[l]$ is indecomposable in
 $\ct_i$ and $\udim{S_j[l]}=m\dz$  an imaginary root. Let
 $\Psi:\cc^*(\llz)_{\cz}\ra\cc^*(\llz)/(t-1)\cc^*(\llz)_{\cz}$
 be the canonical projection. Then one of the main results in [FMV]
 is the following  of which the proof depends on Lemma~7.1.

 \mk\nd{\bf Proposition~7.5} {\sl The vectors $\Psi(u_{[M(\az)]})$ for
 $\az\in \Phi^+_{Prep};$ $\Psi(u_{\az,i})$ for $\az\in\ct_i$ real root,
 $i=1,\cdots,s;$ $\Psi(u_{j,m\dz,i}-u_{j+1,m\dz,i}),$ $m\geq 1,$
 $1\leq j\leq r_i,$ $i=1,\cdots,s;$ $\Psi(\te_{n\dz}), n\geq 1$
 and $\Psi(u_{[M(\bz)]})$ for $\bz\in \Phi^+_{Prei}$ form a $\bbz$-basis of
 $\fn^+.$
}

Note that it is easy to see that all vectors in Proposition~7.5
belong to the Lie algebra $\fn^+,$ and they are linearly
independent over $\bbq.$ For example, $\Psi(\te_{n\dz}), n\geq 1,$
lie in $\fn^+$ by Lemma~3.9. Then by Lemma~7.1, one can prove that
those vectors give rise to a $\bbz$-basis of $\fn^+.$

\mk\nd{\bf 7.4} 
{\bf Proof of Proposition~7.2.} By the definition of
 $\{E^{\bc}|\bc\in\cm\},$ we see that they are linearly
 independent over $\bbq(t).$ For any weight ( or, dimension vector
 ) $w\in\bbn I,$ we define the $\bbq(t)$-space $V_w$ to be spanned
 by those $E^{\bc},$ $\bc\in\cm,$ such that
 $E^{\bc}\in\cc^*(\llz)_w.$ It is well-known from Lusztig that
 $$\dim_{\bbq(t)}\cc^*(\llz)_w=\dim_{\bbq}(\cc^*(\llz)_{\cz}/(t-1)\cc^*(\llz)_{\cz})_w$$
 and the monomials in a fixed order on the basis elements of $\fn^+$  in Proposition~7.5 form   a PBW basis of
$\cc^*(\llz)_{\cz}/(t-1)\cc^*(\llz)_{\cz}$ over $\bbq.$ However,
Lemma~7.3 and~7.4 implies that those PBW basis elements can be obtained by
applying $\Psi$ on $\{E^{\bc}|\bc\in\cm\}.$ Therefore
$\dim_{\bbq(t)}V_w\geq\dim_{\bbq(t)}\cc^*(\llz)_w$ for any
$w\in\bbn I.$ Hence  $\{E^{\bc}|\bc\in\cm\}$ is a
$\bbq(t)$-basis of $\cc^*(\llz).$   $\Box$

As a consequence, the canonical mapping
$$\varphi:\cc^*(Prep)\otimes_{\mathbb{Q}(t)}
\cc^*(\mathcal{T}_1)\otimes_ {\mathbb{Q}(t)}\cdots
\otimes_{\mathbb{Q}(t)} \cc^*(\mathcal{T}_s)\otimes_{\mathbb{Q}(t)}
 \mathcal{K}\otimes_{\mathbb{Q}(t)} \cc^*(Prei)\rightarrow \cc^*(\Lambda)$$
  is an isomorphism
  of $\mathbb{Q}(t)$-spaces.

\mk\nd{\bf 7.5} We may consider the ring $\ca=\bbq[t,t^{-1}],$ and
$\cc^*(\llz)_{\ca}$ is the $ \ca$-subalgebra of the generic composition
algebra $\cc^*(\llz)$  generated by $u_i^{(\ast
m)}=\frac{u_i^{\ast m}}{[m]!},  (i\in I)$.

\mk\nd{\bf Proposition~7.6} {\sl The set $\{E^{\bc}|\bc\in\cm\}$
is an $\ca$-basis of $\cc^*(\llz)_{\ca}.$}

\nd{\bf Proof.} For any monomial $\fk{m}$ on the divided powers of
$u_{[S_i]}$ ($\ i\in I$) by Proposition~7.2,
$$\fk{m}=\sum_{\bc\in\cm}f_{\fk{m}, \bc}(t)E^{\bc}\ \ \ \
(\text {finite sum})$$ in $\cc^*(\llz),$ where $f_{\fk{m},
\bc}(t)\in\bbq(t)$ and $v$ is an indeterminate.
 Note that $E_{\pi_{i \bc}}$ in the definition of $E^\bc$ has the
form (cf.~5.3)
$$E_{\pi}=\lr{M(\pi)}+\sum_{\lz\in\Pi_{\az}\setminus\Pi^a_{\az},
\lz\prec \pi}\eta^{\pi}_{\lz}\lr{M(\lz)}$$
with $\eta^{\pi}_{\lz}\in\cz.$ The formula
$\fk{m}=\sum_{\bc\in\cm}f_{\fk{m}, \bc}(v)E^{\bc}$ still holds in
$\ch^*$ for taking $v=\sqrt{q}.$ Thus, by Lemma~1.2,  for
each $\bc\in\cm,$ there exists $N(\bc)\in\bbn$ such that
$(\sqrt{q})^{N(\bc)}f_{\fk{m},\bc}(\sqrt{q})\in\mathbb{Z}$ for all
$q=p^l$ with $p$ a prime number and $l\geq 1$ in $\bbn.$ It is
easily seen that $t^{N(\bc)}f_{\fk{m},\bc}(t)$ is a polynomial in
$\bbq[t].$ Therefore $f_{\fk{m},\bc}(v)\in\bbq[t,t^{-1}].$   $\Box$

\mk\nd {\bf Corollary 7.7} {\sl The multiplication map
$$\varphi:\cc^*(Prep)_{\ca}\otimes_{\ca}\cc^*(\mathcal{T}_1)_{\ca}\otimes_ {\ca}\cdots \otimes_{\ca}
\cc^*(\mathcal{T}_s)_{\ca}\otimes_{\ca}
 \mathcal{K}_{\ca}\otimes_{\ca} \cc^*(Prei)_{\ca}\rightarrow \cc^*(\Lambda)_{\ca}$$
  is an isomorphism  of $\ca$-modules.} $\Box$

\bigskip
\centerline{\bf 8. A bar-invariant basis of  $\cc^*(\llz)_{\ca}$}
\bigskip

To simplify the notations, in the next rest of the paper, we will use $v$ for the indeterminate $t$. However, we will perform computations over $\bbf_q$ with $v=q^{1/2}$. It should be clear from the context that formulae are independent of $q$ and one can obtain the same formulation as in the generic case as discussed in 1.4.

\nd{\bf 8.1} The first  part of this section is devoted to
finding a monomial basis and a triangular relation with the basis
$\{E^{\bc}|\bc\in\cm\}.$

 We first define the
variety\footnote{Note that the definition of $\co_{\bc}$ here is
different with that in~4.2}

  $$\co_{\bc}=\co_{M(\ba_{\bc})}\star\co_ {M_{\pi_{1 \bc}}}\star\co_  {M_{\pi_{2
\bc}}}\star\cdots\star\co_{M_{\pi_{s
\bc}}}\star\cn_{w_{\bc}\dz}\star\co_{M(\bb_{\bc})}$$ for any
$\bc\in\cm,$ where
$\cn_{w_{\bc}\dz}=\cn_{w_1\dz}\star\cdots\star\cn_{w_t\dz}$ if
$w_{\bc}=(w_1, w_2,\cdots,w_t)$ and $\cn_{w_i\dz}$ are the union
of orbits of regular modules of $\fk{C}(P,L)$ with dimension
vector $w_i\dz.$

   Then  by Proposition~7.6, Lemma~6.6 can be rewritten as follows:

\mk\nd{\bf Lemma~8.1} {\sl For  any $\ba\in\bbn_f^{Prep}$ and $\bb\in\bbn_f^{Prei}$, in  $\cc^*(\llz)$ we have
\setcounter{equation}{0}
\begin{align} \fk{m}_{\ba}=&\lr{M(\ba)}+\sum_{\dim\co_{\bc}<\dim\co_{\ba} }f^{\ba}_{
\bc}E^{\bc},\\
\fk{m}_{\bb}=&\lr{M(\bb)}+\sum_{\dim\co_{\bc}<\dim\co_{\bb} }g^{\bb}_{\bc}E^{\bc},
\end{align}
  where  $f^{\ba}_{\bc}, g^{\bb}_{\bc}\in\bbq[v,v^{-1}]$ and $ \bc \in \cm$.}$\Box$

\mk\nd{\bf Remark.} The conclusion is also true in Lemma~8.1 if we
take $M(\ba)$ to be  finitely many copies of a exceptional module.

\mk\nd{\bf Lemma~8.2} {\sl Let $\pi\in\Pi^a_i$ for some $\ct_i,$
then there exists a monomial $\mathbf{m}_{\pi}$  on the divided
powers of $u_{[S_i]}$ ($ i\in I$) such that

$$\mathbf{m}_{\pi}=E_{\pi}+\sum_{\dim\co_{\bc}<\dim\co_{\pi} }f^{\pi}_{
\bc}E^{\bc},$$ where $f^{\pi}_{ \bc}\in\bbq[v,v^{-1}].$ }

\nd{\bf Proof.} We set $\{\theta_1,\theta_2,\cdots,\theta_{r_i}\}$
 to be a complete set of non-isomorphic quasi-simple modules of
$\ct_i$ in the natural order (see Section~5). By Proposition
5.2, we then have
$${\fk{m}}^{(w_{\pi})}=E_{\pi}+\sum_{\lz\in\Pi^a_{\az}, \lz\prec\pi}
v^{-\dim M(\pi)+\dim\ed M(\pi)+\dim M(\lz)-\dim\ed
M(\lz)}g^{\lz}_{w_{\pi}}(v^2)E_{\lz},$$ where
${\fk{m}}^{(w_{\pi})}=\theta_{j_1}^{(\ast  e_1)}\ast\cdots\ast
\theta_{j_t}^{(\ast  e_t)}.$    Since  each $\theta_{j_i}$ is an exceptional module, we have $\lr{u_{[\theta_{j_p}]}}^{(\ast e_p)}=\lr{e_p\theta_{j_p}}$ (see the proof in~6.2).

   Let $\pi_{j_p}\in\Pi^a_i$ such that
   $M(\pi_{j_p})=e_p\theta_{j_p}$ and $\udim M(\pi_{j_p})=(d_1,\cdots,d_n)$ with $I$ ordered as in~6.4.
By Lemma~8.1 and its remark, we define a monomial
   $\fk{m}_{j_p}$ such that
   $$\fk{m}_{j_p}=\lr{S_1}^{(\ast
   d_1)}\ast\cdots\ast\lr{S_n}^{(\ast d_n)}
   =\lr{M(\pi_{j_p})}+\sum_{\dim\co_{\bc}<\dim\co_{M(\pi_{j_p})}}f_{\bc}^{\pi_{j_p}}E^{\bc}$$
where $f^{\pi_{j_p}}_{\bc}\in\bbq[v,v^{-1}].$

Let $L_0=e_{1}\theta_{j_1},L_1=e_{1}\theta_{j_1}\diamond
e_{2}\theta_{j_2},L_2=L_1\diamond
e_{3}\theta_{j_3},\cdots,L_{t-1}=L_{t-2}\diamond
e_{t}\theta_{j_t}. $ By Lemma~2.3(6), We have $M({\pi})\simeq
L_{t-1}$.  Similar to the argument as in~5.3, we have $g^{L_p}_{L_{p-1},e_{j_{p+1}}\theta_{j_{p+1}}}=1,$ for $1\leq
p\leq t-2.$ Define $\alpha_p=\udim{ L_{p-1}}$ and $\beta_p=\udim
{M(\pi_{j_p})}.$  By Lemma~2.3(6), we have
$$\dim \co_{L_p}=\dim\co_{L_{p-1}}+\dim\co_{e_{p+1}\theta_{j_{p+1}}}+\wm(\alpha_p,\beta_p)$$
 or
 $$\cdim \co_{L_p}=\cdim\co_{L_{p-1}}+\cdim\co_{e_{p+1}\theta_{j_{p+1}}}-\lr{\beta_p,\alpha_p}. $$
 Thus
  $$\dim \co_{M(\pi)}=\dim\co_{L_{t-1}}=\sum^t_{p=1}\dim\co_{e_p\theta_{j_p}}+\sum^{t-1}_{p=1}\wm(\alpha_p,\beta_p).$$
For any $\bc\in\cm$ with
  $\dim\co_{\bc}<\dim\co_ {e_{p+1}\theta_{j_{p+1}}},$ by Lemma~2.2, we have
$$ \aligned \cdim(\overline{\co}_{L_{p-1}} \star
 \overline{\co_{\bc}})
  &=\cdim(\overline{\co}_{L_{p-1}})+\cdim(\overline{\co_{\bc}})-\lr{\beta_p,\alpha_p}+r
\\& > \cdim(\overline{\co}_{L_{p-1}})+\cdim(\co_{
e_{p+1}\theta_{j_{p+1}}})-\lr{\beta_p,\alpha_p}
\\&=\cdim(\co_{L_p}),
\endaligned $$
Thereby, if we take
$\fk{m}_{\pi}=\fk{m}_{\pi_{j_1}}\ast\cdots\ast\fk{m}_{\pi_{j_t}}$,
then
$$ \aligned
\fk{m}_{\pi}&=(\lr{\theta_{j_1}}^{(\ast
e_1)}+\sum_{\dim\co_{\bc_1}<\dim\co_{{e_1}\theta_{j_1}} }
f^{\pi_{j_1}}_{\bc_1}E^{\bc_1})
\ast\cdots\ast(\lr{\theta_{j_t}}^{(\ast
e_t)}+\sum_{\dim\co_{\bc_t}<\dim\co_{e_t\theta_{j_t}} }
f^{\pi_{j_t}}_{ \bc_t}E^{\bc_t})\\
&=(\lr{M(\pi_{j_1})}+\sum_{\dim\co_{\bc_1}<\dim\co_{{e_1}\theta_{j_1}}
} f^{\pi_{j_1}}_{\bc_1}E^{\bc_1})
\ast\cdots\ast(\lr{M(\pi_{j_t})}+\sum_{\dim\co_{\bc_t}<\dim\co_{e_t\theta_{j_t}}
} f^{\pi_{j_t}}_{ \bc_t}E^{\bc_t})
\\&=E_{\pi}+\sum_{\dim\co_{\bc}<\dim\co_{\pi} }f^{\pi}_{
\bc}E^{\bc},
\endaligned $$
where $f^{\pi}_{ \bc}\in\bbq[v,v^{-1}].$ The proof is finished.
$\Box$

\mk\nd{\bf Lemma~8.3} {\sl Let $E_{n\delta}$ be the image embedded
in $\cc^*(\Lambda)$ of the element  $E_{n\delta}$ in $\ck$,then
there exists a monomial $\fk{m}_{n\delta}$  on the divided powers
of $u_{[S_i]}$ ($\ i\in I$)  such that
$$\fk{m}_{n\delta}=E_{n\delta}+\sum_{\dim\co_{\bc}<\dim\co_{n\delta} }h^{n\delta}_{\bc}E^{\bc},$$
where $h^{n\delta}_{ \bc}\in\bbq[v,v^{-1}].$ }

\nd{\bf Proof.} We let $\theta_1,\theta_2$ be the two simple
objects of $\fk{C}(P,L).$ By Lemma~3.13(3), we then have
$$\lr{\theta_2}^{(\ast n)}\ast \lr{\theta_1}^{(\ast n)}=E_{n\delta}+
\sum_{\dim\co_{\bc}<\dim\co_{n\delta} }f^{n\delta}_{\bc}E^{\bc},\
\ \ \ \text{with}\ \ f^{n\delta}_{\bc}\in\bbq[v,v^{-1}].$$

Suppose that $\udim n\theta_1={\bf d}'=(d'_1,\cdots,d'_n)$ and
$\udim n\theta_2={\bf d}''=(d''_1,\cdots,d''_n)$ in $\bbz I.$
Since $\theta_1,\theta_2$ are the exceptional modules, by  the
remark of Lemma~8.1, we then have

$$\aligned
\fk{m}_1&=\lr{S_1}^{(\ast d'_1)}\ast\lr{S_2}^{(\ast
d'_2)}\ast\cdots\ast\lr{S_n}^{(\ast d'_n)}
\\&=\lr{\theta_1}^{(\ast n)}+\sum_{\dim\co_{\bc}<\dim\co_{
n\theta_1 }}f^{n\theta_1}_{\bc}E^{\bc};\endaligned$$ and
$$\aligned
\fk{m}_2&=\lr{S_1}^{(\ast d''_1)}\ast\lr{S_2}^{(\ast
d''_2)}\ast\cdots\ast\lr{S_n}^{(\ast d''_n)}
\\&=\lr{\theta_2}^{(\ast n)}+\sum_{\dim\co_{\bc}<\dim\co_{
n\theta_2}}g^{n\theta_2}_{\bc}E^{\bc};\endaligned,$$ where
$f^{n\theta_1}_{\bc},g^{n\theta_2}_{\bc}\in\bbq[v,v^{-1}].$ By
representations of the Kronecker quiver, we know that
$\cn_{n\delta}$ is open in $\co_{n\theta_2}\star\co_{n\theta_1}.$
Moreover, $\cn_{n\delta}$ is open,
  then dense in $\overline{\co}_{n\theta_2}\star\overline{\co}_{n\theta_1},$
 that is, $\cn_{n\delta}$ is of maximum dimension,  $G$-stable,
 irreducible and open subvariety of
 $\overline{\co}_{n\theta_2}\star\overline{\co}_{n\theta_1}.$
Since
   $\hom(\overline{\co}_{n\theta_1},\overline{\co}_{n\theta_2})=0$,
we then obtain
   $$\cdim \overline{\co}_{n\theta_2}\star\overline{\co}_{n\theta_1}
   =\cdim \overline{\co}_{n\theta_2}+\cdim
   \overline{\co}_{n\theta_1}-\lr{{\bf d}',{\bf d}''}$$ by Lemma~2.2.
If either
$\co_{\bc}\subset\overline{\co}_{n\theta_2}\setminus\co_{n\theta_2}$
or
   $\co_{\bc'}\subset\overline{\co}_{n\theta_1}\setminus\co_{n\theta_1}$, then
   $$\aligned
   \cdim \overline{\co}_{\bc}\star\overline{\co}_{\bc'}
   &=\cdim \overline{\co}_{\bc}+\cdim \overline{\co}_{\bc'}-\lr{{\bf d}',{\bf d}''}+r
   \\&>\cdim \overline{\co}_{n\theta_2}\star\overline{\co}_{n\theta_1}
   =\cdim \cn_{n\delta}.
   \endaligned$$
We now take $\fk{m}_{n\delta}=\fk{m}_2\ast \fk{m}_1$,then
   $$\aligned
   \fk{m}_{n\delta}&=(\lr{\theta_2}^{(\ast n)}+\sum_{\dim\co_{\bc}<\dim\co_{
n\theta_2 }}g^{n\delta}_{\bc}E^{\bc})\ast(\lr{\theta_1}^{(\ast
n)}+\sum_{\dim\co_{\bc'}<\dim\co_{ n\theta_1
}}f^{n\delta}_{\bc'}E^{\bc'})
\\&=E_{n\delta}+\sum_{\dim\co_{\bc}<\dim\co_{n\delta}
}h^{n\delta}_{\bc}E^{\bc},
\endaligned$$
where $h^{n\delta}_{ \bc}\in\bbq[v,v^{-1}].$  $\Box$

\mk\nd{\bf Proposition 8.4} {\sl For any $E^{\bc}, \bc\in\cm,$
there exists a monomial $\fk{m}_{\bc}$ on the divided powers of
$u_{S_i}, i\in I,$  such that
$$\fk{m}_{\bc}=E^{\bc}+\sum_{\bc'\in \cm,\ dim \co_{\bc'}< dim
\co_{\bc}}h^{\bc}_{\bc'}E^{\bc'}, $$where
$h^{\bc}_{\bc'}\in\bbq[v,v^{-1}].$ }

 \nd{\bf Proof.} According to the structure of the Auslander-Reiten quiver of a tame quiver,
 if $P\in Prep, I\in Prei$ and $R$ is a regular module, we then know
that
  $\co_{P\oplus R\oplus I}$ is open in
 $\overline{\co}_P\star\overline{\co}_R\star\overline{\co}_I$ by
 Lemma~2.3(7). So, we need to prove the same property for
 $E_{\pi}\ast E_{n\delta}$ where $\pi\in\Pi^a_i.$ By Lemma~8.2 and~8.3, there exist  $\fk{m}_\pi$ and  $\fk{m}_{n\delta}$  such that
$$\mathbf{m}_{\pi}=E_{\pi}+\sum_{\dim\co_{\bc}<\dim\co_{\pi} }f^{\pi}_{
\bc}E^{\bc},$$ and
$$\mathbf{m}_{n\delta}=E_{n\delta}+\sum_{\dim\co_{\bc'}<\dim\co_{n\delta} }g^{n\delta}_{\bc'}E^{\bc'},$$
where $f^{n\delta}_{ \bc},g^{n\delta}_{ \bc'}\in\bbq[v,v^{-1}].$

  Since we can find smooth points
  $A\in \overline{\co}_\pi $ and $B\in \overline{\co}_{n\delta}$ such
  that $\hom(B,A)=0$,  we have
  $\hom(\overline{\co}_{n\delta},\overline{\co}_\pi)=0.$ Then,
  $$\cdim \overline{\co}_\pi\star\overline{\co}_{n\delta}=\cdim \overline{\co}_\pi
  +\cdim \overline{\co}_{n\delta}-\lr{n\delta,\alpha}.
  $$
If either $\co_{\bc}\subset\overline{\co}_{\pi}\setminus\co_{\pi}$
or
    $\co_{\bc'}\subset\overline{\co}_{n\delta}\setminus\co_{n\delta}$, we
     have again that
    $$\cdim \overline{\co}_{\bc}\star\overline{\co}_{\bc'}
    >\cdim \overline{\co}_{\pi}\star\overline{\co}_{n\delta}
    =\cdim \co_{\pi}\star\co_{n\delta}.$$
So, we get
    $$\aligned
    \fk{m}_{\bc}&=\fk{m}_\pi\ast\fk{m}_{n\delta}
    \\&=(E_{\pi}+\sum_{\dim\co_{\bc}<\dim\co_{\pi} }f^{\pi}_{\bc}E^{\bc})
    \ast (E_{n\delta}+\sum_{\dim\co_{\bc'}<\dim\co_{n\delta}
}g^{n\delta}_{\bc'}E^{\bc'})
 \\&=E^{\bc}+\sum_{\bc'\in\cm, \dim
\co_{\bc'}< \dim \co_{\bc}}h^{\bc}_{\bc'}E^{\bc'},
\endaligned$$
where $h^{\bc}_{\bc'}\in\bbq[v,v^{-1}].$
$\Box$

\mk\nd{\bf 8.2} Let $\mathcal{A}=\bbq[v,v^{-1}].$ We define the
lattice $\cl'$
 to be the $\bbq[v^{-1}]$-submodule of $\cc^{*}(\Lambda)_{\ca}$ with
the basis $ \{E^{\bc}|\bc\in\cm\}.$  By the argument similar to that
in  Section~4, we can use the standard linear algebra method by
Lusztig to get the relation:
$$\ol{E^{\bc}}=\sum_{\bc'\in\cm}\oz^{\bc}_{\bc'}E^{\bc'}\ \ \ \ \
\text{for any}\ \ \bc\in\cm$$ with $\oz^{\bc}_{\bc'}\in\ca$ such
that $\oz^{\bc}_{\bc}=1$ and if $\oz^{\bc}_{\bc'}\neq 0$ and
$\bc\neq\bc$ then $\dim\co_{\bc'}<\dim\co_{\bc}.$ Thus we can
solve the system of equations
$$\zeta^{\bc}_{\bc'}=\sum_{\dim\co_{\bc'}\leq\dim\co_{\bc''}\leq\dim\co_{\bc}}\oz^{\bc''}_{\bc'}\ol{\zeta^{\bc}_{\bc''}}$$
to get a unique solution such that
$$\zeta^{\bc}_{\bc}=1\ \ \ \ \text{and}\ \ \
\zeta^{\bc}_{\bc'}\in v^{-1}\bbq[v^{-1}]\ \ \ \text{if}\ \
\dim\co_{\bc'}<\dim\co_{\bc}.$$ Let
$$\ce^{\bc}=\sum_{\bc'\in\cm}\zeta^{\bc}_{\bc'}E^{\bc'}\ \ \ \
\text{for any}\ \ \ \bc\in\cm.$$ Note that this is a finite sum.
Then we have the following result.

\mk\nd{\bf Theorem 8.5} {\it The set $\{\ce^{\bc}|\bc\in\cm\}$
provides a  basis of $\cc^*(\llz)_{\ca},$ which is characterized by
the two properties: (a) $\ol{\ce^{\bc}}=\ce^{\bc}$ for all
$\bc\in\cm.$ (b) $\pi(\ce^{\bc})=\pi(E^{\bc}),$ where
$\pi:\cl'\ra\cl'/v^{-1}\cl'$ is the canonical projection. $\Box$}

\bigskip
\centerline{\bf 9. Affine canonical bases}
\bigskip

\nd{\bf 9.1} The Ringel-Hall algebra  $\mh^*(\llz)$ is an
associative $\bbq(v)$-algebra with the
    basis $\{ \lr{M} \;|\; M \in
\cp\}.$ Note that $\lr{M}=v^{-\dim
 M+\dim\ed_{\llz}(M)}u_{[M]}.$ According to \cite{G},  the inner product on
$\mh^*(\llz)$ is given by the formula

    $$(\lr{M},\lr{N})=\delta_{MN}\frac{v^{2 \dim \ed(M) }}{a_M}$$
where $a_M=|\aut(M)|.$ It is known that there exists a polynomial,
still denoted by $a_M,$ such that $a_M(|E|)=|\aut_{\llz^E}M^E|$ for
any finite extension $E$ of $k.$ Following Green\cite{G} and
\cite{R1}, we can define a linear operation
$$r(u_{[L]})=\sum_{[M],[N]}v^{\lr{\udim M,\udim
N}}g_{M,N}^L\frac{a_{M}a_{M}}{a_{L}}u_{M}\otimes u_{N}.$$ We have
the following property:

$$(x,y\ast z)=(r(x),y\otimes z).$$ for any $x, y\ \mbox{and}\ z\in
\mh^*(\llz).$

It is not difficult to see the following:
\bigskip

\nd {\bf Proposition 9.1}  {\sl For any preprojective $\llz$-modules
$P, P'\in\cp$ , regular $\llz$-modules $R, R'\in\cp$  and
preinjective $\llz$-modules $I, I'\in\llz$ , in $\ch^*(\llz)$ we
have

$$(\lr{P}*\lr{R}*\lr{I},\lr{P'}*\lr{R'}*\lr{I'})=\delta_{PP'}\delta_{RR'}\delta_{II'}
    \frac{v^{2 (\dim \ed(P)+ \dim \ed(R)+\dim \ed(I))}}{{a_P}a_R{a_I}}.$$}

This inner product is also well-defined on $\cc^*(\llz),$ which
coincides with the paring defined by Lusztig in \cite{L5}.

\bigskip

\nd{\bf 9.2 }\ \ Let $\Lambda=kQ,$  $k$ a finite field. Consider
regular $\Lambda$-modules $M,\cdots, M_t$ and $L,$ it is well known
that there exists the Hall polynomial $\varphi_{M_1\cdots
M_t}^{L}\in \mathbb{Z}[T]$ such that $\varphi_{M_1\cdots
M_t}^{L}(q^n)=g_{M_1\otimes_{k}K,\cdots,M_t\otimes_{k}K}^{L\otimes_{k}K}$
where $K$ is a finite extension of $k$ with $[K:k]=n$ and
$M\otimes_{k} K$ is naturally a $\Lambda\otimes_{k}K$-module.
Similarly, we have the polynomial $a_M$ such that
$a_M(q^n)=|\aut_{\Lambda\otimes_{k}K}(M\otimes_{k}K)|.$

Now we try to calculate the inner product on elements in the
PBW-basis $\{E^c|c\in\cm\}.$ First we consider $E_{m\delta},
m\geq1,$ which are defined in Section 7.1. We have obtained the
following fundamental relations:
$$\te_{n\dz}=E_{(n-1,n)}\ast E_1-v^{-2}E_1\ast E_{(n-1,n)}.$$
$$E_{0\dz}=1,\ \ E_{n\dz}=\frac{1}{[n]}\sum_{s=1}^n
v^{s-n}\te_{s\dz}\ast E_{(n-s)\dz}.$$ According to the
calculations in \cite{BCP}, $E_{n\delta}$ correspond to the
complete symmetric functions $h_{(n)}$ in \cite{M} and
$E_{\omega_{\bc}\delta}$ correspond to the complete symmetric
functions $h_{\omega_{\bc}}$ in \cite{M} (see Section 1 and
Section 3 in \cite{BCP}). Then
$$(E_{n\delta},E_{n\delta})\equiv 1\ \ \ (\mod
v^{-1}\bbq[[v^{-1}]]\cap\bbq(v))$$
$$
(E_{n\delta},E_{\omega_{\bc}\delta})\in
\bbn^{*}+v^{-1}\bbq[[v^{-1}]]\cap\bbq(v)
$$
$$
(E_{\omega_{\bc}\delta},E_{\omega_{\bc}\delta})\in
\bbn^{*}+v^{-1}\bbq[[v^{-1}]]\cap\bbq(v)
$$
for any $n\geq 0$ and any partition $\omega_{\bc}$ of $n.$

Let $F: \mh^*(K)\rightarrow \mh^
   *(\Lambda)$ be the embedding and $\fk{C}(P,L)$ be the full subcategory of $\mod\llz$ with two relative simple objects
   $S_1, S_2$ as in Section 7.1. We denote by $\fk{C}_0$ (resp. $\fk{C}_1$) the full subcategory of $\fk{C}(P,L)$ consisting of the $\Lambda$-modules which belong to
   homogeneous (resp. non-homogeneous) tubes of $\mod\Lambda.$  It is easy to see that this embedding preserves the inner product. According to the Auslander-Reiten quiver of regular $\llz$-modules, we
   have non-homogeneous tubes
   $\ct_i, i=1,\cdots,s,$ $s\leq3$ and others are homogeneous tubes.
   Now we have the decomposition
   $$ E_{n\delta}=E_{n\delta,1}+E_{n\delta,2}+E_{n\delta,3} $$
   where
   $$\hspace{-1cm}E_{n\delta,1}=v^{-n \dim S_1-n \dim S_2}
   \sum_{[M],M\in\fk{C}_1,\udim M=n\delta}u_{[M]},$$

   $$ E_{n\delta,2}=v^{-n \dim S_1-n \dim S_2}
   \sum_{\substack{[M],\udim M=n\delta\\
   M=M_1\oplus M_2, 0\neq M_1\in \fk{C}_1, 0\neq M_2\in \fk{C}_0}}u_{[M]},$$

   $$\hspace{-2cm}E_{n\delta,3}=
   v^{-n \dim S_1-n \dim S_2}\sum_{[M], M\in\fk{C}_0, \udim M=n\delta}u_{[M]}$$
Note that $\dim S_i=\dim_k S_i, i=1,2,$  but the values are
   independent of the choice of finite field $k.$ It is easy to see that
   $(E_{n\delta,i},E_{n\delta,j})=0$ for all $i\neq j.$

\bigskip

 \nd {\bf Lemma 9.2} {\sl We have  the relations:
 $$(E_{n\delta,1},E_{n\delta,1})\equiv 0\ \ \  (\mod
 v^{-1}\bbq[[v^{-1}]]\cap\bbq(v)),$$ $$
 (E_{n\delta,2},E_{n\delta,2})\equiv 0\ \ \ (\mod  v^{-1}\bbq[[v^{-1}]]\cap\bbq(v))$$ and
 $$ (E_{n\delta,3},E_{n\delta,3})\equiv 1\ \ \  (\mod v^{-1}\bbq[[v^{-1}]]\cap\bbq(v))$$}

\bigskip
\nd {\bf Proof}\ \  Since

$$E_{n\delta,1}=v^{-n \dim S_1-n \dim S_2}
   \sum_{[M],M\in\fk{C}_1,\udim M=n\delta}u_{[M]}$$
   $$=\sum_{[M],M\in\fk{C}_1,\udim M=n\delta}
   v^{-\dim \ed(M)}\lr{M}$$
 we have $$\hspace{-0.2cm}(E_{n\delta,1},E_{n\delta,1})=\sum_{[M],M\in\fk{C}_1,\udim M=n\delta}
   v^{-2\dim \ed(M)}(\lr{M},\lr{M}).$$ Note that
   $(\lr{M},\lr{M})=\frac{|\ed(M)| }{a_M}\in \bbq[[v^{-1}]]\cap\bbq(v),$
   $|\ed(M)|=v^{2\dim\ed(M)}$ and $a_M$ is a polynomial on $v$
   with the leading term $v^{2\dim\ed(M)},$ then
$$(\lr{M},\lr{M})\in 1+v^{-1}\bbq[[v^{-1}]]\cap\bbq(v))$$ and $$(E_{n\delta,1},E_{n\delta,1})\in
v^{-1}\bbq[[v^{-1}]]\cap\bbq(v).$$

Obviously, the result is true for $n=1.$ We assume now that it is
also true for all $m$ with $m<n$. Since
\begin{eqnarray}
  E_{n\delta,2} &=& v^{-n \dim S_1-n \dim S_2}\sum_
   {\substack{[M],\udim M=n\delta\\
   M=M_1\oplus M_2, 0\neq M_1\in \fk{C}_1, 0\neq M_2\in \fk{C}_0}}u_{[M]} \nonumber\\
   &=&\sum_{[M_1], M_1\in \fk{C}_1}v^{-\dim \ed(M_1)}\lr{M_1}*E_{n\delta-\udim M_1,3} \nonumber
\end{eqnarray}

   We have
   $$(E_{n\delta,2},E_{n\delta,2})=\sum_{[M_1], M_1\in\fk{C}_1}v^{-2\dim \ed(M_1)}(\lr{M_1},\lr{M_1})
   (E_{n\delta-\udim M_1,3},E_{n\delta-\udim M_1,3})$$
Since $\dim \ed(M_1)\geq 1,$ by the inductive assumption, we have
   $$(E_{n\delta,2},E_{n\delta,2})\equiv 0 \ \ \ (v^{-1}\bbq[[v^{-1}]]\cap\bbq(v)),$$
for all $n>0.$ Since $$(E_{n\delta},E_{n\delta})\equiv 1 \ \ \
(v^{-1}\bbq[[v^{-1}]]\cap\bbq(v))$$
 and   $$\quad\quad\quad(E_{n\delta},E_{n\delta})=(E_{n\delta,1}+E_{n\delta,2}+E_{n\delta,3},E_{n\delta,1}+
   E_{n\delta,2}+E_{n\delta,3})
   ,$$ then

   $$(E_{n\delta,3},E_{n\delta,3})\equiv 1 \ \ \ (\mod
   v^{-1}\bbq[[v^{-1}]]\cap\bbq(v)).$$

   The result is true for all $n.$ \hfill$\Box$

\bigskip

\nd {\bf  9.3}\ \ In the following, we will define a  decomposition
of regular part of $\cc^*(\Lambda)$ with respect to the inner
product $(-,-).$

\bigskip
\nd {\bf Lemma 9.3.1} {\sl Let $M$ be a regular module with $\udim
M=n\delta,$ and $M=\bigoplus_{i=1}^s M_i, M_i\in\ct_i$ for
$i=1,\cdots, s,$ then $\udim M_i=n_i\delta$ and $\sum_{i=1}^s
n_i=n.$}\

\nd {\bf Proof} Let $M=M_1\oplus M'_1, M_1\in\ct_1$ and $M'_1$ has
no nonzero direct summand in $\ct_1.$ Otherwise, we may assume
$\udim M_1=m_1\delta+\beta_1$ with $0<\beta_1<\delta$ and $\udim
M'_1=m'_1\delta+\beta'_1$ with $0<\beta'_1<\delta,$ then
$\beta_1+\beta'_1=\delta.$ Since
$0=(m_1\delta+\beta_1,m'_1\delta+\beta'_1)=(\beta_1,\beta'_1)=(\delta-\beta'_1,\beta'_1)
=-(\beta'_1,\beta'_1),$ we get $\beta'_1=k\delta,k\in\bbn$. This is
a contradiction. \hfill$\Box$
\bigskip

 In Section 7, we have constructed the $\bbq (v)$-basis $\{E^{\bc}\mid \bc\in \mathcal{M}\}$ of $\cc^*(\Lambda).$
 Let $\mr(\cc^*(\Lambda))$
 be the $\bbq (v)$-subspace of $\cc^*(\Lambda)$ with the basis $\{E_{\pi_{1 \bc}}\ast  E_{\pi_{2
\bc}}\ast\cdots\ast  E_{\pi_{s \bc}}\ast
E_{w_{\bc}\dz}|\pi_{\bc}=(\pi_{1\bc}, \dots, \pi_{s\bc})\in
\Pi_{1}^a\times \cdots \times\Pi_{s}^{a}, \text {and }
w_{\bc}=(w_1\geq w_2\geq\cdots\geq w_t)\text{ is a partition}\}.$
Obviously, it is a subalgebra of $\cc^*(\Lambda).$ Naturally, we
take $E_{w_\bc\dz}=1$ if $w_{\bc}=0.$

Let $\mr^a(\cc^*(\Lambda))$ be the subalgebra of
$\mr(\cc^*(\Lambda))$ with the basis $\{E_{\pi_{1 \bc}}\ast
E_{\pi_{2 \bc}}\ast\cdots\ast E_{\pi_{s \bc}}|\pi_{\bc}=(\pi_{1\bc},
\dots, \pi_{s\bc})\in \Pi_{1}^a\times \cdots \times\Pi_{s}^{a}\}.$
For $\alpha,\beta\in\bbn[I],$  we denote $\alpha\leqslant\beta$ if
$\beta-\alpha\in\bbn[I].$ It follows that
$\mr^a(\cc^*(\Lambda))_\beta=\mr(\cc^*(\Lambda))_\beta$ provided
$\beta<\delta.$ We now define
$$\mf_\delta=\{x|(x,\mr^a(\cc^*(\Lambda))_\delta)=0\}.$$

By Proposition 9.1,  $(-,-)$ is   nondegenerate on the
$\mr(\cc^*(\Lambda)).$  According to Lemma 9.3.1, we get

$$\mr(\cc^*(\Lambda))_\delta=\mr^a(\cc^*(\Lambda))_\delta\oplus\mf_\delta
\text{\quad and \quad} \dim(\mf_\delta)=1.$$

By the method of Schmidt orthogonalization, we may set
$$E'_{\delta}= E_{\delta}-\sum_{
M(\pi_{i\bc}),\udim M(\pi_{i \bc})=\delta,1\leqslant i\leqslant
s}a_{\pi_{i \bc}}E_{\pi_{i \bc}}.
$$
satisfying $\mf_\delta=\bbq(v)E'_{\delta}.$ Now let
$\mr(\cc^*(\Lambda))(1)$ be the subalgebra of
 $\mr(\cc^*(\Lambda))$ generated by $\mr^a(\cc^*(\Lambda))$ and $\mf_\delta.$
We have $\mr(\cc^*(\Lambda))(1)_\beta=\mr(\cc^*(\Lambda))_\beta$ if
$\beta<2\delta.$ Define
$$\mf_{2\delta}=\{x|(x,\mr(\cc^*(\Lambda))(1)_{2\delta})=0\}.$$
Then $\dim \mf_{2\delta}=1$ and
$\mr(\cc^*(\Lambda))_{2\delta}=\mr(\cc^*(\Lambda))(1)_{2\delta}\oplus\mf_{2\delta}.$

In general, define
$$ \mf_{n\delta}=\{x\in\mr(\cc^*(\Lambda))_{n\delta}|(x,\mr(\cc^*(\Lambda))(n-1)_{n\delta})=0\}.$$
Let $\mr(\cc^*(\Lambda))(n)$ be the subalgebra of
$\mr(\cc^*(\Lambda))$ generated by $\mr(\cc^*(\Lambda))(n-1)$ and $
\mf_{n\delta}.$  We have
$\mr(\cc^*(\Lambda))_{n\delta}=\mr(\cc^*(\Lambda))(n-1)_{n\delta}\oplus\mf_{n\delta},
\dim\mf_{n\delta}=1.$ Also, we can choose $E'_{n\delta}$ such that $
E_{n\delta}-E'_{n\delta}\in \mr(\cc^*(\Lambda))(n-1)_{n\delta}$ and
$\mf_{n\delta}=\bbq(v) E'_{n\delta}$ for all $n>0.$

We shall need the following facts:

\bigskip
 \nd {\bf Lemma 9.3.2}  {\sl Let $M,N,M_i$ be regular $\Lambda$-modules with
$\udim M,\udim N,\udim M_i$ $\in\bbn\delta,$ Then the degree of the
Hall polynomial $\varphi_{MN}^{L}$ is no more than $\dim
\ed(L)-(\dim \ed(M)+\dim \ed(N)).$}

\nd {\bf Proof} We know the formula
$g^{L}_{MN}=\frac{|\ext^1(M,N)_L||\aut L|}{|\aut M||\aut
N||Hom(M,N)|}.$ We have the Hall polynomial $\varphi_{MN}^{L}$ for
$g_{MN}^{L},$  the polynomials $a_M,a_N,a_L$ and the polynomial for
$|Hom(M,N)|.$ Therefore, we have the rational function, denoted by
$f,$ such that
$f(|K|)=|\ext^1(M\otimes_{k}K,N\otimes_{k}K)_{L\otimes_{k}K}|$ for
any finite extension $K$ of $k.$  Since $f(|K|)$ is an integer for
any finite extension $K$ of $k,$ $f$ is a polynomial with
coefficients in $\mathbb{Q}.$ Since $\langle \udim M,\udim N
\rangle=0,$ $\dim \ext^1(M,N)=\dim \hom(M,N).$ The degree of the
polynomial $f$ is no more than $\dim \ext^1(M,N).$ So,
$$\deg(f)\leq
\deg(a_L)-(\deg(a_M)+\deg(a_N))$$  It is also
known $\deg(a_X)=\dim_{k} \ed(X)$ for any $\Lambda$-module
$X.$ The proof is finished. \hfill$\Box$

\bigskip
Let $w_{\bc}=(w_1,\cdots,w_t)$ be a partition of $n,$ then
$$E_{w_{\bc}\delta}=E_{w_1\delta}*\cdots*E_{w_t\delta}=
(E_{w_1\dz,1}+E_{w_1\dz,2}+E_{w_1\dz,3})*\cdots*(E_{w_t\dz,1}+E_{w_t\dz,2}+E_{w_t\dz,3}).$$
We set $E_{w_{\bc}\dz,3}=E_{w_1\dz,3}*\cdots *E_{w_t\dz,3}$
\\
\\
\nd {\bf Lemma 9.3.3} {\sl Let $w_{\bc}$ be a partition of $n,$ then
$$(E_{n\delta},E_{w_{\bc}\dz} )\equiv(E_{n\delta,3},E_{w_{\bc}\dz,3})
\ \ \ (mod \ v^{-1}\bbq[[v^{-1}]]\cap\bbq(v))$$}\

\nd {\bf Proof}
 When $t=1,$ the result is true by Lemma 9.2. Suppose
$$(E_{n\delta},E_{w_{\bc}\dz}
)\equiv(E_{n\delta,3},E_{w_{\bc}\dz,3})\ \ \ (mod\
v^{-1}\bbq[[v^{-1}]]\cap\bbq(v))$$ for all $|(w_{\bc})|=t<s.$ Now
let $|(w_{\bc})|=s.$ Since
$E_{k\delta}=E_{k\delta,1}+E_{k\delta,2}+E_{k\delta,3}$ for any
$k\in \mathbb{N},$ we get by Section 9.2:
$$E_{m\delta,2}=v^{-m \dim S_1-m \dim S_2}\sum_
   {\substack{[M], M=M_1\oplus M_2\in\fk{C}(P,L),\udim M=m\delta\\
   0\neq M_1\in \fk{C}_1, 0\neq M_2\in \fk{C}_0} }u_{[M]}$$
   $$\hspace{-3.5cm}=\sum_{[M_1], M_1\in\fk{C}_1}v^{-\dim \ed(M_1)}\lr{M_1}*E_{m\delta-\udim M_1,3}$$
for all $m>0.$  Then

$$E_{w_{\bc}\dz}=
(E_{w_1\dz,1}+E_{w_1\dz,2}+E_{w_1\dz,3})*\cdots*(E_{w_t\dz,1}+E_{w_t\dz,2}+E_{w_t\dz,3})$$
$$
=E_{w_1\dz,1}*\cdots *E_{w_t\dz,1}+E_{w_1\dz,3}*\cdots
*E_{w_t\dz,3}+\mbox{ ``the rest part''}
$$
Here
\begin{eqnarray*}
&&E_{w_1\dz,1}*\cdots *E_{w_t\dz,1}\\
 &=&\left(v^{-w_1\dim S_1-w_1\dim S_2}
\sum_{\substack{[M_1], M_1\in \fk{C}_1, \\
\udim M_1=w_1\dz}}u_{[M_1]}\right)*\cdots*\left(v^{-w_t\dim S_1-w_t\dim
S_2}\sum_{\substack{[M_t], M_t\in \fk{C}_1, \\
 \udim M_t=w_t\dz}}u_{[M_t]}\right)\\
&=&v^{-n\dim S_1-n\dim S_2}\sum_{\substack{[M_1],\cdots, [M_t],
M_1\in \fk{C}_1,\cdots, M_t\in \fk{C}_1, \\
\udim M_1=w_1\dz,\cdots,\udim M_t=w_t\dz}}
\sum_{[L],L\in \fk{C}_1}\varphi_{M_1\cdots
M_t}^{L}(v^2)u_{[L]}\\
&=&\sum_{[L],0\neq L\in\fk{C}_1}v^{-\dim \ed
L}\varphi_{M_1\cdots M_t}^{L}(v^2) \lr{L}
\end{eqnarray*}
Using the above expression of $E_{m\delta,2},$ we have
$$\mbox{``the rest part''}
=\sum_{\substack{[M_{i_1}],\cdots,[M_{i_i}];\\ \sum_{j=1}^t\udim M_{i_j}=l\dz, l<n} }\sum_{\substack{[M],w'_{\bc};  0\neq M\in\fk{C}_1, \\1\leq |(w'_{\bc})|< s}} v^{-\dim \ed M}\varphi_{M_{i_1}\cdots M_{i_t}}^{M}(v^2)\lr{M}*E_{w'_{\bc}\dz,3}\\
.$$

Applying the above expansion of $E_{w_{\bc}\dz}$ to $(E_{n\dz},
E_{w_{\bc}\dz}),$ and by Proposition 9.1.1 and Lemma 9.3.2, we have
$$(E_{n\delta},E_{w_{\bc}\dz}
)\equiv(E_{n\delta,3},E_{w_{\bc}\dz,3}) \ \ \ (\mod\
v^{-1}\bbq[[v^{-1}]]\cap\bbq(v)).$$ \hfill$\Box$

\bigskip
\nd {\bf Lemma 9.3.4} {\sl Assume $m+n=p+q, m\geq n\geq 0, p\geq
q\geq 0.$ Then
\begin{enumerate}
    \item If $m\neq p$ or $n\neq q,$ $(E'_{m\delta}\ast E'_{n\delta},E'_{p\delta}\ast
E'_{q\delta})=0.$
    \item $(E'_{p\delta}\ast E'_{q\delta},E'_{p\delta}\ast
E'_{q\delta})=(E'_{p\delta},E'_{p\delta})(E'_{q\delta},E'_{q\delta})$
    \item $(E'_{p\delta}\ast E'_{q\delta},E'_{q\delta}\ast
E'_{p\delta})=(E'_{p\delta},E'_{p\delta})(E'_{q\delta},E'_{q\delta})$
    \item
    $((E'_{m\delta})^n,(E'_{m\delta})^n)=n!(E'_{m\delta},E'_{m\delta})^n.$
\end{enumerate}}
\nd {\bf Proof}  We remark that if $0\rightarrow R\rightarrow
M\rightarrow I\rightarrow 0$ is exact for a regular module $R$ and a
preinjective module $I$, then $M$ contains a nonzero direct summand
being preinjective and contains no nonzero summand being
preprojective. Dually, if $0\rightarrow P\rightarrow N\rightarrow
R\rightarrow 0$ is exact for a preprojective module $P$ and a
regular module $R,$ then $N$ contains a nonzero direct summand being
preprojective and contains no nonzero summand being preinjective.

We may assume that $p\geq m,$ then $n\geq q.$  As we know
$\mr(\cc^*(\Lambda))_{k\delta}=\mr(\cc^*(\Lambda))(k-1)_{k\delta}\oplus
\bbq(v)E'_{k\delta},$ then
$$r(E'_{k\delta})=E'_{k\delta}\otimes 1+1\otimes
E'_{k\delta}+\sum_{i}a_{ki}\otimes
b_{ki}+\sum_{i}x_{ki}\ast\lr{I_{ki}}\otimes \lr{P_{ki}}\ast
y_{ki},$$ where
$a_{ki},b_{ki},x_{ki},y_{ki}\in\mr(\cc^*(\Lambda))(j)$ and $j<k,$
$I_{ki}$ is a nonzero preinjective module and $P_{ki}$ is a nonzero
preprojective module.

We have
$$(E'_{m\delta}\ast E'_{n\delta}, E'_{p\delta}\ast
E'_{q\delta})=(r(E'_{m\delta}\ast E'_{n\delta}),E'_{p\delta}\otimes
E'_{q\delta})=(r(E'_{m\delta})\ast
r(E'_{n\delta}),E'_{p\delta}\otimes E'_{q\delta})$$ By the above
formula of $r(E'_{k\delta})$ and the above remark, we have
$$\hspace{-11cm}(E'_{m\delta}\ast E'_{n\delta}, E'_{p\delta}\ast
E'_{q\delta})$$$$=((E'_{m\delta}\otimes 1+1\otimes
E'_{m\delta}+\sum_{i}a_{mi}\otimes b_{mi})\ast(E'_{n\delta}\otimes
1+1\otimes E'_{n\delta}+\sum_{i}a_{ni}\otimes
b_{ni}),E'_{p\delta}\otimes E'_{q\delta})$$ If $p>m,$ it is easy to
see that it vanishes.  If $p=m>q=n,$ it is easy to see that it
equals to $(E'_{p\delta},E'_{p\delta})(E'_{q\delta},E'_{q\delta}).$
If $p=q=m=n,$ it equals to $2(E'_{p\delta},E'_{p\delta})^2.$ In
general, we have
$((E'_{p\delta})^l,(E'_{p\delta})^l)=l!(E'_{p\delta},E'_{p\delta})^l$
for $l\geq 0.$

Similarly, we can prove $(E'_{p\delta}\ast
E'_{q\delta},E'_{q\delta}\ast
E'_{p\delta})=(E'_{p\delta},E'_{p\delta})(E'_{q\delta},E'_{q\delta})$
for $p>q.$ \hfill$\Box$

\bigskip
\nd {\bf Corollary 9.3.5} {\sl For $m_i,n_i\in \bbn $ $(i=1,\cdots,t)$
satisfying $m_1>\cdots>m_t$ and $l_i,k_i\in\bbn $ $(i=1,\cdots,j)$
satisfying $l_1>\cdots
>l_j,$ we have
$$((E'_{m_1\delta})^{n_1}\ast\cdots \ast
(E'_{m_t\delta})^{n_t},(E'_{l_1\delta})^{k_1}\ast\cdots \ast
(E'_{l_j\delta})^{k_j})\\
=((E'_{m_1\delta})^{n_1},(E'_{m_1\delta})^{n_1})\cdots((E'_{m_t\delta})^{n_t},(E'_{m_t\delta})^{n_t})$$
if $t=j, m_i=l_i$ and $n_i=k_i$ for all $i=1,\cdots t;$
$$((E'_{m_1\delta})^{n_1}\ast\cdots \ast
(E'_{m_t\delta})^{n_t},(E'_{l_1\delta})^{k_1}\ast\cdots \ast
(E'_{l_j\delta})^{k_j})=0$$ otherwise.  }\

\bigskip
For a partition $w=(w_1\geq w_2\geq\cdots\geq w_t),$ we define
$$
E'_{w\delta}=E'_{w_1\delta}\ast\cdots \ast E'_{w_t\delta}
$$
\nd {\bf Lemma 9.3.6} {\sl Let $\{E_{\pi}|\pi\in\Pi^a_i\}$ be the
$\cz$-basis of $\cc^*(\ct_i)_{\cz}$ defined in Section 7.3. We have
$$(E_\pi\ast
E'_{m\delta},E_{\pi'}\ast
E'_{n\delta})=\delta_{mn}(E_\pi,E_{\pi'})(E'_{m\delta},E'_{n\delta})$$
and
$$
(E_{\pi},E'_{w\delta})=0
$$}\

\nd {\bf Proof} We may assume that $m\leq n.$ It follows

$$r(E_{\pi})=E_{\pi}\otimes 1 +1\otimes
E_{\pi}+\sum_{\pi_1,\pi_2}c_{\pi_1,\pi_2}E_{\pi_1}\otimes
E_{\pi_2}+\sum_{\pi_1,\pi_2}d_{\pi_1,\pi_2}E_{\pi_1}\ast\lr{I_{\pi_1,\pi_2}}
\otimes \lr{P_{\pi_1,\pi_2}}\ast E_{\pi_2},$$ where
$c_{\pi_1,\pi_2},d_{\pi_1,\pi_2}\in\cz$ and
$M(\pi),M(\pi_1),M(\pi_2)$ belong to $\mathcal{T}_i,$
$I_{\pi_1,\pi_2}$ is a nonzero preinjective module and
$P_{\pi_1,\pi_2}$ is a nonzero preprojective module for all $i.$ As
we know

  $$r(E'_{m\delta})=E'_{m\delta}\otimes 1+1\otimes
E'_{m\delta}+\sum_{i}a_{mi}\otimes b_{mi}+\sum_{i}x_{mi}\ast
\lr{I_{mi}}\otimes \lr{P_{mi}}\ast y_{mi},$$ where
$a_{mi},b_{mi},x_{mi},y_{mi}\in\mr(\cc^*(\Lambda))(j),j<m,$ and
$I_{mi}$ is a nonzero preinjective module and $P_{mi}$ is a nonzero
preprojective module for all $i.$

The same calculation as that in the proof of Lemma 9.3.4 tells us
$$(E_\pi\ast E'_{m\delta},E_{\pi'}\ast
E'_{n\delta})=\delta_{mn}(E_\pi,E_{\pi'})(E'_{m\delta},E'_{n\delta}).$$
The proof of the second identity is the same. \hfill$\Box$

\bigskip
In Section 7.3, the set $\{\lr{M(\ba_{\bc})}\ast E_{\pi_{1 \bc}}\ast
E_{\pi_{2 \bc}}\ast\cdots\ast E_{\pi_{s \bc}}\lr{M(\bb_{\bc})}\}$ is
the basis of $\mathcal{C}^{*}(\Lambda).$ In the same way, we obtain
\bigskip

 \nd {\bf Lemma 9.3.7} {\sl The following equalities hold
 \begin{enumerate}
    \item $(\lr{M(\ba_{\bc})}\ast E_{\pi_{1 \bc}}\ast E_{\pi_{2
\bc}}\ast\cdots\ast E_{\pi_{s \bc}},E'_{w_{\bc'}\dz})=0$ for
$\ba_{\bc}\neq 0$ and  partition $w_{\bc'}\neq 0.$
    \item $(E_\pi,E_{\pi_{1\bc}}\ast\cdots\ast
E_{\pi_{s\bc}}\ast E'_{w_{\bc}\dz})=0$ for $w_{\bc}\neq 0.$
    \item $(E_{\pi_{1\bc}}\ast\cdots\ast
E_{\pi_{s\bc}},E_{\pi_{1\bc'}}\ast\cdots\ast E_{\pi_{s\bc'}}\ast
E'_{w_{\bc'}\dz})=0$ for $w_{\bc'}\neq 0.$
 \end{enumerate}}

\bigskip
Based on Lemma 9.3.6 and Lemma 9.3.7, we obtain

\bigskip
\nd {\bf Lemma 9.3.8} {\sl The following holds
$$(E_{\pi_{i\bc}}\ast
E'_{w_{\bc}\dz}, E_{\pi_{j \bc'}}\ast E'_{w_{\bc'}\dz})=(E_{\pi_{i
\bc}},E_{\pi_{j \bc'}})(E'_{w_{\bc}\dz},E'_{w_{\bc'}\dz}),1\leq
i,j\leq s.$$}\

\nd {\bf Proof} We know
$$r(E_{\pi_{i\bc}})=E_{\pi_{i\bc}}\otimes 1 +1\otimes
E_{\pi_{i\bc}}+\sum_{\pi_1,\pi_2}c_{\pi_1,\pi_2}E_{\pi_1}\otimes
E_{\pi_2}+\sum_{\pi_1,\pi_2}d_{\pi_1,\pi_2}E_{\pi_1}\ast\lr{I_{\pi_1,\pi_2}}
\otimes \lr{P_{\pi_1,\pi_2}}\ast E_{\pi_2},$$ where
$c_{\pi_1,\pi_2},d_{\pi_1,\pi_2}\in\cz$ and
$M(\pi),M(\pi_1),M(\pi_2)$ belong to $\mathcal{T}_i,$
$I_{\pi_1,\pi_2}$ is a nonzero preinjective module and
$P_{\pi_1,\pi_2}$ is a nonzero preprojective module for all $i.$ Let
$$r^0(E_{\pi_{i\bc}})=E_{\pi_{i\bc}}\otimes 1 +1\otimes
E_{\pi_{i\bc}}+\sum_{\pi_1,\pi_2}c_{\pi_1,\pi_2}E_{\pi_1}\otimes
E_{\pi_2}$$ and
$$
r^1(E_{\pi_{i\bc}})=E_{\pi_{i\bc}}\otimes 1 +1\otimes E_{\pi_{i\bc}}
$$
Also for $w_{\bc}=(w_1,\cdots,w_t)$
\begin{eqnarray*}
&&r(E'_{w_{\bc}\delta})=r(E'_{w_1\delta})\ast\cdots \ast
r(E'_{w_t\delta}) \\
&&\quad =\left(E'_{w_1\delta}\otimes 1+1\otimes
E'_{w_1\delta}+\sum_{i}a_{w_1i}\otimes b_{w_1i}+\sum_{i}x_{w_1i}\ast
\lr{I_{w_1i}}\otimes \lr{P_{w_1i}}\ast
y_{w_1i}\right)\\
&& \quad \quad \ast\cdots\ast \left(E'_{w_t\delta}\otimes 1+1\otimes
E'_{w_t\delta}+\sum_{i}a_{w_ti}\otimes b_{w_ti}+\sum_{i}x_{w_ti}\ast
\lr{I_{w_ti}}\otimes \lr{P_{w_ti}}\ast y_{w_ti}\right)
\end{eqnarray*}
Let
$$
r^0(E'_{w_{\bc}\delta})=(E'_{w_1\delta}\otimes 1+1\otimes
E'_{w_1\delta}+\sum_{i}a_{w_1i}\otimes
b_{w_1i})\ast\cdots\ast(E'_{w_t\delta}\otimes 1+1\otimes
E'_{w_t\delta}+\sum_{i}a_{w_ti}\otimes b_{w_ti})
$$
and
$$
r^1(E'_{w_{\bc}\delta})=(E'_{w_1\delta}\otimes 1+1\otimes
E'_{w_1\delta})\ast\cdots\ast(E'_{w_t\delta}\otimes 1+1\otimes
E'_{w_t\delta})
$$
 It is clear that
$$
(E_{\pi_{i\bc}}\ast E'_{w_{\bc}\dz}, E_{\pi_{j \bc'}}\ast
E'_{w_{\bc'}\dz})=(r(E_{\pi_{i\bc}})\ast
r(E'_{w_{\bc}\dz}),E_{\pi_{j \bc'}}\otimes
E'_{w_{\bc'}\dz})=(r^0(E_{\pi_{i\bc}})\ast
r^0(E'_{w_{\bc}\dz}),E_{\pi_{j \bc'}}\otimes E'_{w_{\bc'}\dz})
$$
Based on Lemma 9.3.6 and Lemma 9.3.7, we use an induction on
$|w_{\bc\delta}|$ to obtain that
$$
(r^0(E_{\pi_{i\bc}})\ast r^0(E'_{w_{\bc}\dz}),E_{\pi_{j
\bc'}}\otimes E'_{w_{\bc'}\dz})=(r^1(E_{\pi_{i\bc}})\ast
r^1(E'_{w_{\bc}\dz}),E_{\pi_{j \bc'}}\otimes E'_{w_{\bc'}\dz})
$$
$$
=((E_{\pi_{i\bc}}\otimes 1 +1\otimes
E_{\pi_{i\bc}})*(E_{w_{\bc\delta}}\otimes 1+1\otimes
E_{w_{\bc\delta}}),E_{\pi_{j \bc'}}\otimes E'_{w_{\bc'}\dz})
=(E_{\pi_{i \bc}},E_{\pi_{j
\bc'}})(E'_{w_{\bc}\dz},E'_{w_{\bc'}\dz})$$ The proof is finished.
\hfill$\Box$

Finally, we have

\bigskip
\nd {\bf Theorem 9.3.9} {\sl With the same notations as above, the
following holds
$$(E_{\pi_{1 \bc}}\ast\cdots\ast E_{\pi_{s \bc}}\ast
E'_{w_{\bc}\dz}, E_{\pi_{1 \bc'}}\ast\cdots\ast
E_{\pi_{s \bc'}}\ast E'_{w_{\bc'}\dz})\\
=(E_{\pi_{1 \bc}},E_{\pi_{1 \bc'}})(E_{\pi_{s \bc}},E_{\pi_{s
\bc'}})\cdots(E'_{w_{\bc}\dz},E'_{w_{\bc'}\dz}).$$}\

\bigskip

 \nd {\bf 9.4.}
 In this subsection, we come to construct the canonical basis.
 Let $\ci\ct$ be the isomorphism classes of indecomposable objects in the non-homogeneous tubes $\ct_1, \dots, \ct_s$
 and $add(\ci\ct)$ be the objects that are direct sums of objects in $ \ci\ct$.

Theorem 9.3.9 and Corollary 9.3.5 imply that $(E'_{n\delta}\ast
E'_{m\delta}-E'_{m\delta}\ast E'_{n\delta},x)=0$ for all
$x\in\mr(\cc^*(\Lambda)).$ Thus $E'_{n\delta}\ast
E'_{m\delta}=E'_{m\delta}\ast E'_{n\delta}$ by the nondegeneracy of
$(-,-)$ on the $\mr(\cc^*(\Lambda)).$


\bigskip
\nd {\bf Lemma 9.4.1} {\sl Assume $\sum_{i=1}^s\dim M(\pi_{i
\bc})+|w_{\bc}|\delta=n\delta.$ Then
\begin{enumerate}
    \item $(E_{n\delta},E_{\pi_{1 \bc}}\ast\cdots\ast
E_{\pi_{s \bc}}\ast E_{w_{\bc}\dz})\in \bbq[[v^{-1}]]\cap\bbq(v).$
    \item moreover, if $|w_{\bc}|\delta<n\delta,$ then $(E_{n\delta},E_{\pi_{1
\bc}}\ast\cdots\ast E_{\pi_{s \bc}}\ast E_{w_{\bc}\dz})\in
v^{-1}\bbq[[v^{-1}]]\cap\bbq(v).$
    \item $(E_{\pi_{1
\bc}}\ast\cdots\ast E_{\pi_{s \bc}}\ast E_{w_{\bc}\dz},E_{\pi_{1
\bc}}\ast\cdots\ast E_{\pi_{s \bc}}\ast E_{w_{\bc}\dz})\in
v^h(\bbn+v^{-1}\bbq[[v^{-1}]]\cap\bbq(v))$ for some $ h\geq 0.$
\end{enumerate}}

\nd {\bf Proof} By the proof of Lemma 9.3.3, we have
$$E_{w_{\bc}\dz}
=E_{w_1\dz,1}*\cdots *E_{w_t\dz,1}+E_{w_1\dz,3}*\cdots
*E_{w_t\dz,3}+\mbox{ ``the rest part''},
$$
where
$$E_{w_1\dz,1}*\cdots *E_{w_t\dz,1}
=\sum_{[L],0\neq L\in\fk{C}_1}v^{-\dim \ed L}\varphi_{M_1\cdots
M_t}^{L}(v^2) \lr{L}
$$
and
$$\mbox{``the rest part''}=\sum_{\substack{[M_{i_1}],\cdots,[M_{i_i}];\\
 \sum_{j=1}^t\udim M_{i_j}=l\dz,\; l<n} }\sum_{\substack{[M],w'_{\bc};  0\neq M\in\fk{C}_1,\\ 1\leq |(w'_{\bc})|< |(w_{\bc})|}}
  v^{-\dim \ed M}\varphi_{M_{i_1}\cdots M_{i_t}}^{M}(v^2)\lr{M}*E_{w'_{\bc}\dz,3}\\
,$$

 If $M(\pi_{i\bc})=0$ for all $1\leq i\leq s,$ then both (1) and (3) are true by
 the property of the complete symmetric function with respect to the
 pair $(-,-),$ see subsection 9.2.

 Suppose $M(\pi_{i \bc})\neq0$ for some $ i.$ By \cite{DDX},
 $$
E_{\pi_{i\bc}}=\lr{M(\pi_{i\bc})}+\sum_{\lz\in\Pi_{i}\setminus\Pi^a_{i},
\lz\prec \pi_{i\bc}}\eta^{\pi_{i\bc}}_{\lz}\lr{M(\lz)}
 $$
where $\eta^{\pi_{i\bc}}_{\lz}\in v^{-1}\bbz[v^{-1}].$
$$
E_{\pi_{1 \bc}}\ast\cdots\ast E_{\pi_{s
\bc}}=\lr{\bigoplus_{i=1}^{s} M(\pi_{i\bc})}+\sum_{N\in
add\{\mathcal{T}_i\mid i=1,\cdots,s\}}\eta_N \lr{N}
$$
where $\eta_N\in v^{-1}\bbz[v^{-1}].$
$$
E_{\pi_{1 \bc}}\ast\cdots\ast E_{\pi_{s \bc}}\ast
E_{w_1\dz,1}\ast\cdots \ast E_{w_t\dz,1}=\sum_{[L],0\neq
L\in\fk{C}_1}v^{-\dim \ed L}\varphi_{M_1\cdots
M_t}^{L}(v^2)\lr{\bigoplus_{i=1}^{s} M(\pi_{i\bc})}\ast \lr{L}
$$
$$+\sum_{N\in add(\mathcal{IT}),\; [L],0\neq
L\in\fk{C}_1}v^{-\dim \ed L}\varphi_{M_1\cdots M_t}^{L}(v^2)\eta_N
\lr{N}\ast \lr{L}
$$
Here,
$$
\lr{\bigoplus_{i=1}^{s} M(\pi_{i\bc})}\ast
\lr{L}=\sum_{[U]}v^{\dim\ed \bigoplus_{i=1}^{s}
M(\pi_{i\bc})+\dim\ed L-\dim\ed U}\varphi_{\oplus_{i=1}^{s}
M(\pi_{i\bc}),L}^U(v^2)\lr{U}
$$
and
$$
\lr{N}\ast\lr{L}=\sum_{[V]}v^{\dim\ed N+\dim\ed L-\dim\ed
V}\varphi_{NL}^{V}(v^2)\lr{V}
$$

 By Lemma 9.3.2, we know
 $$\deg_q\varphi_{\oplus_{i=1}^{s}
M(\pi_{i\bc}),L}^U\leq \dim\ed U-(\dim\ed L+\dim\ed \oplus_{i=1}^{s}
M(\pi_{i\bc}) )$$
$$
\deg_q\varphi_{ M_1\cdots M_t}^L\leq \dim\ed L-\sum_{i=1}^{t} \dim\ed
M_i
$$
and
$$
\deg_q\varphi_{NL}^{V}\leq \dim\ed V-(\dim\ed N+\dim\ed L)
$$
Hence,
\begin{eqnarray*}
&&E_{\pi_{1 \bc}}\ast\cdots\ast E_{\pi_{s \bc}}\ast
E_{w_1\dz,1}\ast\cdots \ast
E_{w_t\dz,1}\\
&& =\sum_{\substack{0\neq[U]}\in \mathcal{IT}}v^{\dim\ed
U-\dim\ed\oplus_{i=1}^sM(\pi_{i\bc})-2\sum_{i}\dim\ed
M_i}f_U(v^{-1})\lr{U}\\
&&\quad +\sum_{\substack{0\neq[V]}\in \mathcal{IT}}v^{\dim\ed V-\dim\ed N-2\sum_{i}\dim\ed M_i}f_V(v^{-1})\lr{V}
\end{eqnarray*}
where $f_U(v^{-1}),f_V(v^{-1})\in \bbq[v^{-1}].$

In general,
$$\hspace{-9cm}E_{\pi_{1 \bc}}\ast\cdots\ast E_{\pi_{s \bc}}\ast
E_{w_{\bc}\dz}$$
$$=\sum_{\substack{0\neq[L]}\in \mathcal{IT}}f_L\lr{L}\\
+\sum_{\substack{0\neq[M]\in \mathcal{IT},\; 1\leq
|(w'_{\bc})|< |w_{\bc}|}} f_M\lr{M}*E_{w'_{\bc}\dz,3}\\
$$$$+\lr{M(\pi_{1 \bc})\oplus\cdots\oplus M(\pi_{s
\bc})}*E_{w_{\bc}\dz,3}+\sum_{N\neq 0 \in \mathcal{IT}}f_N\lr{N}*E_{w_{\bc}\dz,3},$$ where $v^{-\dim \ed(L)}f_L,
v^{-\dim \ed(M)}f_M, f_N\in v^{-1}\bbq[v^{-1}].$

Using expressions of $E_{n\dz,1},E_{n\dz,2},E_{n\dz,3}$, it is easy
to see that
$$(E_{n\delta},E_{\pi_{1 \bc}}\ast\cdots\ast E_{\pi_{s \bc}}\ast
E_{w_{\bc}\dz})\in v^{-1}\bbq[[v^{-1}]]\cap\bbq(v).$$ Then the
conclusion (2) is proved.

By Lemma 9.3.3 and the property of complete symmetric function, we
have
$$(\lr{M(\pi_{1 \bc})\oplus\cdots\oplus M(\pi_{s
\bc})}*E_{w_{\bc}\dz,3},\lr{M(\pi_{1 \bc})\oplus\cdots\oplus
M(\pi_{s \bc})}*E_{w_{\bc}\dz,3})\in
\bbn^{*}+v^{-1}\bbq[[v^{-1}]]\cap\bbq(v))$$

Thus $$(E_{\pi_{1 \bc}}\ast\cdots\ast E_{\pi_{s \bc}}\ast
E_{w_{\bc}\dz},E_{\pi_{1 \bc}}\ast\cdots\ast E_{\pi_{s \bc}}\ast
E_{w_{\bc}\dz})\in v^h(\bbn+v^{-1}\bbq[[v^{-1}]]\cap\bbq(v))$$ for
some $h\geq 0.$ Then the conclusion (3) is proved. \hfill$\Box$

 Let $\mathcal{A}=\bbq[v,v^{-1}].$ The lattice
$\mathcal{L}$ defined in [L5] is the $\bbq[v^{-1}]$-submodule
of $\cc^{*}(\Lambda)_{\ca}$ which is characterized by
$$
\mathcal{L}=\{x\in \cc^{*}(\Lambda) \mid (x,x)\in
\bbq[[v^{-1}]]\cap\bbq(v)\}
$$

\bigskip
\nd {\bf Lemma 9.4.2} {\sl We have $E'_{n\delta} \in \mathcal{L}$
and $(E'_{n\delta},E'_{n\delta})
 \equiv 1/n\ \ \ (\mod v^{-1}\bbq[[v^{-1}]]\cap\bbq(v)).$}\

\nd {\bf Proof} We know $\{E_{\pi_{i\bc}}\mid \udim
M(\pi_{i\bc})=\delta, 1\leq i\leq s\}$ is the basis of
$\mathcal{R}^{a}(\mathcal{C}^{*}(\Lambda))_{\delta}.$ Via Schmidt
orthogonalization, we define $E'_{\pi_{i\bc}}$ to be the orthogonal
element corresponding to $E_{\pi_{i\bc}}$ for $i=1,\cdots,s.$

It is easy to see that $(E'_\pi,E'_\pi)\equiv  1$ and $(E_\delta,E'_\pi)\equiv  0\ \
\ (\mod \ v^{-1}\bbq[[v^{-1}]]\cap\bbq(v)).$ Thus

$$E'_\delta=E_\delta-\sum_\pi\frac{(E_\delta,E'_\pi)}{(E'_\pi,E'_\pi)}E'_\pi$$
Hence, $(E'_{\dz},E'_{\dz})\equiv 1\ \ \ (\mod\
v^{-1}\bbq[[v^{-1}]]\cap\bbq(v))).$

Now suppose $E'_{n\delta} \in \mathcal{L}$ and
$(E'_{n\delta},E'_{n\delta})
 \equiv 1/n\ \ \ (\mod v^{-1}\bbq[[v^{-1}]]\cap\bbq(v))$ for all $n$ with $n<k.$

By the definition of $E'_{k\delta}$, we also have $\{E_{\pi_{1
\bc}}\ast\cdots\ast E_{\pi_{s \bc}}\ast
E'_{w_{\bc}\dz},E_{n\delta}|\sum_{1\leq i\leq s}\udim M(\pi_{1
\bc})+|w_{\bc}|\dz=n\dz ,\text{ for some } \udim M(\pi_{1 \bc})\neq
0\}$ is a basis of $\mr(\cc^*(\Lambda))_{n\dz}.$

By Lemma 9.3.9 and the induction hypothesis, $\{E_{\pi_{1
\bc}}\ast\cdots\ast E_{\pi_{s \bc}}\ast
E'_{w_{\bc}\dz},E_{n\delta}\}\subset \mathcal{L}.$ Similar as in
Lemma 9.4.1, we may get $(E_{n\dz},E_{\pi_{1 \bc}}\ast\cdots\ast
E_{\pi_{s \bc}}\ast E'_{w_{\bc}\dz})\in
v^{-1}\bbq[[v^{-1}]]\cap\bbq(v)$ if there exists $i$ such that
$M(\pi_{i\bc})\neq 0.$

 Thus

$$E'_{n\dz}\equiv E_{n\dz}-\sum_{w_{\bc}\vdash n,w_{\bc}\neq
(n)}\frac{(E_{n\dz},E'_{w_{\bc}\dz})}{(E'_{w_{\bc}\dz},E'_{w_{\bc}\dz})}E'_{w_{\bc}\dz}
\ \ \ (\mod\ v^{-1}\mathcal{L}).$$

In addition, $(E'_{n\dz},E'_{n\dz})=(E_{n\dz},E'_{n\dz})$ by Schmidt
orthogonalization.

We  now claim that
$$(E_{n\dz},E'_{w_{\bc}\dz})=(E'_{w_1},E'_{w_1})^{k_1}(E'_{w_2},E'_{w_2})^{k_2}\cdots(E'_{w_t},E'_{w_t})^{k_t}$$
if $w_{\bc}=(w_1^{k_1},w_2^{k_2},\cdots,w_t^{k_t}),
w_1>w_2>\cdots>w_t.$

 As we known
 $$r(E_{n\dz})=\sum_{0\leq i\leq n}E_{i\dz}\otimes E_{(n-i)\dz}+\mbox{ ``the rest part''}.$$

 Let
 $r^0(E_{n\dz})=\sum_{0\leq i\leq n}E_{i\dz}\otimes E_{(n-i)\dz}$  and
$w'_{\bc}=(w_1^{k_1-1},w_2^{k_2},\cdots,w_t^{k_t}).$ Then
$$(E_{n\dz},E'_{w_{\bc}\dz})=(r(E_{n\dz}),E'_{w_1\dz}\otimes
E'_{w'_{\bc}\dz})=(r^0(E_{n\dz}),E'_{w_1\dz}\otimes
E'_{w'_{\bc}\dz}).$$

 Based on the definition of $E'_{w_1\dz}$ and Lemma 9.3.9, we have
 $$(E_{n\dz},E'_{w_{\bc}\dz})=
 (E_{w_1\dz},E'_{w_1\dz})(E_{(n-w_1)\dz},E'_{w'_{\bc}\dz}).$$

 By the induction hypothesis, we obtain
 $$(E_{n\dz},E'_{w_{\bc}\dz})
 =(E'_{w_1},E'_{w_1})^{k_1}(E'_{w_2},E'_{w_2})^{k_2}\cdots(E'_{w_t},E'_{w_t})^{k_t}.$$

 Since $n!=\sum_{(1^{r_1}2^{r_2}\cdots)\vdash n}\frac{n!}{\prod_{i\geq 1}r_i!i^{r_i}},$
 we have

 $(E'_{n\dz},E'_{n\dz})
 \equiv (E_{n\dz},E_{n\dz})-\sum_{w_{\bc}\vdash n,w_{\bc}\neq
(n)}\frac{(E_{n\dz},E'_{w_{\bc}\dz})^2}{(E'_{w_{\bc}\dz},E'_{w_{\bc}\dz})}\\
\equiv (E_{n\dz},E_{n\dz})-\sum_{(n)\neq(1^{r_1}2^{r_2}\cdots)\vdash
n}\frac{\prod_{i\geq 1}(E'_{i\dz},E'_{i\dz})^{r_i}}{\prod_{i\geq
1}r_i!}\\
\equiv 1-\sum_{(n)\neq(1^{r_1}2^{r_2}\cdots)\vdash
n}\frac{1}{\prod_{i\geq 1}r_i!i^{r_i}}\text{\quad ( by the induction
hypothesis )}\\
\equiv\frac{1}{n!}(n!-(\sum_{(1^{r_1}2^{r_2}\cdots)\vdash
n}\frac{n!}{\prod_{i\geq 1}r_i!i^{r_i}})+(n-1)!) \equiv 1/n\ \ \
(mod\ v^{-1}\bbq[[v^{-1}]]\cap\bbq(v)).$

So the proof is completed. \hfill$\Box$

Let $P_{n\dz}=n E'_{n\dz}.$  For a partition
$w_{\bc}=(1^{r_1}2^{r_2}\cdots t^{r_t}),$  let
$z_{w_{\bc}}=\prod_{i\geq 1}i^{r_i}r_i!$ and
$P_{w_{\bc}\dz}=P_{1\dz}^{\ast r_1}\ast\cdots\ast P_{t\dz}^{\ast
r_t}.$

\bigskip
\nd {\bf Corollary  9.4.3} Let $w_{\bc}=(1^{r_1}2^{r_2}\cdots),
w_{\bc'}=(1^{r'_1}2^{r'_2}\cdots)$ be partitions.  Then
\begin{enumerate}
    \item   $$(E_{\pi_{1
\bc}}\ast\cdots\ast E_{\pi_{s \bc}}\ast E'_{w_{\bc}\dz}, E_{\pi_{1
\bc'}}\ast\cdots\ast E_{\pi_{s \bc'}}\ast E'_{w_{\bc'}\dz})$$$$
\equiv\dz_{\pi_{1 \bc},\pi_{1 \bc'}}\cdots\dz_{\pi_{1 \bc},\pi_{1
\bc'}}\dz_{w_{\bc},w_{\bc'}}\prod_{i}r_i!
(E'_{i\dz},E'_{i\dz})^{r_i}\ \ \ (mod\
v^{-1}\bbq[[v^{-1}]]\cap\bbq(v))
$$
    \item  $$(P_{n\dz},P_{n\dz})\equiv n\ \ \ (mod\
v^{-1}\bbq[[v^{-1}]]\cap\bbq(v))$$ and
$$(P_{w_{\bc}\dz},P_{w_{\bc}\dz})\equiv\dz_{w_{\bc}w_{\bc'}}z_{w_{\bc}}\
\ \ (mod\ v^{-1}\bbq[[v^{-1}]]\cap\bbq(v)).$$

\end{enumerate}

\bigskip

By this property of $P_{w_{\bc}\dz}$, it is easy to see that
$P_{w_{\bc}\dz}$ corresponds to Newton symmetric functions.

Let $S_{{w_{\bc}\dz}}$ be the Schur functions corresponding to
$P_{w_{\bc}\dz},$ and $e^{\bc}=\lr{M(\ba_{\bc})}\ast E_{\pi_{1
\bc}}\ast E_{\pi_{2 \bc}}\ast\cdots\ast E_{\pi_{s \bc}}\ast
S_{{w_{\bc}\dz}}*\lr{M(\bb_{\bc})}$ for $\bc\in\cm.$

By Theorem 9.3.9, Lemma 9.4.2 Corollary 9.4.3 and the Nakayama
Lemma, we have the following corollary:

\bigskip
\nd {\bf Corollary 9.4.4} {\sl $\{e^{\bc}|\bc\in\cm\}$ is an almost
orthonormal basis of $\cl,$ that is, $(e^{\bc},
e^{\bc'})\in\delta_{\bc,\bc'}+v^{-1}\bbq[[v^{-1}]]\cap\bbq(v)$ for
$\bc$ and $\bc'\in\cm.$}\

We have defined the constructible set in Section 8.1
$$\co_{\bc}=\co_{M(\ba_{\bc})}\star\co_ {M_{\pi_{1 \bc}}}\star\co_
{M_{\pi_{2 \bc}}}\star\cdots\star\co_{M_{\pi_{s
\bc}}}\star\cn_{w_{\bc}\dz}\star\co_{M(\bb_{\bc})}$$ for any
$\bc\in\cm.$

 Now we define a new partial order
$\prec$ for those $e^{\bc}, \bc\in\cm$ with the same weight
(dimension vector) as follow:

 (1)$\hspace{10mm} e^{\bc}\prec e^{\bc'}\text{ if }\dim \co_{\bc}<\dim
\co_{\bc'}.$

 (2)$\hspace{10mm} e^{\bc}\prec e^{\bc'}\text{ if }\dim \co_{\bc}=\dim
\co_{\bc'}$ but $w_{\bc}>w_{\bc'}.$

 Base on the definition of
$E'_{n\dz},$ we  have

$$E'_{n\dz}= E_{n\dz}-\sum_{w_{\bc}\vdash n,w_{\bc}\neq
(n)}\frac{(E_{n\dz},E'_{w_{\bc}\dz})}{(E'_{w_{\bc}\dz},E'_{w_{\bc}\dz})}E'_{w_{\bc}\dz}+\sum_{\dim\co_{\bc'}<\dim\co_{n\dz}}a_{n\dz,\bc'}E_{\pi_{1
\bc'}}\ast E_{\pi_{2 \bc'}}\ast\cdots\ast E_{\pi_{s \bc'}}\ast
S_{{w_{\bc'}\dz}}$$ where $a_{n\dz,\bc'}\in\bbq(v),$ in fact, by
Corollary 9.4.4, we have
$a_{n\dz,\bc'}\in\bbq[[v^{-1}]]\cap\bbq(v).$  Thus

$$E_{n\dz}=E'_{n\dz}+\sum_{w_{\bc}\vdash n,w_{\bc}\neq
(n)}\frac{(E_{n\dz},E'_{w_{\bc}\dz})}{(E'_{w_{\bc}\dz},E'_{w_{\bc}\dz})}E'_{w_{\bc}\dz}
+\sum_{\dim\co_{\bc'}<\dim\co_{n\dz}}a_{n\dz,\bc'}E_{\pi_{1
\bc'}}\ast E_{\pi_{2 \bc'}}\ast\cdots\ast E_{\pi_{s \bc'}}\ast
S_{{w_{\bc'}\dz}},$$

$$E_{n\dz}=\frac{1}{n}P_{n\dz}+\sum_{w_{\bc}\vdash n,w_{\bc}\neq
(n)}\frac{1}{z_{w_{\bc}}}P_{w_{\bc}\dz}+\sum_{\dim\co_{\bc'}<\dim\co_{n\dz}}a_{n\dz,\bc'}E_{\pi_{1
\bc'}}\ast E_{\pi_{2 \bc'}}\ast\cdots\ast E_{\pi_{s \bc'}}\ast
S_{{w_{\bc'}\dz}}.$$

Let $H_{n\dz}$ be the $n$th complete symmetric function
corresponding to $P_{n\dz}.$  From \cite{M}, p25, we have
$$E_{n\dz}=H_{n\dz}+\sum_{\dim\co_{{\bc'}}<\dim\co_{n\dz}}a_{n\dz,\bc'}E_{\pi_{1
\bc'}}\ast E_{\pi_{2 \bc'}}\ast\cdots\ast E_{\pi_{s \bc'}}\ast
S_{{w_{\bc'}\dz}}.$$ Let $w_{\bc}$ be  a partition of $n,$ according
to Lemma 2.2 and the above formula, we have
$$E_{w_{\bc}\dz}=H_{w_{\bc}\dz}+\sum_{\dim\co_{{\bc'}}<\dim\co_{n\dz}}a_{n\dz,\bc'}E_{\pi_{1
\bc'}}\ast E_{\pi_{2 \bc'}}\ast\cdots\ast E_{\pi_{s \bc'}}\ast
S_{{w_{\bc'}\dz}}.$$ We have the monomial
$\mathfrak{m}_{w_{\bc}\dz}$ on the divided powers of $u_{S_i}, i\in
I,$ in Proposition 8.4, corresponding to $E_{w_{\bc}\dz}$, such that

$$\mathfrak{m}_{w_{\bc}\dz}=H_{w_{\bc}\dz}+\sum_{\dim\co_{{\bc'}}<\dim\co_{n\dz}}b_{n\dz,\bc'}E_{\pi_{1
\bc'}}\ast E_{\pi_{2 \bc'}}\ast\cdots\ast E_{\pi_{s \bc'}}\ast
S_{{w_{\bc'}\dz}}+
\sum_{\substack{\dim\co_{{\bc'}}<\dim\co_{n\dz};\\
{M(\ba_{\bc'})\neq 0}\ \mbox{or}\ M(\bb_{\bc'})\neq
0}}c_{n\dz,\bc'}e^{\bc'}$$
$$=S_{w_{\bc}\dz}+\sum_{w_{\bc''}>w_{\bc}}K_{w_{\bc''}w_{\bc}}S_{w_{\bc''}\dz}+
\sum_{\dim\co_{{\bc'}}<\dim\co_{n\dz}}b_{n\dz,\bc'}E_{\pi_{1
\bc'}}\ast E_{\pi_{2 \bc'}}\ast\cdots\ast E_{\pi_{s \bc'}}\ast
S_{{w_{\bc'}\dz}}$$$$+\sum_{\substack{\dim\co_{{\bc'}}<\dim\co_{n\dz};\\
{M(\ba_{\bc'})\neq 0}\ \mbox{or}\ M(\bb_{\bc'})\neq
0}}c_{n\dz,\bc'}e^{\bc'}$$ where $K_{\mu\lambda}$ are Kostka numbers
and $b_{n\dz,\bc'}, c_{n\dz,\bc'}\in\bbq(v).$  Furthermore, for
$\bc\in\cm$ and the monomials $\mathfrak{m}_{\bc}$ given in
Proposition 8.4, we have
$$\mathfrak{m}_{\bc}=e^{\bc}+\sum_{e^{\bc'}\prec
e^{\bc}}a_{\bc'\bc}e^{\bc'}$$ where $a_{\bc'\bc}\in\bbq(v).$ However
Proposition 8.4 and the above formulae tell us that the transition
matrix between $\{E^{\bc}|\bc\in\cm\}$ and $\{e^{\bc}|\bc\in\cm\}$
is triangular with diagonal entries equal to $1,$ and
$\{E^{\bc}|\bc\in\cm\}$ is an $\ca$-basis of
$\cc^{*}(\Lambda)_{\ca},$  $\{e^{\bc}|\bc\in\cm\}\subset\mathcal{L}$
and $\{\mathfrak{m}_{\bc}|\bc\in\cm\}\subset\cc^{*}(\Lambda)_{\ca}.$
Then the constants $a_{\bc'\bc}$ in the above formulae must lie in
$\mathcal{A}.$

By applying the same argument as in Section 8 to
$\{e^{\bc}|\bc\in\cm\},$ we obtain an $\mathcal{A}$-basis of
$\cc^{*}(\Lambda)_{\ca}$ which is denoted by
$\{\ce'^{\bc}|\bc\in\cm\}$ satisfying that
$$\ce'^{\bc}=\sum_{\bc'\in\cm}\zeta^{\bc}_{\bc'}e^{\bc'}\ \ \ \
\text{for any}\ \ \ \bc\in\cm$$ where $\zeta^{\bc}_{\bc}=1$ and $
\zeta^{\bc}_{\bc'}\in v^{-1}\bbq[v^{-1}]$ if $ e^{\bc'}\prec
e^{\bc}.$

 Finally,  we have
following theorem:

\mk\nd{\bf Theorem 9.4.5} {\it The set
$\{\ce'^{\bc}|\bc\in\cm\}\subset\cl$ provides an
$\mathcal{A}$-basis of $\cc^{*}(\Lambda)_{\ca}$ which is
characterized by the following three properties:

(a)\  $\ol{\ce'^{\bc}}=\ce'^{\bc}$ for all $\bc\in\cm.$

(b)\  $\pi(\ce'^{\bc})=\pi(e^{\bc}),$ where
$\pi:\mathcal{L}\ra\mathcal{L}/v^{-1}\mathcal{L}$ is the canonical
projection.

(c)\ $(\ce'^{\bc},\ce'^{\bc'})\equiv\dz_{\bc\bc'}\ \ \ (\mod
 v^{-1}\bbq[[v^{-1}]]\cap\bbq(v))$.}

According to Lusztig [L5], we obtain the signed canonical basis
$(\ce'^{\bc})$ of $\mathcal{L}.$  From the above formulae, we have
the relations

$$\mathfrak{m}_{\bc}=\ce'^{\bc}+\sum_{e^{\bc'}\prec
e^{\bc}}d_{\bc'\bc}\ce'^{\bc'}$$ where $d_{\bc'\bc}\in\mathcal{A}.$
By the total positivity of the canonical bases, we have

\mk\nd{\bf Theorem 9.4.7} {\sl The set $\{\ce'^{\bc}|\bc\in\cm\}$ is
the canonical basis of $\mathcal{L}$ in the sense of  Lusztig.}

This answers a question raised by  Nakajima in [N].

\end{document}